\documentclass[12pt]{article}
\usepackage[leqno]{amsmath}
\usepackage{amssymb,mathrsfs,bm}
\usepackage{hyperref}
\usepackage{mathtools}
\usepackage{fullpage}
\usepackage{epic}
\usepackage{eepic}
\usepackage[english]{babel}
\usepackage[all]{xy}
\usepackage{enumerate}
\usepackage[T1]{fontenc}
\usepackage[latin1]{inputenc}

\usepackage{soul}
\usepackage{mathabx}
\usepackage{dsfont}
\newtheorem{rem}{Remark}
\newtheorem{conj}{Conjecture}
\newtheorem{theointro}{Theorem}
\newtheorem{theo}{Theorem}[subsection]
\newtheorem{lem}{Lemma}[subsection]
\newtheorem{prop}{Proposition}[subsection]

\def\lv{\lvert}
\def\rv{\rvert}

\def\oF{\overline{F}}

\DeclareMathOperator{\Hom}{Hom}
\DeclareMathOperator{\Tr}{Trace}

\DeclareMathOperator{\No}{Norm}

\DeclareMathOperator{\Ker}{Ker}

\DeclareMathOperator{\Ad}{Ad}

\DeclareMathOperator{\vol}{vol}

\DeclareMathOperator{\Gal}{Gal}
\DeclareMathOperator{\op}{op}
\DeclareMathOperator{\cusp}{cusp}
\DeclareMathOperator{\sqr}{sqr}
\DeclareMathOperator{\scusp}{scusp}
\DeclareMathOperator{\Res}{Res}
\DeclareMathOperator{\St}{St}
\DeclareMathOperator{\ker1}{ker}
\DeclareMathOperator{\Cent}{Cent}
\DeclareMathOperator{\Irr}{Irr}
\DeclareMathOperator{\geom}{geom}
\DeclareMathOperator{\elli}{ell}
\DeclareMathOperator{\reg}{reg}

\numberwithin{equation}{subsection}
\newcounter{keepeqno}
\newenvironment{num}
 {\setcounter{keepeqno}{\value{equation}}%
  \begin{list}{(\theequation)}{\usecounter{equation}}%
  \setcounter{equation}{\value{keepeqno}}}
 {\end{list}}

\title{On distinguished square-integrable representations for Galois pairs and a conjecture of Prasad}
\author{Rapha\"el Beuzart-Plessis \protect\footnote{Universit\'e d'Aix-Marseille, I2M-CNRS(UMR 7373), Campus de Luminy, 13288 Marseille C\'edex 9, France rbeuzart@gmail.com}}

\begin{document}

\maketitle

\begin{abstract}
We prove an integral formula computing multiplicities of square-integrable representations relative to Galois pairs over $p$-adic fields and we apply this formula to verify two consequences of a conjecture of Dipendra Prasad. One concerns the exact computation of the multiplicity of the Steinberg representation and the other the invariance of multiplicities by transfer among inner forms.
\end{abstract}

\tableofcontents

\addcontentsline{toc}{section}{Introduction}
\section*{Introduction}

Let $F$ be a $p$-adic field (that is a finite extension of $\mathbf{Q}_p$ for a certain prime number $p$) and $H$ be a connected reductive group over $F$. Let $E/F$ be a quadratic extension and set $G:=R_{E/F}H_E$ where $R_{E/F}$ denotes Weil's restriction of scalars (so that $G(F)=H(E)$). To every complex smooth irreducible representation $\pi$ of $G(F)$ and every continuous character $\chi$ of $H(F)$ we associate a multiplicity $m(\pi,\chi)$ (which is always finite by \cite{Del} Theorem 4.5) defined by

$$\displaystyle m(\pi,\chi):=\dim \Hom_H(\pi,\chi)$$

\noindent where $\Hom_H(\pi,\chi)$ stands for the space of $(H(F),\chi)$-equivariant linear forms on (the space of) $\pi$. Recently, Dipendra Prasad \cite{Pras15} has proposed very general conjectures describing this multiplicity in terms of the Langlands parameterization of $\pi$, at least for representations belonging to the so-called `generic' $L$-packets. These predictions, which generalize earlier conjectures of Jacquet (\cite{JY1}, \cite{JY2}), are part of a larger stream that has come to be called the `local relative Langlands program' and whose main aim is roughly to describe the `spectrum' of general homogeneous {\it spherical} varieties $X=H\backslash G$ over local fields in terms of Langlands dual picture and correspondence. In the paper \cite{SV}, and under the assumption that $G$ is {\it split}, Sakellaridis and Venkatesh set up a very general framework to deal with these questions by introducing a certain complex reductive group $\check{G}_X$ associated to the variety $X$, which generalizes Langlands construction of a dual group, together with a morphism $\check{G}_X\to \check{G}$ (actually, in the most general case, this should also include an extra $SL_2$ factor) which, according to them, should govern a great part of the spectral decomposition of $L^2(H\backslash G)$ (see \cite{SV} Conjecture 16.2.2). In a similar way, in the case where $G=R_{E/F}H_E$ as above (note that such a group is never split) Prasad introduces a certain $L$-group ${}^L H^{\op}$ (further explanations on this notation are given below) and a morphism ${}^L H^{\op}\to {}^L G$ which should govern, on the dual side, the behavior of the multiplicities $m(\pi,\chi)$ for a very particular quadratic character $\chi$ (denoted by $\omega_{H,E}$ below) that has also been defined by Prasad. The main goal of this paper is to present some coarse results supporting Prasad's very precise conjectures in the particular case of stable (essentially) square-integrable representations. In the rest of this introduction we will recall the part of Prasad's conjecture that we are interested in as well as the two consequences of it that we have been able to verify. We will also say some words on the proofs which are based on a certain simple local trace formula adapted to the situation and which takes its roots in Arthur's local trace formula (\cite{Art}) as well as in Waldspurger's work on the Gross-Prasad conjecture for orthogonal groups (\cite{Wal2}, \cite{Wal3}).

\vspace{2mm}

Prasad associates a number of invariants to the situation at hand. First, he constructs a certain quadratic character $\omega_{H,E}:H(F)\to \{\pm 1 \}$ as well as a certain quasi-split group $H^{\op}$ over $F$ which is an $E/F$ form of the quasi-split inner form of $H$. We refer the reader to \cite{Pras15}\S 7-8 for precise constructions of those and content ourself to give three examples here:

\begin{itemize}
\item If $H=GL_n$, then $H^{\op}=U(n)_{qs}$ (quasi-split form) and $\omega_{H,E}=(\eta_{E/F} \circ \det)^{n+1}$ where $\eta_{E/F}$ is the quadratic character associated to $E/F$;

\item If $H=U(n)$ (a unitary group of rank $n$), then $H^{\op}=GL_n$ and $\omega_{H,E}=1$;

\item If $H=SO(2n+1)$ (any odd special orthogonal group), then $H^{\op}=SO(2n+1)_{qs}$ (the quasi-split inner form) and $\omega_{H,E}=\eta_{E/F}\circ N_{spin}$ where $N_{spin}:SO(2n+1)(F)\to F^\times/(F^{\times})^2$ denotes the spin norm.
\end{itemize}

\noindent To continue we need to restrict slightly the generality by only considering characters $\chi$ that are of `Galois type' i.e. which are in the image of a map constructed by Langlands

$$\displaystyle H^1(W_F,Z(\check{H}))\to \Hom_{cont}(H(F),\mathbf{C}^\times)$$

\noindent This map is always injective (because $F$ is $p$-adic) but not always surjective (although it is most of the time, e.g. if $H$ is quasi-split). We refer the reader to \cite{LL} for further discussion on these matters. The character $\omega_{H,E}$ is always of Galois type and, to every character $\chi$ of Galois type of $H(F)$, Prasad associates a certain `Langlands dual group' $\mathcal{H}^{\op}_\chi$ which sits in a short exact sequence

$$\displaystyle 1\to \check{H}^{\op}\to\mathcal{H}^{\op}_\chi\to W_F\to 1$$

\noindent together with a group embedding $\iota:\mathcal{H}^{\op}_\chi\hookrightarrow {}^L G$ (where ${}^L G$ denotes the $L$-group of $G$) compatible with the projections to $W_F$ and algebraic when restricted to $\check{H}^{\op}$. In the particular case where $\chi=\omega_{H,E}$, we have $\mathcal{H}^{\op}_\chi={}^L H^{\op}$ and $\iota$ is the homomorphism of quadratic base-change.

\begin{rem}
Although the short exact sequence above always splits, there does not necessarily exist a splitting preserving a pinning of $\check{H}^{\op}$ and hence $\mathcal{H}^{\op}_\chi$ is not always an $L$-group in the usual sense.
\end{rem}

Let $WD_F:=W_F\times SL_2(\mathbf{C})$ be the Weil-Deligne group of $F$. An `$L$-parameter' taking values in $\mathcal{H}_\chi^{\op}$ is defined as usual: a continuous Frobenius semi-simple morphism $WD_F\to \mathcal{H}_\chi^{\op}$ which commutes with the projections to $W_F$ and is algebraic when restricted to $SL_2(\mathbf{C})$. We are now ready to state (a slight generalization of) the stable version of Prasad's conjecture for square-integrable representations:

\begin{conj}\label{Prasad conj}

Let $\phi: WD_F\to {}^L G$ be a discrete L-parameter, $\Pi^G(\phi)\subseteq \Irr(G)$ the corresponding $L$-packet and $\displaystyle \Pi_\phi=\sum_{\pi\in \Pi^G(\phi)} d(\pi)\pi$ the stable representation associated to $\phi$. Then, we have

$$\displaystyle m(\Pi_\phi,\chi)=\lvert \ker1^1(F;H,G)\rvert^{-1} \sum_{\psi}\frac{\lvert Z(\phi)\rvert}{\lvert Z(\psi)\rvert}$$

\noindent where

\begin{itemize}
\item The sum is over the set of `$L$-parameters' $\psi:WD_F\to \mathcal{H}^{\op}_\chi$ (taken up to $\check{H}^{\op}$-conj) making the following diagram commute up to $\check{G}$-conj, i.e. there exists $g\in \check{G}$ such that $\iota\circ \psi=g\phi g^{-1}$,

\begin{center}
$\xymatrix{ & \mathcal{H}_\chi^{\op} \ar@{^{(}->}[d]^{\iota} \\ WD_F \ar@{->}[ru]^{\psi} \ar@{->}[r]^{\phi} & {}^LG}$
\end{center}
\item $\ker1^1(F;H,G):=\Ker\left(H^1(F,H)\to H^1(F,G)\right)$ (corresponds to certain twists of the parameter $\psi$ that become trivial in ${}^L G$);

\item $Z(\phi):=\Cent_{\check{G}}(\phi)/ Z(\check{G})^{W_F}$ and $Z(\psi):=\Cent_{\check{H}^{\op}}(\psi)/ Z(\check{H}^{\op})^{W_F}$
\end{itemize}
\end{conj}

\vspace{1mm}

\noindent As we said, this is only part of Prasad's general conjectures which aim to describe (almost) all the multiplicities $m(\pi,\chi)$ explicitly. This version of the conjecture (and far more) is known in few particular cases: for $H=GL(n)$ by Kable and Anandavardhanan-Rajan (\cite{Kable}, \cite{AR}), for $H=U(n)$ by Feigon-Lapid-Offen (\cite{FLO}) and for $H=GSp(4)$ by Hengfei Lu (\cite{Lu}). The following theorems are both formal consequences of Conjecture \ref{Prasad conj} and are the main results of this paper (see Theorem \ref{theo inner forms} and Theorem \ref{theo Steinberg}):

\begin{theointro}\label{inner forms}
Let $H$, $H'$ be inner forms over $F$, $G:=R_{E/F} H$, $G':=R_{E/F}H'$ and $\chi$, $\chi'$ characters of Galois type of $H(F)$ and $H'(F)$ corresponding to each other (i.e. coming from the same element in $H^1(W_F,Z(\check{H}))=H^1(W_F,Z(\check{H}'))$). Let $\Pi$, $\Pi'$ be (essentially) square-integrable representations of $G(F)$ and $G'(F)$ respectively which are stable (but not necessarily irreducible) and transfer of each other (i.e. $\Theta_\Pi(x)=\Theta_{\Pi'}(y)$ for all stably conjugate regular elements $x\in G_{\reg}(F)$ and $y\in G'_{\reg}(F)$ where $\Theta_\Pi$, $\Theta_{\Pi'}$ denote the Harish-Chandra characters of $\Pi$ and $\Pi'$ respectively). Then, we have

$$\displaystyle m(\Pi,\chi)=m(\Pi',\chi')$$
\end{theointro}

\begin{theointro}\label{Steinberg}
For $\pi=\St(G)$ the generalized Steinberg representation of $G(F)$ and $\chi$ a character of Galois type we have

$$\displaystyle m(\St(G),\chi)=\left\{
    \begin{array}{ll}
        1 & \mbox{if } \chi=\omega_{H,E} \\
        0 & \mbox{otherwise}
    \end{array}
\right.
$$
\end{theointro}

\vspace{1mm}

\noindent Theorem \ref{Steinberg} also confirms an older conjecture of Prasad (\cite{Pras01}, Conjecture 3) which was already known for split groups and tamely ramified extensions by work of Broussous-Court\`es and Court\`es (\cite{Brous}, \cite{Cour1}, \cite{Cour2}) and for inner forms of $GL_n$ by work of Matringe (\cite{Mat}). The proof of Broussous and Court\`es is mainly based on a careful study of the geometry of the building whereas Matringe's work uses some Mackey machinery. Our approach is completely orthogonal to theirs and is based on a certain integral formula computing the multiplicity $m(\pi,\chi)$ in terms of the Harish-Chandra character of $\pi$. This formula is reminiscent and inspired by a similar result of Waldspurger in the context of the so-called Gross-Prasad conjecture (\cite{Wal2}, \cite{Wal3}). It can also be seen as a `twisted' version (`twisted' with respect to the non-split extension $E/F$) of the orthogonality relations between characters of discrete series due Harish-Chandra (\cite{Clo}, Theorem 3). It can be stated as follows (see Theorem \ref{theo formula multiplicity}):

\begin{theointro}\label{formula}
Let $\pi$ be a square-integrable representation of $G(F)$ and $\chi$ be a continuous character of $H(F)$. Assume that $\chi$ and the central character of $\pi$ coincide on $A_H(F)$ (the maximal split central torus in $H(F)$). Then, we have

$$\displaystyle m(\pi,\chi)=\int_{\Gamma_{\elli}(\overline{H})}D^H(x)\Theta_\pi(x)\chi(x)^{-1}dx$$

\noindent where $\Theta_\pi$ denotes the Harish-Chandra character of $\pi$ (a locally constant function on $G_{\reg}(F)$), $D^H$ is the usual Weyl discriminant and $\Gamma_{\elli}(\overline{H})$ stands for the set of regular elliptic conjugacy classes in $\overline{H}(F):=H(F)/A_H(F)$ equipped with a suitable measure $dx$.
\end{theointro}

\vspace{1mm}

\noindent Theorem \ref{inner forms} is an easy consequence of this formula and Theorem \ref{Steinberg} also follows from it with some extra work. Let us give an outline of the proof of Theorem \ref{Steinberg} assuming Theorem \ref{formula}. For notational simplicity we will assume that $H$ is semi-simple. We have the following explicit formula for the character of the Steinberg representation (see \S \ref{section Steinberg} for a reminder)

$$\displaystyle D^G(x)^{1/2}\Theta_{\St(G)}(x)=\sum_{P_0\subseteq P=MU} (-1)^{a_P-a_{P_0}} \sum_{\{y\in M(F);\; y\sim_{conj} x \}/M-conj}D^M(y)^{1/2} \delta_P(y)^{1/2}$$

\noindent where $P_0$ is a minimal parabolic subgroup of $G$ and we refer the reader to the core of the paper for other unexplained notations which are however pretty standard. Plugging this explicit formula in Theorem \ref{formula} and rearranging somewhat the terms we get 

$$\displaystyle m(\St(G),\chi)=\sum_{(\mathcal{M},P)/conj} (-1)^{a_{P}-a_{P_0}} \int_{\Gamma_{\elli}(\mathcal{M})} D^{\mathcal{M}}(x) \chi(x)^{-1}dx$$

\noindent where the sum runs over the $H(F)$-conjugacy classes of pairs $(\mathcal{M},P)$ with

\begin{itemize}
\item $\mathcal{M}$ an {\it elliptic twisted Levi subgroup} of $H$ by which we mean an algebraic subgroup of $H$ with trivial split center such that $R_{E/F}\mathcal{M}_E$ is a Levi subgroup of $G$;
\item $P$ a parabolic subgroup of $G$ with Levi component $R_{E/F}\mathcal{M}_E$.
\end{itemize}

\noindent Using a particular case of Harish-Chandra orthogonality relations between characters of discrete series (\cite{Clo} Theorem 3), we can show that (see \S \ref{section orthogonal relations})

$$\displaystyle \int_{\Gamma_{\elli}(\mathcal{M})} D^{\mathcal{M}}(x) \chi(x)^{-1}dx=(\chi_{\mid \mathcal{M}},\mathbf{1})$$

\noindent where $\chi_{\mid \mathcal{M}}$ denotes the restriction of $\chi$ to $\mathcal{M}(F)$, $\mathbf{1}$ the trivial character of $\mathcal{M}(F)$ and $(.,.)$ denotes the natural scalar product on the space of virtual characters of $\mathcal{M}(F)$. Then, in the above expression for $m(\St(G),\chi)$, we can group together pairs $(\mathcal{M},P)$ according to their {\it stable} conjugacy classes ending up with an equality

$$\displaystyle m(\St(G),\chi)=\sum_{(\mathcal{M},P)/stab} (-1)^{a_{P}-a_{P_0}} \lvert \ker1^1(F;\mathcal{M},H)\rvert (\chi_{\mid \mathcal{M}},\mathbf{1})$$

\noindent where $\ker1^1(F;\mathcal{M},H):=\Ker\left(H^1(F,\mathcal{M})\to H^1(F,H) \right)$, a set which naturally parametrizes conjugacy classes inside the stable conjugacy class of $(\mathcal{M},P)$. Set $H_{ab}$ for the quotient of $H(F)$ by the common kernel of all the characters of Galois type (in case $H$ is quasi-split it is just the {\it abelianization} of $H(F)$) and let $\mathcal{M}_{ab}$ denote, for all elliptic twisted Levi $\mathcal{M}$, the image of $\mathcal{M}(F)$ in $H_{ab}$. Then, using Frobenius reciprocity, the last identity above can be rewritten as the equality between $m(\St(G),\chi)$ and

$$\displaystyle \left(\chi, \sum_{(\mathcal{M},P)/stab} (-1)^{a_{P}-a} \left\lvert \ker1^1(F;\mathcal{M},H)\right\rvert Ind_{\mathcal{M}_{ab}}^{H_{ab}}(\mathbf{1}) \right)$$

\noindent and thus Theorem \ref{Steinberg} is now equivalent to the following identity in the Grothendieck group of $H_{ab}$:

\begin{align}\label{eq 0}
\displaystyle \sum_{(\mathcal{M},P)/stab} (-1)^{a_{P}-a} \left\lvert \ker1^1(F;\mathcal{M},H)\right\rvert Ind_{\mathcal{M}_{ab}}^{H_{ab}}(\mathbf{1})=\omega_{H,E}
\end{align}

\noindent The proof of this identity in general is rather long and technical (see Proposition \ref{prop twisted Levi}), so we content ourself (again) with giving two examples here:

\begin{itemize}
\item If $H=GL_n$, we have $H_{ab}=F^\times$ and $\omega_{H,E}=\eta_{E/F}^{n+1}$. If $n$ is odd, $\mathcal{M}=H$ is the only elliptic twisted Levi and then \ref{eq 0} reduces to $\mathbf{1}=\mathbf{1}$. On the other hand, if $n$ is even there are two (stable) conjugacy classes of pairs $(\mathcal{M},P)$:

$$\displaystyle \mathcal{M}_0=P_0=H \mbox{ and } \mathcal{M}_1=GL_{n/2}(E)\subseteq P_1=\begin{pmatrix}GL_{n/2}(E) & \ast \\ & GL_{n/2}(E)\end{pmatrix}$$

we have $\mathcal{M}_{0,ab}=H_{ab}=F^\times\supset \mathcal{M}_{1,ab}=N(E^\times)$ and \ref{eq 0} reduces to the identity

$$\displaystyle Ind_{N(E^\times)}^{F^\times} \mathbf{1}-\mathbf{1}=\eta_{E/F}$$

\item For $H=U(n)$ (a unitary group of rank $n$) we have $H_{ab}=\Ker N_{E/F}$ and $\omega_{H,E}=\mathbf{1}$. In this case, stable conjugacy classes of pairs $(\mathcal{M},P)$ are parametrized by (ordered) partitions $(n_1,\ldots,n_k)$ of $n$ as follows:

$$\displaystyle (n_1,\ldots,n_k)\mapsto \mathcal{M}=U(n_1)\times\ldots\times U(n_k)\subseteq P=\begin{pmatrix} GL_{n_1}(E) & \ast & \ast \\  & \ddots & \ast \\ & & GL_{n_k}(E)\end{pmatrix}$$

Moreover, $\lvert \ker1^1(F;\mathcal{M},H)\rvert=2^{k-1}$ and $\mathcal{M}_{ab}=H_{ab}$ for all $\mathcal{M}$ as above. Thus, in this case \ref{eq 0} reduces to the following combinatorial identity

$$\displaystyle \sum_{\substack{(n_1,\ldots,n_k) \\ n_1+\ldots+n_k=n}} (-1)^{n-k}2^{k-1}=1$$

\end{itemize}

\vspace{2mm}

As we said, for its part, Theorem \ref{formula} is a consequence of a certain simple local trace formula adapted to the situation and to the proof of which most of the paper is devoted. Let us state briefly the content of this formula by assuming again, for simplicity, that the group $H$ is semi-simple. Starting with a function $f\in C_c^\infty(G(F))$, we consider the following expression in two variables

$$\displaystyle K^\chi_f(x,y):=\int_{H(F)} f(x^{-1}hy) \chi(h)^{-1} dh,\;\;\; x,y\in G(F)$$

\noindent This function is precisely the kernel of the operator on $L^2(H(F)\backslash G(F),\chi)$ given by convolution by $f$. Formally, the trace of such an operator should be given by the integral of this kernel over the diagonal that is

$$\displaystyle J^\chi(f):=\int_{H(F)\backslash G(F)} K_f^\chi(x,x)dx$$

\noindent Unfortunately, in general the convolution operator given by $f$ isn't of trace-class and the above expression diverges. Nevertheless, we can still restrict our attention to some `good' space of test functions for which the above integral is absolutely convergent. Recall, following Waldspurger \cite{Wal2}, that the function $f$ is said to be {\it strongly cuspidal} if for all proper parabolic subgroups $P=MU\subsetneq G$ we have

$$\displaystyle \int_{U(F)} f(xu)du=0$$

\noindent for all $x\in M(F)$. We also say, following Harish-Chandra \cite{HCVd}, that $f$ is a {\it cusp form} if the above kind of integrals vanish for all $x\in G(F)$ (thus, and contrary to what we might guess, being a cusp form is stronger than being strongly cuspidal). Actually, it will be more convenient for us to work with functions that are not necessarily compactly supported: we will take $f$ in the so-called {\it Harish-Chandra-Schwartz space} (see \S \ref{function spaces} for a reminder) denoted by $\mathcal{C}(G(F))$. The notions of strong cuspidality and of cusp forms extend verbatim to this bigger space. The following theorem, whose proof is scattered all over this paper (see Theorem \ref{theo conv}, Theorem \ref{theo spec} and Theorem \ref{theo geom side}), is our main technical result:

\begin{theointro}
Let $f\in \mathcal{C}(G(F))$ be a strongly cuspidal function. Then, the expression defining $J^\chi(f)$ is absolutely convergent (see Theorem \ref{theo conv}) and we have:
\begin{enumerate}[(i)]
\item (see Theorem \ref{theo geom side}) A geometric expansion

$$\displaystyle J^\chi(f)=\int_{\Gamma_{\elli}(H)}D^H(x)\Theta_f(x)\chi(x)^{-1}dx$$

\noindent where the function $\Theta_f$ is defined using weighted orbital integrals of Arthur (see \S \ref{strongly cuspidal});
\item (see Theorem \ref{theo spec}) If $f$ is moreover a cusp form, a spectral expansion

$$\displaystyle J^\chi(f)=\sum_{\pi\in \Irr_{\sqr}(G)}m(\pi,\chi) \Tr \pi^\vee(f)$$

\noindent where $\Irr_{\sqr}(G)$ denotes the set of (equivalence classes of) irreducible square-integrable representations of $G(F)$ and for all $\pi\in \Irr_{\sqr}(G)$, $\pi^\vee$ is the smooth contragredient of $\pi$.
\end{enumerate}
\end{theointro}

We prove this theorem by following closely the general method laid down by \cite{Wal2}, \cite{Wal3} and \cite{B1}. In particular, a crucial point to get the spectral expansion in the above theorem is to show that for $\pi$ square-integrable the abstract multiplicity $m(\pi,\chi)$ is also the multiplicity of $\pi$ in the discrete spectrum of $L^2(H(F)\backslash G(F))$. This fact is established in the course of the proof of Proposition \ref{prop intertwinings} using the simple adaptation of an idea that goes back to Sakellaridis-Venkatesh (\cite{SV} Theorem 6.4.1) and Waldspurger (\cite{Wal3} Proposition 5.6).

Here is an outline of the contents of the different parts of the paper. In the first part, we set up the main notations and conventions as well as collect different results that will be needed in the subsequent sections. It includes in particular a discussion of a natural generalization of Arthur's $(G,M)$-families to symmetric pairs that we call $(G,M,\theta)$-families. The second part contains the proof of the absolute convergence of $J^\chi(f)$ for strongly cuspidal functions $f$ and in the third part  we establish a spectral expansion of this distribution when $f$ is a cusp form. These two parts are actually written in the more general setting of {\it tempered} symmetric pairs $(G,H)$ (which were called {\it strongly discrete} in \cite{GO}) to which the proofs extend verbatim. The fourth part deals with the geometric expansion of $J^\chi(f)$. There we really have to restrict ourself to the setting of Galois pairs (that is when $G=R_{E/F}H_E$) since a certain equality of Weyl discriminants (see \ref{eq 4.1.1}), which is only true in this particular case, plays a crucial role in allowing to control the uniform convergence of certain integrals. Finally, in the last part of this paper we prove the formula for the multiplicity (Theorem \ref{formula}) and give two applications of it towards Prasad's conjecture (Theorem \ref{inner forms} and Theorem \ref{Steinberg}). 

\section{Preliminaries}

\subsection{Groups, measures, notations}\label{Groups, measures, notations}

Throughout this paper we will let $F$ be a $p$-adic field (i.e. a finite extension of $\mathbf{Q}_p$ for a certain prime number $p$) for which we will fix an algebraic closure $\overline{F}$. We will denote by $\lvert .\rvert$ the canonical absolute value on $F$ as well as its unique extension to $\overline{F}$. Unless specified otherwise, all groups and varieties that we consider in this paper will be tacitly assumed to be defined over $F$ and we will identify them with their points over $\overline{F}$. Moreover for every finite extension $K$ of $F$ and every algebraic variety $X$ defined over $K$ we will denote by $R_{K/F}X$ Weil's restriction of scalars (so that in particular we have a canonical identification $(R_{K/F}X)(F)=X(K)$).

\vspace{2mm}

Let $G$ be a connected reductive group over $F$ and $A_G$ be its maximal central split torus. We set $\overline{G}:=G/A_G$ and

$$\displaystyle \mathcal{A}_G:=X_*(A_G)\otimes \mathbf{R}$$

\noindent where $X_*(A_G)$ denotes the abelian group of cocharacters of $A_G$. If $V$ is a real vector space we will always denote by $V^*$ its dual. The space $\mathcal{A}_G^*$ can naturally be identified with $X^*(A_G)\otimes \mathbf{R}=X^*(G)\otimes \mathbf{R}$ where $X^*(A_G)$ and $X^*(G)$ stand for the abelian groups of algebraic characters of $A_G$ and $G$ respectively. More generally, for every extension $K/F$ we will denote by $X^*_K(G)$ the group of characters of $G$ defined over $K$. There is a natural morphism $H_G:G(F)\to \mathcal{A}_G$ characterized by

$$\displaystyle \langle \chi,H_G(g)\rangle=\log(\lvert \chi(g)\rvert)$$

\noindent for all $\chi\in X^*(G)$. We set $\mathcal{A}_{G,F}:=H_G(A_G(F))$. It is a lattice in $\mathcal{A}_G$. The same notations will be used for the Levi subgroups of $G$ (i.e. the Levi components of parabolic subgroups of $G$): if $M$ is a Levi subgroup of $G$ we define similarly $A_M$, $\mathcal{A}_M$, $H_M$ and $\mathcal{A}_{M,F}$. We will also use Arthur's notations: $\mathcal{P}(M)$, $\mathcal{F}(M)$ and $\mathcal{L}(M)$ will stand for the sets of parabolic subgroups with Levi component $M$, parabolic subgroups containing $M$ and Levi subgroups containing $M$ respectively. Let $K$ be a maximal special compact subgroup of $G(F)$. Then, for all parabolic subgroups $P$ with Levi decomposition $P=MU$ the Iwasawa decomposition $G(F)=M(F)U(F)K$ allows to extend $H_M$ to a map $H_P:G(F)\to \mathcal{A}_M$ defined by $H_P(muk):=H_M(m)$ for all $m\in M(F)$, $u\in U(F)$ and $k\in K$. For all Levi subgroups $M\subset L$ there is a natural decomposition

$$\displaystyle \mathcal{A}_M=\mathcal{A}_M^L\oplus \mathcal{A}_L$$ 

\noindent where $\mathcal{A}_M^L$ is generated by $H_M(\Ker(H_{L\mid M(F)}))$ and we will set $a_M^L:=\dim(\mathcal{A}_M^L)$. The Lie algebra of $G$ will be denoted by $\mathfrak{g}$ and more generally for any algebraic group we will denote its Lie algebra by the corresponding Gothic letter. We will write $\Ad$ for the adjoint action of $G$ on $\mathfrak{g}$. We denote by $\exp$ the exponential map which is an $F$-analytic map from an open neighborhood of $0$ in $\mathfrak{g}(F)$ to $G(F)$. For all subsets $S\subset G$, we write $\Cent_G(S)$, $\Cent_{G(F)}(S)$ and $\No_{G(F)}(S)$ for the centralizer of $S$ in $G$, resp. the centralizer of $S$ in $G(F)$, resp. the normalizer of $S$ in $G(F)$. If $S=\{x\}$ we will denote by $G_x$ the neutral connected component of $\Cent_G(x):=\Cent_G(\{x\})$. We define $G_{\reg}$ as the open subset of regular semisimple elements of $G$ and for all subgroups $H$ of $G$ we will write $H_{\reg}:=H\cap G_{\reg}$. Recall that a regular element $x\in G_{\reg}(F)$ is said to be {\it elliptic} if $A_{G_x}=A_G$. We will denote by $G(F)_{\elli}$ the set of regular elliptic elements in $G(F)$. The Weyl discriminant $D^G$ is defined by

$$\displaystyle D^G(x):=\left\lvert det(1-\Ad(x)_{\mid \mathfrak{g}/\mathfrak{g}_x})\right\rvert$$

\noindent For every subtorus $T$ of $G$ we will write

$$\displaystyle W(G,T):=\No_{G(F)}(T)/\Cent_{G(F)}(T)$$

\noindent for its Weyl group. If $A\subset G$ is a split subtorus we will denote by $R(A,G)$ the set of roots of $A$ in $G$ i.e. the set of nontrivial characters of $A$ appearing in the action of $A$ on $\mathfrak{g}$. More generally, if $H$ is a subgroup of $G$ and $A\subset H$ is a split subtorus we will denote by $R(A,H)$ the set of roots of $A$ in $H$.

\vspace{2mm}

In all this paper we will assume that all the groups that we encounter have been equipped with Haar measures (left and right invariant as we will only consider measures on unimodular groups). In the particular case of tori we normalize these Haar measure by requiring that they give mass $1$ to their maximal compact subgroups. For any Levi subgroup $M$ of $G$ we equip $\mathcal{A}_M$ with the unique Haar measure such that $\vol(\mathcal{A}_M/\mathcal{A}_{M,F})=1$. If $M\subset L$ are two Levi subgroups then we give $\mathcal{A}_M^L\simeq \mathcal{A}_M/\mathcal{A}_L$ the quotient measure.

\vspace{2mm}

Finally, we will adopt the following slightly imprecise but convenient notations. If $f$ and $g$ are positive functions on a set $X$, we will write

\begin{center}
$f(x)\ll g(x)$ for all $x\in X$
\end{center}

\noindent and we will say that $f$ is essentially bounded by $g$, if there exists a $c>0$ such that

\begin{center}
$f(x)\leqslant cg(x)$, for all $x\in X$
\end{center}

\noindent We will also say that $f$ and $g$ are equivalent and we will write

\begin{center}
$f(x)\sim g(x)$ for all $x\in X$
\end{center}

\noindent if both $f$ is essentially bounded by $g$ and $g$ is essentially bounded by $f$.

\subsection{log-norms}\label{log-norms}

All along this paper, we will assume that $\mathfrak{g}(F)$ has been equipped with a (classical) norm $\lvert .\rvert_{\mathfrak{g}}$, that is a map $\lvert .\rvert_{\mathfrak{g}}:\mathfrak{g}(F)\to \mathbb{R}_+$ satisfying $\lvert \lambda X\rvert_{\mathfrak{g}}=\lvert \lambda\rvert. \lvert X\rvert_{\mathfrak{g}}$, $\lvert X+Y\rvert_{\mathfrak{g}}\leqslant \lvert X\rvert_{\mathfrak{g}} +\lvert Y\rvert_{\mathfrak{g}}$ and $\lvert X\rvert_{\mathfrak{g}}=0$ if and only if $X=0$ for all $\lambda\in F$ and $X,Y\in \mathfrak{g}(F)$. For any $R>0$, we will denote by $B(0,R)$ the closed ball of radius $R$ centered at the origin in $\mathfrak{g}(F)$.

\vspace{2mm}

In this paper we will freely use the notion of log-norms on varieties over $F$. The concept of norm on varieties over local fields has been introduced by Kottwitz in \cite{Kott1} \S 18. A log-norm is essentially just the logarithm of a Kottwitz's norm and we refer to \cite{B1} \S 1.2 for the basic properties of these log-norms. For convenience, we collect here the definition and basic properties of these objects.

\noindent First, an abstract log-norm on a set $X$ is just a real-valued function $x\mapsto \sigma(x)$ on $X$ such that $\sigma(x)\geqslant 1$, for all $x\in X$. For two abstract log-norms $\sigma_1$ and $\sigma_2$ on $X$, we will say that $\sigma_2$ dominates $\sigma_1$ and we will write $\sigma_1\ll \sigma_2$ if
 
$$\sigma_1(x)\ll \sigma_2(x)$$

\noindent for all $x\in X$. We will say that $\sigma_1$ and $\sigma_2$ are equivalent if each of them dominates the other.

\vspace{2mm}

\noindent For an affine algebraic variety $X$ over $\overline{F}$, choosing a set of generators $f_1,\ldots,f_m$ of its $\overline{F}$-algebra of regular functions $\mathcal{O}(X)$, we can define an abstract log-norm $\sigma_X$ on $X$ by setting

$$\sigma_X(x)=1+\log\left(max\{1,\lvert f_1(x)\rvert,\ldots,\lvert f_m(x)\rvert\}\right)$$

\noindent for all $x\in X$. The equivalence class of $\sigma_X$ doesn't depend on the particular choice of $f_1,\ldots,f_m$ and by a log-norm on $X$ we will mean any abstract log-norm in this equivalence class. Note that if $U$ is the principal Zariski open subset of $X$ defined by the non-vanishing of $Q\in \mathcal{O}(X)$, then we have

$$\sigma_U(x)\sim \sigma_X(x)+\log\left(2+\lvert Q(x)\rvert^{-1}\right)$$

\noindent for all $x\in U$. More generally, for $X$ any algebraic variety over $\overline{F}$, choosing a finite covering $\left(U_i\right)_{i\in I}$ of $X$ by open affine subsets and fixing log-norms $\sigma_{U_i}$ on each $U_i$, we can define an abstract log-norm on $X$ by setting

$$\sigma_X(x)=\inf\{\sigma_{U_i}(x);i\in I \; \mbox{ such that } x\in U_i\}$$

\noindent Once again, the equivalence class of $\sigma_X$ doesn't depend on the various choices and by a log-norm on $X$ we will mean any abstract log-norm in this equivalence class.

\vspace{2mm}

\noindent We will assume that all varieties that we consider in this paper are equipped with log-norms and we will set $\sigma:=\sigma_G$ and $\overline{\sigma}:=\sigma_{\overline{G}}$.

\vspace{2mm}

\noindent Let $p:X\to y$ be a regular map between algebraic varieties then we have

$$\displaystyle \sigma_Y(p(x))\ll \sigma_X(x)$$

\noindent for all $x\in X$. If $p$ is a closed immersion or more generally if $p$ is a finite morphism (\cite{Kott1} Proposition 18.1(1)) we have

$$\displaystyle \sigma_Y(p(x))\sim \sigma_X(x)$$

\noindent for all $x\in X$. We say that $p$ has the {\it norm descent property} (with respect to $F$) if, denoting by $p_F$ the induced map on $F$-points, we have

$$\displaystyle \sigma_Y(y)\sim \inf_{x\in p_F^{-1}(y)}\sigma_X(x)$$

\noindent for all $y\in p_F(X(F))$. By Proposition 18.3 of \cite{Kott1}, if $T$ is a subtorus of $G$ then the projection $G\twoheadrightarrow T\backslash G$ has the norm descent property i.e. we have

\begin{align}\label{eq 1.2.1}
\displaystyle \sigma_{T\backslash G}(g)\sim \inf_{t\in T(F)} \sigma(tg)
\end{align}

\noindent for all $g\in G(F)$. In section \ref{Estimates} we will prove that for $H$ an $F$-spherical subgroup of $G$ (i.e. a subgroup such that there exists a minimal parabolic subgroup $P_0$ of $G$ with $HP_0$ open) the projection $G\twoheadrightarrow H\backslash G$ also has the norm descent property.

\vspace{2mm}

Let $T\subset G$ again be a maximal subtorus. As the regular map $T\backslash G\times T_{\reg}\to G_{\reg}$, $(g,t)\mapsto g^{-1}tg$ is finite we get

\begin{align}\label{eq 1.2.2}
\displaystyle \sigma_{T\backslash G}(g)\ll \sigma(g^{-1}tg)\log\left(2+D^G(t)^{-1}\right)
\end{align}

\noindent for all $g\in G$ and all $t\in T_{\reg}$.

\vspace{2mm}

For every variety $X$ defined over $F$, equipped with a log-norm $\sigma_X$, and all $M>0$ we will denote by $X[<M]$, resp. $X[\geqslant M]$, the set of all $x\in X(F)$ such that $\sigma_X(x)<M$, resp. $\sigma_X(x)\geqslant M$. With this notation, if $T$ is a torus over $F$ and $k:=dim(A_T)$ we have

\begin{align}\label{eq 1.2.3}
\displaystyle \vol\left(T[<M]\right)\ll M^k
\end{align}

\noindent for all $M>0$.

\subsection{Function spaces}\label{function spaces}

Let $\omega$ be a continuous character of $A_G(F)$. We define $\mathcal{S}_\omega(G(F)):=C_c^\infty(A_G(F)\backslash G(F),\omega)$ as the space of functions $f:G(F)\to \mathbf{C}$ which are smooth (i.e. locally constant), satisfy $f(ag)=\omega(a)f(g)$ for all $(a,g)\in A_G(F)\times G(F)$ and are compactly supported modulo $A_G(F)$. 

Assume moreover that $\omega$ is unitary and let $\Xi^G$ be Harish-Chandra Xi function associated to a special maximal compact subgroup $K$ of $G(F)$ (see \cite{Wal1} \S II.1). Then, we define the Harish-Chandra-Schwartz space $\mathcal{C}_\omega(G(F))$ as the space of functions $f:G(F)\to \mathbf{C}$ which are biinvariant by an open subgroup of $G(F)$, satisfy $f(ag)=\omega(a)f(g)$ for all $(a,g)\in A_G(F)\times G(F)$ and such that for all $d>0$ we have $\lvert f(g)\rvert\ll \Xi^G(g)\overline{\sigma}(g)^{-d}$ for all $g\in G(F)$.

\subsection{Representations}\label{section representations}

In this paper all representations we will consider are smooth and we will always use the slight abuse of notation of identifying a representation $\pi$ with the space on which it acts. We will denote by $\Irr(G)$ the set of equivalence classes of smooth irreducible representations of $G(F)$ and by $\Irr_{\cusp}(G)$, $\Irr_{\sqr}(G)$ the subsets of supercuspidal and essentially square-integrable representations respectively. If $\omega$ is a continuous unitary character of $A_G(F)$ we will also write $\Irr_\omega(G)$ (resp. $\Irr_{\omega,\cusp}(G)$, $\Irr_{\omega,\sqr}(G)$) for the sets of all $\pi\in \Irr(G)$ (resp. $\pi\in \Irr_{\cusp}(G)$, $\pi\in \Irr_{\sqr}(G)$) whose central character restricted to $A_G(F)$ equals $\omega$. For all $\pi\in \Irr(G)$ we will denote by $\pi^\vee$ its contragredient and by $\langle .,.\rangle:\pi\times \pi^\vee\to \mathbf{C}$ the canonical pairing. For all $\pi\in \Irr_{\omega,\sqr}(G)$, $d(\pi)$ will stand for the formal degree of $\pi$. Recall that it depends on the Haar measure on $G(F)$ and that it is uniquely characterized by the following identity (Schur orthogonality relations)

$$\displaystyle \int_{A_G(F)\backslash G(F)} \langle \pi(g)v_1,v^\vee_1\rangle \langle v_2,\pi^\vee(g)v_2^\vee\rangle dg=\frac{1}{d(\pi)} \langle v_1,v_2^\vee\rangle \langle v_2,v_1^\vee\rangle$$

\noindent for all $v_1,v_2\in \pi$ and all $v_1^\vee,v_2^\vee\in \pi^\vee$. From this, we easily infer that for every coefficient $f$ of $\pi$ we have

\begin{align}\label{eq 1.4.1}
\displaystyle \Tr(\pi^\vee(f))=\frac{1}{d(\pi)} f(1)
\end{align}

\vspace{2mm}

Let $\pi\in \Irr(G)$ and let $\omega$ be the inverse of the restriction of the central character of $\pi$ to $A_G(F)$. Then, for all $f\in \mathcal{S}_\omega(G(F))$ we can define an operator $\pi(f)$ on $\pi$ by

$$\displaystyle \langle \pi(f)v,v^\vee\rangle:=\int_{A_G(F)\backslash G(F)} f(g)\langle \pi(g)v,v^\vee\rangle dg$$

\noindent for all $(v,v^\vee)\in \pi\times \pi^\vee$. For all $f\in \mathcal{S}_\omega(G(F))$ this operator is of finite rank and a very deep theorem of Harish-Chandra (\cite{HCDeBS} Theorem 16.3) asserts that the distribution

$$\displaystyle f\in \mathcal{S}_\omega(G(F))\mapsto \Tr(\pi(f))$$

\noindent is representable by a locally integrable function which is locally constant on $G_{\reg}(F)$. This function, the Harish-Chandra character of $\pi$, will be denoted $\Theta_\pi$. It is characterized by

$$\displaystyle \Tr(\pi(f))=\int_{A_G(F)\backslash G(F)} \Theta_\pi(g) f(g)dg$$

\noindent for all $f\in \mathcal{S}_\omega(G(F))$. If moreover the representation $\pi$ is square-integrable (or more generally {\it tempered}) with unitary central character, then the integral defining $\pi(f)$ still makes sense for all $f\in \mathcal{C}_\omega(G(F))$, the resulting operator is again of finite rank and the above equality continues to hold.

\subsection{Weighted orbital integrals}

Let $M$ be a Levi subgroup and fix a maximal special compact subgroup $K$ of $G(F)$. Using $K$ we can define maps $H_P:G(F)\to\mathcal{A}_M$ for all $P\in \mathcal{P}(M)$ (cf \S \ref{Groups, measures, notations}). Let $g\in G(F)$. The family

$$\displaystyle \left\{-H_P(g); P\in \mathcal{P}(M) \right\}$$

\noindent is a positive $(G,M)$-orthogonal set in the sense of Arthur (see \cite{Art0} \S 2). In particular, following {\it loc. cit.} using this family we can define a weight $v_M^Q(g)$ for all $Q\in \mathcal{F}(M)$. Concretely, $v_M^Q(g)$ is the volume of the convex hull of the set $\{ H_P(g);\; P\in \mathcal{P}(M), P\subset Q\}$ (this convex hull belongs to a certain affine subspace of $\mathcal{A}_M$ with direction $\mathcal{A}^L_M$ where $Q=LU$ with $M\subset L$ and we define the volume with respect to the fixed Haar measure on $\mathcal{A}^L_M$). If $Q=G$ we set $v_M(g):=v^G_M(g)$ for simplicity. For every character $\omega$ of $A_G(F)$, every function $f\in \mathcal{C}_\omega(G(F))$ and all $x\in M(F)\cap G_{\reg}(F)$, we define, again following Arthur, a weighted orbital integral by

$$\displaystyle \Phi_M^Q(x,f):=\int_{G_x(F)\backslash G(F)} f(g^{-1}xg)v_M^Q(g)dg$$

\noindent The integral is absolutely convergent by the following lemma which is an immediate consequence of \ref{eq 1.2.2} and Lemma \ref{lemma HC Clo} (which will be proved later).

\begin{lem}
Let $x\in M(F)\cap G_{\reg}(F)$. Then, for all $d>0$ there exists $d'>0$ such that the integral
$$\displaystyle \int_{G_x(F)\backslash G(F)} \Xi^G(g^{-1}xg)\overline{\sigma}(g^{-1}xg)^{-d'}\sigma_{G_x\backslash G}(x)^ddg$$
converges.
\end{lem}

\noindent Once again if $Q=G$, we will set $\Phi_M(x,f):=\Phi^G_M(x,f)$ for simplicity. If $M=G$ (so that necessarily $Q=G$), $\Phi_G(x,f)$ reduces to the usual orbital integral.

\subsection{Strongly cuspidal functions}\label{strongly cuspidal}

Let $\omega$ be a continuous unitary character of $A_G(F)$. Following \cite{Wal2}, we say that a function $f\in \mathcal{C}_\omega(G(F))$ is {\it strongly cuspidal} if for every proper parabolic subgroup $P=MU$ of $G$ we have

$$\displaystyle \int_{U(F)} f(mu)du=0,\;\;\;\forall m\in M(F)$$

\noindent (the integral is absolutely convergent by \cite{Wal1} Proposition II.4.5). By a standard change of variable, $f$ is strongly cuspidal if and only if for every proper parabolic subgroup $P=MU$ and for all $m\in M(F)\cap G_{\reg}(F)$ we have

$$\displaystyle \int_{U(F)} f(u^{-1}mu)du=0$$

\noindent We will denote by $\mathcal{C}_{\omega,\scusp}(G(F))$ the subspace of strongly cuspidal functions in $\mathcal{C}_\omega(G(F))$ and we will set $\mathcal{S}_{\omega,\scusp}(G(F)):=\mathcal{S}_\omega(G(F))\cap \mathcal{C}_{\omega,\scusp}(G(F))$. Let $K$ be a maximal special compact subgroup of $G(F)$. For $x\in G_{\reg}(F)$ set $M(x):=\Cent_G(A_{G_x})$ (it is the smallest Levi subgroup containing $x$). Then, by \cite{Wal2} Lemme 5.2, for all $f\in \mathcal{C}_{\omega,\scusp}(G(F))$, all Levi subgroups $M$, all $Q\in \mathcal{F}(M)$ and all $x\in M(F)\cap G_{\reg}(F)$ we have $\Phi_M^Q(x,f)=0$ unless $Q=G$ and $M=M(x)$. For all $x\in G_{\reg}(F)$ we set

$$\displaystyle \Theta_f(x):=(-1)^{a_{M(x)}^G}\Phi_{M(x)}(x,f)$$

\noindent Then the function $\Theta_f$ is independent of the choice of $K$ and invariant by conjugation (\cite{Wal2} Lemme 5.2 and Lemme 5.3). Also by \cite{Wal2} Corollaire 5.9, the function $(D^G)^{1/2}\Theta_f$ is locally bounded on $G(F)$.

\vspace{2mm}

We say that a function $f\in \mathcal{C}_\omega(G(F))$ is a {\it cusp form} if it satisfies one of the following equivalent conditions (see \cite{Wal1} Th\'eor\`eme VIII.4.2 and Lemme VIII.2.1 for the equivalence between these two conditions):

\begin{itemize}
\item For every proper parabolic subgroup $P=MU$ and all $x\in G(F)$ we have

$$\displaystyle \int_{U(F)} f(xu)du=0;$$

\item $f$ is a sum of matrix coefficients of representations in $\Irr_{\omega,\sqr}(G)$.
\end{itemize}

\noindent We will denote by ${}^0\mathcal{C}_\omega(G(F))$ the space of cusp forms. Let $f\in {}^0\mathcal{C}_\omega(G(F))$ and set $f_\pi(g):=\Tr(\pi^\vee(g^{-1})\pi^\vee(f))$ for all $\pi\in \Irr_{\omega,\sqr}(G(F))$ and all $g\in G(F)$. Then, $f_\pi$ belongs to ${}^0\mathcal{C}_\omega(G(F))$ for all $\pi\in \Irr_{\omega,\sqr}(G(F))$ (\cite{HC} Theorem 29) and we have an equality

\begin{align}\label{eq 1.6.1}
\displaystyle f=\sum_{\pi\in \Irr_{\omega,\sqr}(G)}d(\pi)f_\pi
\end{align}

\noindent (This is a special case of Harish-Chandra-Plancherel formula, see \cite{Wal1} Theorem VIII.4.2).

Let ${}^0\mathcal{S}_\omega(G(F)):=\mathcal{S}_\omega(G(F))\cap {}^0\mathcal{C}_\omega(G(F))$ be the space of compactly supported cusp forms. Similar to the characterization of ${}^0\mathcal{C}_\omega(G(F))$, a function $f\in \mathcal{S}_\omega(G(F))$ belongs to ${}^0\mathcal{S}_\omega(G(F))$ if and only if it satisfies one the the following equivalent conditions:

\begin{itemize}
\item For every proper parabolic subgroup $P=MU$ and all $x\in G(F)$ we have

$$\displaystyle \int_{U(F)} f(xu)du=0;$$

\item $f$ is a sum of matrix coefficients of representations in $\Irr_{\omega,\cusp}(G)$.
\end{itemize}

\noindent Moreover, for $f\in {}^0\mathcal{S}_{\omega}(G(F))$, we have $f_\pi\in {}^0\mathcal{S}_{\omega}(G(F))$ for all $\pi\in \Irr_{\omega,\cusp}(G)$ and a spectral decomposition

\begin{align}\label{eq 1.6.2}
\displaystyle f=\sum_{\pi\in \Irr_{\omega,\cusp}(G)}d(\pi)f_\pi
\end{align}

Finally, we will need the following proposition.

\begin{prop}\label{prop coeff char}
Let $\pi\in \Irr_{\sqr}(G)$ and let $f$ be a matrix coefficient of $\pi$. Then, we have
$$\displaystyle \Theta_f(x)=\frac{1}{d(\pi)} f(1)\Theta_\pi(x)$$
for all $x\in G_{\reg}(F)$.
\end{prop}

\noindent\ul{Proof}: Unfortunately, the author has been unable to find a suitable reference for this probably well-known statement (however see \cite{Clo} Proposition 5 for the case where $x$ is elliptic and \cite{Art0} for the case where $\pi$ is supercuspidal). Let us say that it follows from a combination of Arthur's noninvariant local trace formula (\cite{Art2}, Proposition 4.1) applied to the case where one of the test functions is our $f$ and of Schur orthogonality relations. Note that Arthur's local trace formula was initially only proved for compactly supported test functions, but see \cite{Art3} Corollary 5.3 for the extension to Harish-Chandra Schwartz functions. $\blacksquare$

\subsection{Tempered pairs}

Let $H$ be a unimodular algebraic subgroup of $G$ (e.g. a reductive subgroup). We say that the pair $(G,H)$ is {\it tempered} if there exists $d>0$ such that the integral

$$\displaystyle \int_{H(F)} \Xi^G(h)\sigma(h)^{-d}dh$$

\noindent is convergent. This notion already appeared in \cite{GO} under the name of {\it strongly discrete} pairs. Following the referee suggestion we have decided to call these pairs {\it tempered} instead so that it is more in accordance with the notion of {\it strongly tempered} pairs introduced by Sakellaridis-Venkatesh in \cite{SV} \S 6 (since the latter implies the former but not conversely). This terminology is also justified by the fact that $(G,H)$ is tempered if and only if the Haar measure on $H(F)$ defines a tempered distribution on $G(F)$ i.e. it extends to a continuous linear form on $\mathcal{C}(G(F))$. Moreover, by a result of Benoist and Kobayashi \cite{BK}, when $H$ is reductive and in the case where $F=\mathbb{R}$ (which is not properly speaking included in this paper) a pair $(G,H)$ is tempered if and only if $L^2(H(F)\backslash G(F))$ is tempered as a unitary representation of $G(F)$. Although the author has not checked all the details, the proof of Benoist and Kobayashi seems to extend without difficulties to the $p$-adic case. However, we propose here a quick proof of one of the implications (but we won't use it in this paper).

\begin{prop}
Assume that the pair $(G,H)$ is tempered. Then, the unitary representation of $G(F)$ on $L^2(H(F)\backslash G(F))$ given by right translation is tempered i.e. the Plancherel measure of $L^2(H(F)\backslash G(F))$ is supported on tempered representations.
\end{prop}

\noindent\ul{Proof}: We will use the following criterion for temperedness due to Cowling-Haagerup-Howe \cite{CHH}:

\begin{num}
\item\label{eq 1.9.1} Let $(\Pi,\mathcal{H})$ be a unitary representation of $G(F)$. Then $(\Pi,\mathcal{H})$ is tempered if and only if there exists $d>0$ and a dense subspace $V\subset \mathcal{H}$ such that for all $u,v\in V$ we have
$$\displaystyle \lvert\left(\Pi(g)u,v\right)\rvert\ll \Xi^G(g) \sigma(g)^d$$
for all $g\in G(F)$ where $(.,.)$ denotes the scalar product on $\mathcal{H}$.
\end{num}

We will check that this criterion is satisfied for $V=C_c^\infty(H(F)\backslash G(F))\subset \mathcal{H}=L^2(H(F)\backslash G(F))$. Let $\varphi_1,\varphi_2\in C_c^\infty(H(F)\backslash G(F))$ and choose $f_1,f_2\in C_c^\infty(G(F))$ such that
$$\displaystyle \varphi_i(x)=\int_{H(F)} f_i(hx)dh$$
for $i=1,2$ and all $x\in H(F)\backslash G(F)$. Then, denoting by $R(g)$ the operator of right translation by $g$ and by $(.,.)$ the $L^2$-inner product on $L^2(H(F)\backslash G(F))$, we have
\[\begin{aligned}
\displaystyle \left(R(g)\varphi_1,\varphi_2\right) & =\int_{H(F)\backslash G(F)} \varphi_1(xg) \overline{\varphi_2(x)} dx \\
 & =\int_{H(F)\backslash G(F)} \int_{H(F)\times H(F)} f_1(h_1xg) \overline{f_2(h_2x)} dh_2dh_1dx \\
 & =\int_{H(F)\backslash G(F)} \int_{H(F)\times H(F)} f_1(h_1xg) \overline{f_2(h_2h_1x)} dh_2dh_1dx \\
 & =\int_{G(F)}\int_{H(F)} f_1(\gamma g)\overline{f_2(h\gamma)}dhd\gamma
\end{aligned}\]
for all $g\in G(F)$. Let $d>0$ that we will assume sufficiently large in what follows. As $f_1$ and $f_2$ are compactly supported, there obviously exist $C_1>0$ and $C_2>0$ such that
$$\displaystyle \lvert f_1(\gamma)\rvert\leqslant C_1\Xi^G(\gamma)\sigma(\gamma)^{-2d} \mbox{ and } \lvert f_2(\gamma)\rvert\leqslant C_2\Xi^G(\gamma)\sigma(\gamma)^{-d}$$
for all $\gamma\in G(F)$. It follows that for all $g\in G(F)$ we have
$$\displaystyle \left\lvert \left(R(g)\varphi_1,\varphi_2\right)\right\rvert\leqslant C_1C_2\int_{G(F)}\int_{H(F)}\Xi^G(\gamma g)\Xi^G(h\gamma)\sigma(h\gamma)^{-d}dh\sigma(\gamma g)^{-2d}d\gamma$$
Since $\sigma(\gamma_1\gamma_2)^{-1}\ll \sigma(\gamma_1)^{-1}\sigma(\gamma_2)$ for all $\gamma_1,\gamma_2\in G(F)$, this last expression is essentially bounded by
$$\displaystyle \sigma(g)^{2d}\int_{G(F)}\int_{H(F)}\Xi^G(\gamma g)\Xi^G(h\gamma)\sigma(h)^{-d}dh\sigma(\gamma)^{-d}d\gamma$$
for all $g\in G(F)$. Let $K$ be the special maximal compact subgroup used to define $\Xi^G$. Then $\Xi^G$ is invariant both on the left and on the right by $K$ and since $\sigma(k_1\gamma k_2)^{-1}\ll \sigma(\gamma)^{-1}$ for all $\gamma\in G(F)$ and $k_1,k_2\in K$, we see that the last integral above is essentially bounded by
$$\displaystyle \sigma(g)^{2d}\int_{G(F)}\int_{H(F)}\int_{K\times K}\Xi^G(\gamma k_1g)\Xi^G(hk_2\gamma)dk_1dk_2\sigma(h)^{-d}dh\sigma(\gamma)^{-d}d\gamma$$
for all $g\in G(F)$. By the `doubling principle' (\cite{Wal1} Lemme II.1.3), it follows that
$$\displaystyle \left\lvert \left(R(g)\varphi_1,\varphi_2\right)\right\rvert\ll \Xi^G(g)\sigma(g)^{2d} \int_{G(F)} \Xi^G(\gamma)^2\sigma(\gamma)^{-d}d\gamma \times \int_{H(F)} \Xi^G(h)\sigma(h)^{-d}dh$$
for all $g\in G(F)$. By \cite{Wal1} Lemme II.1.5 and the assumption that $(G,H)$ is tempered, for $d$ sufficiently large the two integrals above are absolutely convergent. Thus, the criterion of Cowling-Haagerup-Howe is indeed satisfied for $V=C_c^\infty(H(F)\backslash G(F))$ and consequently $L^2(H(F)\backslash G(F))$ is tempered. $\blacksquare$

\vspace{2mm}

Finally, we include the following easy lemma which gives an alternative characterization of tempered pairs because it is how the tempered condition will be used in this paper.

\begin{lem}\label{lem tempered}
Set $A_G^H=(A_G\cap H)^0$. The pair $(G,H)$ is tempered if and only if there exists $d>0$ such that the integral
$$\displaystyle \int_{A_G^H(F)\backslash H(F)} \Xi^G(h)\overline{\sigma}(h)^{-d}dh$$
converges.
\end{lem}

\noindent\ul{Proof}: As $\Xi^G$ is $A_G(F)$ invariant, it clearly suffices to show:

\begin{num}
\item\label{eq 1.9.2} For $d>0$ sufficiently large, we have
$$\displaystyle \overline{\sigma}(h)^{-3d}\ll \int_{A_G^H(F)} \sigma(ah)^{-3d}da\ll \overline{\sigma}(h)^{-d}$$
for all $h\in H(F)$.
\end{num}

\noindent For this, we need first to observe that

\begin{align}\label{eq 1.9.3}
\displaystyle \overline{\sigma}(h)\sim \inf_{a\in A_G^H(F)} \sigma(ah)
\end{align}

\noindent for all $h\in H(F)$. Indeed, as $A_G^H\backslash H$ is a closed subgroup of $A_G\backslash G$, this is equivalent to the fact that the projection $H\twoheadrightarrow A_G^H\backslash H$ has the norm descent property and this can be easily deduce from the existence of an algebraic subgroup $H'$ of $H$ such that the multiplication morphism $A_G^H\times H'\to H$ is surjective and finite (so that in particular $H'(F)$ has a finite number of orbits in $A_G^H(F)\backslash H(F)$).

\vspace{2mm}

By the inequalities $\overline{\sigma}(h)\ll \sigma(ah)$ and $\sigma(a)\ll \sigma(ah)\sigma(h)$ for all $a\in A_G^H(F)$ and all $h\in H(F)$, for any $d>0$  we have
\[\begin{aligned}
\displaystyle \int_{A_G^H(F)} \sigma(ah)^{-3d}da \ll   \overline{\sigma}(h)^{-2d}\int_{A_G^H(F)} \sigma(ah)^{-d}da \ll \sigma(h)^d\overline{\sigma}(h)^{-2d}\int_{A_G^H(F)} \sigma(a)^{-d}da
\end{aligned}\]
for all $h\in H(F)$. For $d$ sufficiently large, the last integral above is absolutely convergent. Moreover, as the left hand side of the above inequality is clearly invariant by $h\mapsto ah$ for any $a\in A_G^H(F)$, by \ref{eq 1.9.3}, for $d$ sufficiently large we get
$$\displaystyle \displaystyle \int_{A_G^H(F)} \sigma(ah)^{-3d}da\ll \left(\inf_{a\in A_G^H(F)} \sigma(ah) \right)^d \overline{\sigma}(h)^{-2d}\ll \overline{\sigma}(h)^{-d}$$
for all $h\in H(F)$ and this shows one half of \ref{eq 1.9.2}. On the other hand, by the inequality $\sigma(ah)\ll \sigma(a)\sigma(h)$ for all $a\in A_G^H(F)$ and $h\in H(F)$, for any $d>0$ we have
$$\displaystyle \sigma(h)^{-3d}\int_{A_G^H(F)} \sigma(a)^{-3d}da\ll \int_{A_G^H(F)} \sigma(ah)^{-3d}da$$
for all $h\in H(F)$. Once again, for $d$ sufficiently large the two integrals above are absolutely convergent and, as the right hand side of the inequality is invariant by $h\mapsto ah$ for any $a\in A_G^H(F)$, by \ref{eq 1.9.3} for $d$ sufficiently large we get
$$\displaystyle \overline{\sigma}(h)^{-3d}\ll \left(\inf_{a\in A_G^H(F)} \sigma(ah) \right)^{-3d}\ll \int_{A_G^H(F)} \sigma(ah)^{-3d}da$$
for all $h\in H(F)$ and this proves the second half of \ref{eq 1.9.2}. $\blacksquare$

\subsection{Symmetric varieties}

\subsubsection{Basic definition, $\theta$-split subgroups}\label{symmetric subgroups}

Let $H$ be an algebraic subgroup of $G$. Recall that $H$ is said to be {\it symmetric} if there exists an involutive automorphism $\theta$ of $G$ (defined over $F$) such that

$$\displaystyle (G^\theta)^0\subset H\subset G^\theta$$

\noindent where $G^\theta$ denotes the subgroup of $\theta$-fixed elements. If this is the case, we say that $H$ and $\theta$ are associated. The involution $\theta$ is not, in general, determined by $H$ but by \cite{HW} Proposition 1.2, its restriction to the derived subgroup of $G$ is. From now on and until the end of \S \ref{GMT families} we fix a symmetric subgroup $H$ of $G$ and we will denote by $\theta$ an associated involutive automorphism.

\vspace{2mm}

Let $T\subset G$ be a subtorus. We say that $T$ is $\theta$-split if $\theta(t)=t^{-1}$ for all $t\in T$ and we say that it is $(\theta,F)$-split if it is $\theta$-split as well as split as a torus over $F$. For every $F$-split subtorus $A\subset G$ we will denote by $A_\theta$ the maximal $(\theta,F)$-split subtorus of $A$. A parabolic subgroup $P\subset G$ is said to be $\theta$-split if $\theta(P)$ is a parabolic subgroup opposite to $P$. If this is the case, $HP$ is open, for the Zariski topology, in $G$ (this is because $\mathfrak{h}+\mathfrak{p}=\mathfrak{g}$) and similarly $H(F)P(F)$ is open, for the analytic topology, in $G(F)$. If $P$ is a $\theta$-split parabolic subgroup, we will say that the Levi component $M:=P\cap \theta(P)$ of $P$ is a $\theta$-split Levi subgroup. Note that this terminology can be slightly confusing since a torus can be a $\theta$-split Levi without being $\theta$-split as a torus (e.g. for $G=GL(2)$, $T$ the standard maximal torus and $\theta$ given by $\theta(g)=\begin{pmatrix} & 1 \\ 1 & \end{pmatrix}g\begin{pmatrix} & 1 \\ 1 & \end{pmatrix}$). Nevertheless, the author believe that no confusion should arise in this paper as the context will clarify which notion is being used.

Actually, a Levi subgroup $M\subset G$ is $\theta$-split (i.e. it is the $\theta$-split Levi component of a $\theta$-split parabolic) if and only if $M$ is the centralizer of a $(\theta,F)$-split subtorus if and only if $M$ is the centralizer of $A_{M,\theta}$. We will adapt Arthur's notation to $\theta$-split Levi and parabolic subgroups as follows: if $M$ is a $\theta$-split Levi subgroup we will denote by $\mathcal{P}^\theta(M)$, resp. $\mathcal{F}^\theta(M)$, resp. $\mathcal{L}^\theta(M)$ the set of all $\theta$-split parabolic subgroups with Levi component $M$, resp. containing $M$, resp. the set of all $\theta$-split Levi subgroups containing $M$.

\vspace{2mm}

Let $M\subset G$ be a $\theta$-split Levi subgroup. We set

$$\displaystyle \mathcal{A}_{M,\theta}:=X_*(A_{M,\theta})\otimes \mathbf{R}$$

\noindent and $a_{M,\theta}:=\dim \mathcal{A}_{M,\theta}$. Note that we have

$$\displaystyle \mathcal{A}_M=\mathcal{A}_{M,\theta}\oplus \mathcal{A}_M^\theta$$

\noindent where as before a $\theta$ superscript indicates the subset of $\theta$-fixed points. This decomposition is compatible with the decompositions

$$\displaystyle \mathcal{A}_M=\mathcal{A}_M^L\oplus \mathcal{A}_L$$

\noindent for all $L\in \mathcal{L}^\theta(M)$. Hence, we also have

$$\displaystyle \mathcal{A}_{M,\theta}=\mathcal{A}_{M,\theta}^L\oplus \mathcal{A}_{L,\theta}$$

\noindent for all $L\in \mathcal{L}^\theta(M)$, where we have set $\mathcal{A}_{M,\theta}^L:=\mathcal{A}_{M,\theta}\cap \mathcal{A}_{M}^L$. Also, we let $a_{M,\theta}^L:=\dim \mathcal{A}_{M,\theta}^L=a_{M,\theta}-a_{L,\theta}$. We define an homomorphism $H_{M,\theta}:M(F)\to \mathcal{A}_{M,\theta}$ as the composition of the homomorphism $H_M$ with the projection $\mathcal{A}_M\twoheadrightarrow \mathcal{A}_{M,\theta}$. For all $P\in \mathcal{P}^\theta(M)$, the roots $R(A_{M,\theta},U_P)$ of $A_{M,\theta}$ in the unipotent radical $U_P$ of $P$ can be considered as elements of the dual space $\mathcal{A}_{M,\theta}^*$ of $\mathcal{A}_{M,\theta}$. There is a unique subset $\Delta_{P,\theta}\subset R(A_{M,\theta},U_P)$ such that every element of $(A_{M,\theta},U_P)$ is in an unique way a nonnegative integral linear combination of elements of $\Delta_{P,\theta}$. The set $\Delta_{P,\theta}$ is the image of $\Delta_P$ by the natural projection $\mathcal{A}_M^*\twoheadrightarrow \mathcal{A}_{M,\theta}^*$ and it forms a basis of $(\mathcal{A}_{M,\theta}^G)^*$. We call it the set of simple roots of $A_{M,\theta}$ in $P$. Define

$$\displaystyle \mathcal{A}_{P,\theta}^+:=\{X\in \mathcal{A}_M;\; \langle \alpha,X\rangle>0 \; \forall \alpha\in \Delta_{P,\theta} \}$$

\noindent Then, we have the decomposition

$$\displaystyle \mathcal{A}_{M,\theta}=\bigsqcup_{Q\in \mathcal{F}^\theta(M)} \mathcal{A}_{Q,\theta}^+$$

\noindent More precisely the set $R(A_{M,\theta},G)$ of roots of $A_{M,\theta}$ in $G$ divides $\mathcal{A}_{M,\theta}$ into certain facets which are exactly the cones $\mathcal{A}_{Q,\theta}^+$ where $Q\in \mathcal{F}^\theta(M)$. In a similar way, the subspaces supporting the facets of this decomposition are precisely the subspaces of the form $\mathcal{A}_{L,\theta}$, $L\in \mathcal{L}^\theta(M)$, whereas the chambers (i.e. the open facets) are precisely the cones $\mathcal{A}_{P,\theta}^+$ for $P\in \mathcal{P}^\theta(M)$. Fixing a maximal special compact subgroup $K$ of $G(F)$, for all $P\in \mathcal{P}^\theta(M)$, we define a map

$$\displaystyle H_{P,\theta}:G(F)\to \mathcal{A}_{M,\theta}$$

\noindent as the composition of $H_P$ with the projection $\mathcal{A}_M\twoheadrightarrow \mathcal{A}_{M,\theta}$ i.e. we have $H_{P,\theta}(muk)=H_{M,\theta}(m)$ for all $m\in M(F)$, $u\in U_P(F)$ and $k\in K$.

\vspace{2mm}

Let $A_0$ be a maximal $(\theta,F)$-split subtorus and let $M_0$ be its centralizer in $G$. For simplicity we set $\mathcal{A}_0:=\mathcal{A}_{M_0,\theta}$. It is known that the set of roots $R(A_0,G)$ of $A_0$ in $G$ forms a root system in the dual space $\left(\mathcal{A}_0^G \right)^*$ to $\mathcal{A}_0^G$ (Proposition 5.9 of \cite{HW}). The Weyl group associated to this root system is naturally isomorphic to

$$\displaystyle W(G,A_0):=\No_{G(F)}(A_0)/M_0(F)$$

\noindent and is called the little Weyl group (associated to $A_0$) (again Proposition 5.9 of \cite{HW}). Two maximal $(\theta,F)$-split subtori are not necessarily $H(F)$-conjugate (e.g. for $G=GL_n$ and $H=O(n)$) but they are always $G(F)$-conjugate \footnote{Indeed if $A_0$ and $A'_0$ are two maximal $(\theta,F)$-split tori, $M_0:=\Cent_G(A_0)$,$M'_0:=\Cent_G(A'_0)$, $P_0\in \mathcal{P}^\theta(M_0)$ and $P'_0\in \mathcal{P}^\theta(M'_0)$ then by \cite{HW} Proposition 4.9, $P_0$ and $P'_0$ are conjugated by an element of $g\in G(F)\cap HP_0$ and it suffices to show that we can take $g$ in $G(F)\cap HM_0$ but this follows from the fact that $H\cap P_0=H\cap M_0$}.

\vspace{2mm}

Let $M$ be a $\theta$-split Levi subgroup and let $\alpha\in R(A_{M,\theta},G)$. Then we define a `coroot' $\alpha^\vee\in \mathcal{A}^G_{M,\theta}$ as follows. First assume that $\alpha$ is a reduced root (i.e. $\frac{\alpha}{2}\notin R(A_{M,\theta},G)$). Let $M_\alpha$ be the unique $\theta$-split Levi containing $M$ such that $\mathcal{A}_{M_\alpha,\theta}=\Ker(\alpha)$. Let $Q_\alpha$ be the unique $\theta$-split parabolic subgroup of $M_\alpha$ with $\theta$-split Levi $M$ such that $\Delta_{Q_\alpha}=\{\alpha \}$. Let $P_0^{M_\alpha}$ be a minimal $\theta$-split parabolic subgroup of $M_\alpha$ contained in $Q_\alpha$ and set $M_0:=P_0^{M_\alpha}\cap \theta(P_0^{M_\alpha})$, $A_0:=A_{M_0,\theta}$. Let $\Delta_0^{M_\alpha}$ be the set of simple roots of $A_0$ in $P_0^{M_\alpha}$. Then there is an unique simple root $\beta\in \Delta_0^{M_\alpha}$ whose projection to $\mathcal{A}_{M,\theta}^*$ equals $\alpha$. Let $\beta^\vee\in \mathcal{A}_0$ be the corresponding coroot. Then we define $\alpha^\vee$ as the image of $\beta^\vee$ by the projection $\mathcal{A}_0\twoheadrightarrow \mathcal{A}_{M,\theta}$. We easily check that this construction does not depend on the choice of $P_0^{M_\alpha}$ since for another choice ${P'_0}^{,M_\alpha}$ with $M_0':={P_0'}^{,M_\alpha}\cap \theta({P'_0}^{,M_\alpha})$ and $A_0':=A_{M_0',\theta}$ there exists $m\in M(F)$ with $mA_0'm^{-1}=A_0$ and $m{P'_0}^{,M_\alpha}m^{-1}=P_0^{M_\alpha}$. If $\alpha$ is nonreduced, there exists $\alpha_0\in R(_{M,\theta},G)$ such that $\alpha=2\alpha_0$ and we simply set $\alpha^\vee=\frac{\alpha_0^\vee}{2}$.

Let $\widetilde{\alpha}\in R(A_M,G)$ be a root extending $\alpha$ and $\widetilde{\alpha}^\vee\in \mathcal{A}_M$ the corresponding coroot. Then, in general the projection of $\widetilde{\alpha}^\vee$ to $\mathcal{A}_{M,\theta}$ does not coincide with $\alpha^\vee$ as defined above but, however, the two are always positively proportional. Finally, we remark that when $M$ is a minimal $\theta$-split Levi subgroup, so that $R(A_{M,\theta},G)$ is a root system, then for all $\alpha\in R(A_{M,\theta},G)$, $\alpha^\vee$ coincides with the usual coroot defined using this root system. 

\vspace{2mm}

Let $P$ be a $\theta$-split parabolic subgroup. We set $\mathcal{A}_{P,\theta}:=\mathcal{A}_{M,\theta}$ and $a_{P,\theta}=a_{M,\theta}$ where $M:=P\cap\theta(P)$ and we let $\Delta_{P,\theta}^\vee\subseteq \mathcal{A}^G_{M,\theta}$ be the set of simple coroots corresponding to $\Delta_{P,\theta}\subseteq (\mathcal{A}^{G}_{M,\theta})^*$ and $\widehat{\Delta}_{P,\theta}\subseteq (\mathcal{A}^{G}_{M,\theta})^*$ be the basis dual to $\Delta_{P,\theta}^\vee$. More generally, let $Q\supset P$ be another $\theta$-split Levi subgroup. We set $\mathcal{A}_{P,\theta}^Q:=\mathcal{A}_{M,\theta}^L$ and $a_{P,\theta}^Q:=a_{M,\theta}^L$ where $L:=Q\cap \theta(Q)$ and we let $\Delta^Q_{P,\theta}\subseteq (\mathcal{A}^Q_{P,\theta})^*$ be the set of simple roots of $A_{M,\theta}$ in $P\cap L$, $(\Delta^Q_{P,\theta})^\vee\subseteq \mathcal{A}^Q_{P,\theta}$ be the corresponding set of simple coroots and $\widehat{\Delta}_{P,\theta}^Q\subseteq (\mathcal{A}_{P,\theta}^Q)^*$ be the basis dual to $(\Delta_{P,\theta}^Q)^\vee$. We have decompositions

$$\displaystyle \mathcal{A}_{P,\theta}=\mathcal{A}_{P,\theta}^Q\oplus \mathcal{A}_{Q,\theta},\;\;\; \mathcal{A}^*_{P,\theta}=(\mathcal{A}_{P,\theta}^Q)^*\oplus \mathcal{A}_{Q,\theta}^*$$

\noindent for which $\Delta_{P,\theta}^Q\subseteq \Delta_{P,\theta}$, $(\Delta_{P,\theta}^Q)^\vee\subseteq \Delta^\vee_{P,\theta}$ and moreover $\Delta_{Q,\theta}$ (resp. $\Delta_{Q,\theta}^\vee$) is the image of $\Delta_{P,\theta}-\Delta_{P,\theta}^Q$ (resp. $\Delta_{P,\theta}^\vee-(\Delta_{P,\theta}^Q)^\vee$) by the projection $\mathcal{A}^*_{P,\theta}\twoheadrightarrow \mathcal{A}_{Q,\theta}^*$ (resp. $\mathcal{A}_{P,\theta}\twoheadrightarrow \mathcal{A}_{Q,\theta}$). We define the following functions:

\begin{itemize}
\item $\tau_{P,\theta}^Q$: characteristic function of the set of $X\in \mathcal{A}_{M,\theta}$ such that $\langle \alpha,X\rangle>0$ for all $\alpha\in \Delta_{P,\theta}^Q$;

\item $\widehat{\tau}_{P,\theta}^Q$: characteristic function of the set of $X\in \mathcal{A}_{M,\theta}$ such that $\langle \varpi_\alpha,X\rangle>0$ for all $\varpi_\alpha\in \widehat{\Delta}_{P,\theta}^Q$;

\item $\delta^Q_{M,\theta}$: characteristic function of the subset $\mathcal{A}_{L,\theta}$ of $\mathcal{A}_{M,\theta}$.
\end{itemize}

\noindent We also define a function $\Gamma_{P,\theta}^Q$ on $\mathcal{A}_{M,\theta}\times \mathcal{A}_{M,\theta}$, whose utility will be revealed in the next section, by

$$\displaystyle \Gamma_{P,\theta}^Q(H,X):=\sum_{R\in \mathcal{F}^\theta(M); P\subseteq R\subseteq Q} (-1)^{a_{R,\theta}-a_{Q,\theta}} \tau_{P,\theta}^R(H)\widehat{\tau}_{R,\theta}^Q(H-X)$$

\vspace{2mm}

Let $M$ be a $\theta$-split Levi subgroup. Then, for all $P\in \mathcal{P}^\theta(M)$ we set

$$\displaystyle \mathcal{A}^{+,*}_{P,\theta}:=\{\lambda\in \mathcal{A}_{P,\theta}^*;\; \langle \alpha^\vee,\lambda\rangle>0\; \forall \alpha^\vee\in \Delta_{P,\theta}^\vee \}$$

\noindent As for $\mathcal{A}_{M,\theta}$, the set of coroots $R(A_{M,\theta},G)^\vee$ divides $\mathcal{A}^*_{M,\theta}$ into facets which are exactly the cones $\mathcal{A}_{Q,\theta}^{+,*}$ for $Q\in \mathcal{F}^\theta(M)$ and the chambers for this decomposition are the $\mathcal{A}_{P,\theta}^{+,*}$, $P\in \mathcal{P}^\theta(M)$. As usual, we say that two parabolics $P,P'\in \mathcal{P}^\theta(M)$ are adjacent if the intersection of the closure of their corresponding chambers contains a facet of codimension one. If this is the case, the hyperplane generated by this intersection is called the wall separating the two chambers.

\subsubsection{$(G,M,\theta)$-families and orthogonal sets}\label{GMT families}

As we have recalled in the previous section, the combinatorics of $\theta$-split Levi and parabolic subgroups is entirely governed, as is the case for classical Levi and parabolic subgroups, by a root system. As a consequence, for $M$ a $\theta$-split Levi subgroup of $G$ the classical theory of $(G,M)$-families due to Arthur extends without difficulty to a theory of $(G,M,\theta)$-families indexed by $\theta$-split parabolics that we now introduce. By definition, a $(G,M,\theta)$-family is a family $(\varphi_{P,\theta})_{P\in \mathcal{P}^\theta(M)}$ of $C^\infty$ functions on $i\mathcal{A}_{M,\theta}^*$ such that for any two adjacent parabolic subgroups $P,P'\in \mathcal{P}^\theta(M)$, the functions $\varphi_{P,\theta}$ and $\varphi_{P',\theta}$ coincide on the wall separating the chambers $i\mathcal{A}_{P,\theta}^{+,*}$ and $i\mathcal{A}_{P',\theta}^{+,*}$. To a $(G,M,\theta)$-family $(\varphi_{P,\theta})_{P\in \mathcal{P}^\theta(M)}$ we can associate a scalar $\varphi_{M,\theta}$ as follows: the function

$$\displaystyle \varphi_{M,\theta}(\lambda):=\sum_{P\in \mathcal{P}^\theta(M)} \varphi_{P,\theta}(\lambda) \varepsilon_{P,\theta}(\lambda),\;\;\;\lambda\in i\mathcal{A}_{M,\theta}^*$$

\noindent where we have set

$$\displaystyle \varepsilon_{P,\theta}(\lambda):=meas\left( \mathcal{A}_{M,\theta}^G/\mathbf{Z}[\Delta_{P,\theta}^\vee]\right)\prod_{\alpha^\vee\in \Delta_{P,\theta}^\vee} \langle \lambda,\alpha^\vee\rangle^{-1}$$

\noindent is $C^\infty$ and we define $\varphi_{M,\theta}:=\varphi_{M,\theta}(0)$. Note that we need a Haar measure on $\mathcal{A}_{M,\theta}^G$ for the definition of the functions $\varepsilon_{P,\theta}$ ($P\in \mathcal{P}^\theta(M)$) to make sense. We fix one as follows. Let $\overline{A}_{M,\theta}$ be the image of $A_{M,\theta}$ in $\overline{G}:=G/A_G$ and let $A_{M,\theta}'$ be the inverse image of $\overline{A}_{M,\theta}$ in $G$. Let $\mathcal{A}_{M,\theta,F}^G\subseteq \mathcal{A}_{M,\theta}^G$ denote the image of $A'_{M,\theta}(F)$ by $H_{M,\theta}^G$. It is a lattice of $\mathcal{A}_{M,\theta}^G$ and we choose our measure so that the quotient $\mathcal{A}_{M,\theta}^G/\mathcal{A}_{M,\theta,F}^G$ is of measure one.

\vspace{2mm}

We will actually only need $(G,M,\theta)$-families of a very particular shape obtained as follows. We say that a family of points $\mathcal{Y}_{M,\theta}=(Y_{P,\theta})_{P\in \mathcal{P}^\theta(M)}$ in $\mathcal{A}_{M,\theta}$ is a $(G,M,\theta)$-{\it orthogonal set} if for all adjacent $P,P'\in \mathcal{P}^\theta(M)$ we have

$$\displaystyle Y_{P,\theta}-Y_{P',\theta}=r_{P,P'}\alpha^\vee$$

\noindent where $r_{P,P'}\in \mathbf{R}$ and $\{\alpha^\vee\}=\Delta^\vee_{P,\theta}\cap -\Delta^\vee_{P',\theta}$. We say that the family is a {\it positive $(G,M,\theta)$-orthogonal set} if it is a $(G,M,\theta)$-orthogonal set and moreover $r_{P,P'}\geqslant 0$ for all adjacent $P,P'\in \mathcal{P}^\theta(M)$. To a $(G,M,\theta)$-orthogonal set $\mathcal{Y}_{M,\theta}=(Y_{P,\theta})_{P\in \mathcal{P}^\theta(M)}$ we associate the $(G,M,\theta)$-family $(\varphi_{P,\theta}(.,\mathcal{Y}_{M,\theta}))_{P\in \mathcal{P}^\theta(M)}$ defined by

$$\displaystyle \varphi_{P,\theta}(\lambda,\mathcal{Y}_{M,\theta}):=e^{\langle \lambda, Y_{P,\theta}\rangle},\;\;\; \lambda\in i\mathcal{A}_{M,\theta}^*$$

\noindent and we let $v_{M,\theta}(\mathcal{Y}_{M,\theta}):=\varphi_{M,\theta}(0,\mathcal{Y}_{M,\theta})$ be the scalar associated to this $(G,M,\theta)$-family. If $\mathcal{Y}_{M,\theta}$ is a positive $(G,M,\theta)$-orthogonal set then $v_{M,\theta}(\mathcal{Y}_{M,\theta})$ is just the volume of the convex hull of the elements in the family $\mathcal{Y}_{M,\theta}$ (with respect to the fixed Haar measure on $\mathcal{A}_{M,\theta}^G$). For any $(G,M,\theta)$-orthogonal set $\mathcal{Y}_{M,\theta}=(Y_{P,\theta})_{P\in \mathcal{P}^\theta(M)}$ and all $Q\in \mathcal{F}^\theta(M)$ we define $Y_{Q,\theta}$ to be the projection of $Y_{P,\theta}$ to $\mathcal{A}_{Q,\theta}$ for any $P\in \mathcal{P}^\theta(M)$ with $P\subset Q$ (the result is independent of the choice of $P$) and more generally for any $Q,R\in \mathcal{F}^\theta(M)$ with $Q\subset R$ we let $Y_{Q,\theta}^R$ be the projection of $Y_{P,\theta}$ to $\mathcal{A}_{Q,\theta}^R$ for any $P\in \mathcal{P}^\theta(M)$ with $P\subset Q$. Then, for any $L\in \mathcal{L}^\theta(M)$ the family $\mathcal{Y}_{L,\theta}:=(Y_{Q,\theta})_{Q\in \mathcal{P}^\theta(L)}$ forms a $(G,L,\theta)$-orthogonal set.

Fixing a maximal special compact subgroup $K$ of $G(F)$ to define maps $H_{P,\theta}$ ($P\in \mathcal{P}^\theta(M)$), for all $g\in G(F)$ the family $\mathcal{Y}_{M,\theta}(g):=(-H_{P,\theta}(g))_{P\in \mathcal{P}^\theta(M)}$ is a positive $(G,M,\theta)$-orthogonal set. Indeed, the family $(-H_P(g))_{P\in \mathcal{P}(M)}$ is a positive $(G,M)$-orthogonal set in the classical sense of Arthur (see \cite{Art0} \S 2) and thus for all $P=MU_P,P'=MU_{P'}\in \mathcal{P}^\theta(M)$ we have
$$\displaystyle -H_P(g)+H_{P'}(g)\in \sum_{\alpha\in R(A_M,U_P)\cap -R(A_M,U_{P'})}\mathbf{R}_+ \alpha^\vee.$$
As the projection of $R(A_M,U_P)$ (resp. $R(A_M,U_{P'})$) to $\mathcal{A}_{M,\theta}^*$ is $R(A_{M,\theta},U_P)$ (resp. $R(A_{M,\theta},U_{P'})$) and for all $\alpha\in R(A_M,G)$ the projection of the coroot $\alpha^\vee$ to $\mathcal{A}_{M,\theta}*$ is positively proportional to $\overline{\alpha}^\vee$, where $\overline{\alpha}\in R(A_{M,\theta},G)$ denotes the projection of $\alpha$ to $\mathcal{A}_{M,\theta}^*$, it follows that
$$\displaystyle -H_{P,\theta}(g)+H_{P',\theta}(g)\in \sum_{\alpha\in R(A_{M,\theta},U_P)\cap -R(A_{M,\theta},U_{P'})}\mathbf{R}_+ \alpha^\vee$$
for all $P,P'\in \mathcal{P}^\theta(M)$ i.e. $\mathcal{Y}_{M,\theta}(g)$ is a positive $(G,M,\theta)$-orthogonal set.

We define

$$\displaystyle v_{M,\theta}(g):=v_{M,\theta}(\mathcal{Y}_{M,\theta}(g))$$

\noindent There is another easier way to obtain $(G,M,\theta)$-orthogonal sets. It is as follows. Let $M_0\subset M$ be a minimal $\theta$-split Levi subgroup with little Weyl group $W_0$. Fix $P_0\in \mathcal{P}^\theta(M_0)$. Then, for all $X\in \mathcal{A}_0:=\mathcal{A}_{M_0,\theta}$ we define a $(G,M_0,\theta)$-orthogonal set $\mathcal{Y}[X]_0:=(Y[X]_{P_0',\theta})_{P_0'\in \mathcal{P}^\theta(M_0)}$ by setting $Y[X]_{P_0',\theta}:=w_{P_0'}X$ for all $P_0'\in \mathcal{P}^\theta(M_0)$ where $w_{P_0'}\in W_0$ is the unique element such that $w_{P_0'}P_0=P_0'$. By the general construction explained above this also yields a $(G,M,\theta)$-orthogonal set $\mathcal{Y}[X]_{M,\theta}=(Y[X]_{P,\theta})_{P\in \mathcal{P}^\theta(M)}$.

\vspace{2mm}

Let $\mathcal{Y}_{M,\theta}=(Y_{P,\theta})_{P\in \mathcal{P}^\theta(M)}$ be a $(G,M,\theta)$-orthogonal set. For $Q\in \mathcal{F}^\theta(M)$ we define a function $\Gamma_{M,\theta}^Q(.,\mathcal{Y}_{M,\theta})$ on $\mathcal{A}_{M,\theta}$ by

$$\displaystyle \Gamma_{M,\theta}^Q(H,\mathcal{Y}_{M,\theta}):=\sum_{R\in \mathcal{F}^\theta(M); R\subset Q} \delta_{M,\theta}^R(H)\Gamma_{R,\theta}^Q(H,Y_{R,\theta})$$

\noindent where the functions $\delta_{M,\theta}^R$ and $\Gamma_{R,\theta}^Q$ have been defined in the previous section. Let $L=Q\cap \theta(Q)$. Fixing a norm $\lvert .\rvert$ on $\mathcal{A}_{M,\theta}^L$, we have the following basic property concerning the support of this function (see \cite{LW} Corollaire 1.8.5):

\vspace{3mm}

\begin{num}
\item\label{eq 1.7.1} There exists $c>0$ independent of $\mathcal{Y}_{M,\theta}$ such that for all $H\in \mathcal{A}_{M,\theta}$ with $\Gamma_{M,\theta}^Q(H,\mathcal{Y}_{M,\theta})\neq 0$ we have $\lvert H^Q\rvert\leqslant c\sup_{P\in \mathcal{P}^\theta(M); P\subset Q} \lvert Y_{P,\theta}^Q\rvert$ where $H^Q$ denotes the projection of $H$ to $\mathcal{A}_{M,\theta}^L$.
\end{num}

\vspace{3mm}

\noindent Moreover, if the $(G,M,\theta)$-orthogonal set $\mathcal{Y}_{M,\theta}$ is positive then $\Gamma_{M,\theta}^Q(.,\mathcal{Y}_{M,\theta})$ is just the characteristic function of the set of $H\in \mathcal{A}_{M,\theta}$ such that $H^Q$ belongs to the convex hull of $(Y_{P,\theta}^Q)_{P\in \mathcal{P}^\theta(M); P\subset Q}$ (\cite{LW} Proposition 1.8.7). Without assuming the positivity of our $(G,M,\theta)$-orthogonal set, we have the identity (\cite{LW} Lemme 1.8.4(3))

\begin{align}\label{eq 1.7.2}
\displaystyle \sum_{Q\in \mathcal{F}^\theta(M)} \Gamma_{M,\theta}^{Q}(H,\mathcal{Y}_{M,\theta})\tau_{Q,\theta}^G(H-Y_{Q,\theta})=1
\end{align}

\noindent for all $H\in \mathcal{A}_{M,\theta}$.

Let $\mathcal{R}$ be a free $\mathbf{Z}$-module of finite type. Recall that a {\it exponential-polynomial} on $\mathcal{R}$ is a function on $\mathcal{R}$ of the following form

$$\displaystyle f(Y)=\sum_{\chi\in \widehat{\mathcal{R}}}\chi(Y) p_\chi(Y)$$

\noindent where $\widehat{\mathcal{R}}$ denotes the group of complex (not necessarily unitary) characters of $\mathcal{R}$ and for all $\chi\in \widehat{\mathcal{R}}$, $p_\chi$ is a `complex polynomial' function on $\mathcal{R}$, i.e. an element of $Sym((\mathbf{C}\otimes \mathcal{R})^*)$,which is zero for all but finitely many $\chi\in \widehat{\mathcal{R}}$. If $f$ is an exponential-polynomial on $\mathcal{R}$ then a decomposition as above is unique, the set of characters $\chi\in \widehat{\mathcal{R}}$ such that $p_\chi\neq 0$ is called the {\it set of exponents} of $f$ and $p_1$ (corresponding to $\chi=1$ the trivial character) is called the {\it purely polynomial part} of $f$. Finally, we define the {\it degree} of $f$ as the maximum, over all $\chi\in \widehat{\mathcal{R}}$, of the degree of the polynomials $p_\chi$. We record the following lemma whose proof is elementary:

\begin{lem}\label{lemma polynomial-exponential}
Let $\mathcal{R}$ be a free $\mathbf{Z}$-module of finite type, let $f$ be a exponential-polynomial on $\mathcal{R}$ and let $C\subset \mathbf{R}\otimes \mathcal{R}$ be an open cone (with a vertex possibly different from the origin). Then,if the limit

$$\displaystyle \lim\limits_{\substack{Y\in \mathcal{R}\cap C \\ \lvert Y\rvert\to \infty}} f(Y)$$

\noindent exists it equals the constant term of the purely polynomial part of $f$.
\end{lem}

Let $M_0\subset M$ be a minimal $\theta$-split Levi subgroup, $\mathcal{A}_0:=\mathcal{A}_{M_0,\theta}$ and fix $P_0\in \mathcal{P}^\theta(M_0)$. For all $X\in \mathcal{A}_0$ we dispose of the $(G,M,\theta)$-orthogonal set $\mathcal{Y}[X]_{M,\theta}=(Y[X]_{P,\theta})_{P\in \mathcal{P}^\theta(M)}$ defined above. We let $\mathcal{Y}_{M,\theta}+\mathcal{Y}[X]_{M,\theta}:=(Y_{P,\theta}+Y[X]_{M,\theta})_{P\in \mathcal{P}^\theta(M)}$ be the sum of the two $(G,M,\theta)$-orthogonal sets $\mathcal{Y}_{M,\theta}$ and $\mathcal{Y}[X]_{M,\theta}$. Obviously, it is also a $(G,M,\theta)$-orthogonal set. We let

$$\displaystyle \widetilde{v}_{M,\theta}(\mathcal{Y}_{M,\theta}+\mathcal{Y}[X]_{M,\theta}):=\int_{A_G(F)\backslash A'_{M,\theta}(F)} \Gamma_{M,\theta}^G(H_{M,\theta}(a),\mathcal{Y}_{M,\theta}+\mathcal{Y}[X]_{M,\theta})da$$

\noindent where we recall that $A'_{M,\theta}$ is the subtorus generated by $A_G$ and $A_{M,\theta}$. Let $\mathcal{A}_{0,F}$ denote the image of $A_{M_0,\theta}(F)$ by $H_{M_0,\theta}$. It is a lattice in $\mathcal{A}_0$ and we have the following lemma (combine \cite{MWl4} Lemme 1.7(ii) with equalities 1.5(2) and 1.3(7) of {\it loc.cit.}):

\begin{lem}\label{lem GMT families}
For every lattice $\mathcal{R}\subset \mathcal{A}_{0,\theta,F}\otimes \mathbf{Q}$ the function $X\in \mathcal{R}\mapsto \widetilde{v}_{M,\theta}(\mathcal{Y}_{M,\theta}+\mathcal{Y}[X]_{M,\theta})$ is an exponential-polynomial whose degree and exponents belong to finite sets which are independent of $\mathcal{Y}_{M,\theta}$. Moreover, if we denote by $\widetilde{v}_{M,\theta,0}(\mathcal{Y}_{M,\theta},\mathcal{R})$ the constant coefficient of the purely polynomial part of this exponential-polynomial there exists $c>0$ depending only on $\mathcal{R}$ such that for all $k\geqslant 1$ we have

$$\displaystyle \left\lvert \widetilde{v}_{M,\theta,0}(\mathcal{Y}_{M,\theta},\frac{1}{k}\mathcal{R})-v_{M,\theta}(\mathcal{Y}_{M,\theta})\right\rvert\leqslant c k^{-1} \left(\sup_{P\in \mathcal{P}^\theta(M)} \lvert Y_{P,\theta}\rvert\right)^{a^G_{M,\theta}}$$
\end{lem}

\vspace{2mm}

Let now $\mathcal{Y}_M=(Y_P)_{P\in \mathcal{P}(M)}$ be a usual $(G,M)$-orthogonal set. This induces a $(G,M,\theta)$-orthogonal set $\mathcal{Y}_{M,\theta}:=(Y_{P,\theta})_{P\in \mathcal{P}^\theta(M)}$ where, for all $P\in \mathcal{P}^\theta(M)$, we denote by $Y_{P,\theta}$ the projection of $Y_P$ to $\mathcal{A}_{M,\theta}$. The subspace $\mathcal{A}_{M,\theta}+\mathcal{A}_G$ of $\mathcal{A}_M$ being special in the sense of \cite{Art1}\S 7 \footnote{Indeed with the notations of {\it loc.cit.} we need to check that for every root $\beta\in R(A_{M,\theta},G)$ the sum $\sum_{\alpha\in \Sigma(\beta)} m_\alpha \alpha$ is trivial on $\mathcal{A}_M^\theta$ but this is trivial since for all $\alpha\in \Sigma(\beta)$ we have $\iota(\alpha):=-\theta(\alpha)\in \Sigma(\beta)$ and $m_{\iota(\alpha)}=m_\alpha$}, we have a descent formula (Proposition 7.1 of {\it loc.cit.})

\begin{align}\label{eq 1.7.3}
\displaystyle v_{M,\theta}(\mathcal{Y}_{M,\theta})=\sum_{L\in \mathcal{L}(M)} d^G_{M,\theta}(L)v^Q_M(\mathcal{Y}_M)
\end{align}

\noindent where for all $L\in \mathcal{L}(M)$, $Q$ is a parabolic with Levi component $L$ which depends on the choice of a generic point $\xi\in \mathcal{A}_M$ and $d^G_{M,\theta}(L)$ is a coefficient which is nonzero only if $\mathcal{A}_M^G=\mathcal{A}_{M}^{G,\theta}\oplus\mathcal{A}^L_M$. Moreover if $\mathcal{A}_M^{G,\theta}=0$ then we have $d^G_{M,\theta}(G)=1$. Let $K$ be a special maximal compact subgroup of $G(F)$ that we use to define the maps $H_P$ for $P\in \mathcal{P}(M)$ and $H_{P,\theta}$ for $P\in \mathcal{P}^\theta(M)$. Then, the formula \ref{eq 1.7.3} applied to the particular case where $Y_P=H_P(g)$ for all $P\in \mathcal{P}(M)$ and for some $g\in G(F)$ yields

\begin{align}\label{eq 1.7.4}
\displaystyle v_{M,\theta}(g)=\sum_{L\in \mathcal{L}(M)} d^G_{M,\theta}(L)v^Q_M(g)
\end{align}

\subsection{Estimates}\label{Estimates}

In this section we collect some estimates that we will need in the core of the paper. We start with four lemmas concerning maximal tori and integrals over regular orbits in $G$.

\begin{lem}\label{lemma max torus}
Let $T\subset G$ be a maximal torus . Then, we have

$$\displaystyle \sigma(t)\ll \sigma(g^{-1}tg)$$

\noindent for all $t\in T$ and all $g\in G$.
\end{lem}

\noindent\ul{Proof}: Let $W:=W(G_{\overline{F}},T_{\oF})$ be the absolute Weyl group of $T$ and set $\mathcal{B}:=G//G-\Ad$ (i.e. the GIT quotient of $G$ acting on itself by the adjoint action). Let $p:G\to \mathcal{B}$ be the natural projection. By Chevalley theorem, the inclusion $T\hookrightarrow G$ induces an isomorphism $T//W\simeq \mathcal{B}$ and thus the restriction of $p$ to $T$ is a finite morphism. Hence, we have

$$\displaystyle \sigma(t)\sim \sigma_{\mathcal{B}}(p(t))$$

\noindent for all $t\in T$ and it follows that

$$\sigma(t)\sim \sigma_{\mathcal{B}}(p(t))=\sigma_{\mathcal{B}}(p(g^{-1}tg))\ll \sigma(g^{-1}tg)$$

\noindent for all $t\in T$ and all $g\in G$. $\blacksquare$

\begin{lem}[Harish-Chandra, Clozel]\label{lemma HC Clo}
Let $T\subset G$ be a maximal torus. Then, for all $d>0$ there exists $d'>0$ such that

$$\displaystyle D^G(t)^{1/2} \int_{T(F)\backslash G(F)} \Xi^G(g^{-1}tg)\sigma(g^{-1}tg)^{-d'}dg\ll \sigma(t)^{-d}$$

\noindent for all $t\in T_{\reg}(F)$.
\end{lem}

\noindent\ul{Proof}: By Corollary 2 of \cite{Clo} there exists $d_0>0$ such that

$$\displaystyle \sup_{t\in T_{\reg}(F)}D^G(t)^{1/2}\int_{T(F)\backslash G(F)} \Xi^G(g^{-1}tg)\sigma(g^{-1}tg)^{-d_0}dg<\infty$$

\noindent Thus by Lemma \ref{lemma max torus}, for all $d>0$ we have

$$\displaystyle D^G(t)^{1/2}\int_{T(F)\backslash G(F)} \Xi^G(g^{-1}tg)\sigma(g^{-1}tg)^{-d_0-d}dg\ll \sigma(t)^{-d}$$

\noindent for all $t\in T_{\reg}(F)$. $\blacksquare$

\begin{lem}\label{lemma conv int tori}
Let $T\subset G$ be a subtorus such that $T_{\reg}:=T\cap G_{\reg}$ is nonempty (i.e. $T$ contains nonsingular elements). Then, for all $k>0$, there exists $d>0$ such that the integral

$$\displaystyle \int_{T(F)} \log(2+D^G(t)^{-1})^k \sigma(t)^{-d} dt$$

\noindent converges.
\end{lem}

\noindent\ul{Proof}: We denote by $X^*_{\oF}(T)$ the group of regular characters of $T$ defined over $\oF$. There exists a multiset $\Sigma$ of nontrivial elements in $X^*_{\oF}(T)$ such that

$$\displaystyle D^G(t)=\prod_{\alpha\in \Sigma} \lvert \alpha(t)-1\rvert$$

\noindent for all $t\in T_{\reg}(F)$ where we have denoted by $\lvert .\rvert$ the unique extension of the absolute value over $F$ to $\overline{F}$. We have

$$\displaystyle \log(2+D^G(t)^{-1})\ll \prod_{\alpha\in \Sigma}\log(2+\lvert \alpha(t)-1\rvert^{-1})$$

\noindent Thus, by Cauchy-Schwartz, it suffices to prove the following claim:

\vspace{3mm}

\begin{num}
\item\label{eq 1.8.1} For all $\alpha\in X^*_{\oF}(T)-\{ 1\}$ and all $k>0$ there exists $d>0$ such that the integral
$$\displaystyle \int_{T(F)} \log\left(2+\lvert \alpha(t)-1\rvert^{-1}\right)^k \sigma(t)^{-d}dt$$
converges.
\end{num}

\vspace{3mm}

\noindent Let $\alpha\in X^*_{\oF}(T)-\{ 1\}$ and let $\Gamma_\alpha\subset \Gamma_F$ be the stabilizer of $\alpha$ for the natural Galois action. Write $\Gamma_\alpha=\Gal(\oF/F_\alpha)$ where $F_\alpha/F$ is a finite extension. By the universal property of restriction of scalars, $\alpha$ induces a morphism $\widetilde{\alpha}:T\to R_{F_\alpha/F}\mathbf{G}_m$. Denoting by $\Ker(\widetilde{\alpha})$ the kernel of $\widetilde{\alpha}$, for all $k>0$ and all $d>0$ we have

$$\displaystyle \int_{T(F)} \log\left(2+\lvert \alpha(t)-1\rvert^{-1}\right)^k \sigma(t)^{-k}dt=\int_{T(F)/\Ker(\widetilde{\alpha})(F)} \log\left(2+\lvert \alpha(t)-1\rvert^{-1}\right)^k\int_{\Ker(\widetilde{\alpha})(F)} \sigma(tt')^{-d}dt'dt$$

\noindent As there exists a subtorus $T'\subset T$ such that the multiplication map $T'\times \Ker(\widetilde{\alpha})\to T$ is an isogeny (so that $\sqrt{\sigma(t)\sigma(t')}\ll \sigma(tt')$ and $\sigma(t)\sim \sigma_{T/\Ker(\widetilde{\alpha})}(t)$ for all $(t,t')\in \Ker(\widetilde{\alpha})\times T'$), we see that for all $d$ sufficiently large (i.e. so that the integral below converges) we have

$$\displaystyle \int_{\Ker(\widetilde{\alpha})(F)} \sigma(tt')^{-d}dt'\ll \sigma_{T/\Ker(\widetilde{\alpha})}(t)^{-d/2}$$

\noindent for all $t\in T(F)/\Ker(\widetilde{\alpha})(F)$. Since $T(F)/\Ker(\widetilde{\alpha})(F)$ is an open subset of $\left(T/\Ker(\widetilde{\alpha}) \right)(F)$ we are thus reduced to the case where $\widetilde{\alpha}$ is an embedding.

\vspace{2mm}

\noindent Define $N(\alpha)\in X^*(T)$ by

$$\displaystyle N(\alpha):=\prod_{\sigma\in \Gamma_F/\Gamma_\alpha} \sigma(\alpha)$$

\noindent We distinguish two cases. First, if $N(\alpha)\neq 1$ then we have an inequality

$$\displaystyle \log\left(2+\lvert \alpha(t)-1\rvert^{-1}\right)\ll \log\left(2+\lvert N(\alpha)(t)-1\rvert^{-1}\right)$$

\noindent for all $t\in T(F)$ with $N(\alpha)(t)\neq 1$. Hence, up to replacing $\alpha$ by $N(\alpha)$ we may assume that $\alpha\in X^*(T)$ in which case by the previous reduction we are left to prove \ref{eq 1.8.1} in the particular case where $T=\mathbf{G}_m$ and $\alpha=Id$ in which case it is easy to check. Assume now that $N(\alpha)=1$. Since $\widetilde{\alpha}$ is an embedding this implies that $T$ is anisotropic and we just need to prove that for all $k>0$ the function

$$\displaystyle t\in T(F)\mapsto \log\left(2+\lvert \alpha(t)-1\rvert^{-1}\right)^k$$

\noindent is locally integrable. Using the exponential map we are reduced to proving a similar statement for vector spaces where we replace $\alpha(t)-1$ by a linear form which is easy to check directly. $\blacksquare$

\vspace{2mm}

Combining \ref{eq 1.2.2} with Lemma \ref{lemma HC Clo} and Lemma \ref{lemma conv int tori} we get the following:

\begin{lem}\label{lemma conv abs woi}
Let $T\subset G$ be a subtorus such that $T_{\reg}:=T\cap G_{\reg}$ is nonempty and let $T^G$ be the centralizer of $T$ in $G$ (a maximal torus). Then, for all $k>0$ there exists $d>0$ such that the integral
$$\displaystyle \int_{T(F)} D^G(t)^{1/2}\int_{T^G(F)\backslash G(F)} \Xi^G(g^{-1}tg)\sigma(g^{-1}tg)^{-d}\sigma_{T^G\backslash G}(g)^kdgdt$$
converges.
\end{lem}

\vspace{2mm}

The following lemma will be needed in the proof of the next proposition. As it might be of independent interest we present it separately.

\begin{lem}\label{lem anisotropic action}
Let $G$ be an anisotropic group over $F$ and $Y$ an affine $G$-variety. Set $Y'=Y/G$ for the GIT quotient (it is an affine algebraic variety over $F$) and denote by $p:Y\to Y'$ the natural projection. Then we have

$$\displaystyle \sigma_Y(y)\sim \sigma_{Y'}(p(y))$$

\noindent for all $y\in Y(F)$.
\end{lem}

\vspace{3mm}

\noindent\ul{Proof}: First, we have $\sigma_{Y'}(p(y))\ll \sigma_Y(y)$ for all $y\in Y$ since $p$ is a morphism of algebraic varieties. Let $f\in F[Y]$, we need to show that

$$\displaystyle \log(2+\lv f(y)\rvert)\ll \sigma_{Y'}(p(y))$$

\noindent for all $y\in Y(F)$. Let $W$ be the $G$-submodule of $F[Y]$ generated by $f$ and $V$ be its dual. There is a natural morphism $\varphi:Y\to V$ and we have a commutative diagram

$$\xymatrix{ Y \ar[d]^{\varphi} \ar[r] & Y'=Y/G \ar[d]^{\varphi'} \\ V \ar[r] & V':=V/G}$$

\noindent By definition there is a function $f_V\in F[V]$ such that $f=f_V\circ \varphi$ and moreover $\sigma_{V'}(y')\ll \sigma_{Y'}(y')$ for all $y'\in Y'$. Hence, we are reduced to the case where $Y=V$ and we may assume that $f$ is homogeneous. By Kempf's extension of the stability criterion of Mumford over any perfect field (\cite{Kem}, Corollary 5.1) and since $G$ is anisotropic, for every $v\in V(F)$ the $G$-orbit $G.v\subset V$ is closed. It follows that there exist homogeneous polynomials $P_1,\ldots,P_N\in F[V']=F[V]^G$ whose only common zero in $V(F)$ is $0$. We only need to show that for some $R,C>0$ we have

\begin{align}\label{eq 1.8.3}
\displaystyle \max(1,\lv f(v)\rv)\leqslant C\max(1,\lv P_1(v)\rv,\ldots, \lv P_N(v)\rv)^R
\end{align}

\noindent for all $v\in V(F)$. Up to replacing $f,P_1,\ldots,P_N$ by some powers, we may assume that they are all of the same degree. Then, for every $1\leqslant i\leqslant N$, $f/P_i$ is a rational function on the projective space $\mathbf{P}(V)$ and the map

$$\displaystyle [v]\in \mathbf{P}(V)(F)\mapsto \min\left(\lv (f/P_1)([v])\rvert,\ldots, \lv (f/P_N)([v])\rvert\right)\in \mathbf{R}_+$$

\noindent is continuous for the analytic topology hence bounded (as $\mathbf{P}(V)(F)$ is compact) and this proves that inequality \ref{eq 1.8.3} is true for $R=1$ and some constant $C$. $\blacksquare$

\vspace{2mm}

Following \cite{KK} Definition 4.9, we say that a subgroup $H\subset G$ is {\it $F$-spherical} if there exists a minimal parabolic subgroup $P_0$ of $G$ such that $HP_0$ is open, in the Zariski topology, in $G$. For example symmetric subgroups (see \S \ref{symmetric subgroups}) are $F$-spherical. Recall that in \S \ref{log-norms} we have defined a `norm descent property' for regular maps between $F$-varieties. 

\begin{prop}\label{prop norm descent}
Let $H\subset G$ be an $F$-spherical subgroup. Then, the natural projection $p:G\to H\backslash G$ has the norm descent property.
\end{prop}

\vspace{2mm}

\noindent\ul{Proof}: Set $X:=H\backslash G$. By \cite{Kott1} Proposition 18.2 (1), it suffices to show that $X$ can be covered by Zariski open subsets over which the projection $p$ has the norm descent property. Since $G$ acts transitively on $X$ it even suffices to construct only one such open subset (because its $G$-translates will have the same property). By the local structure theorem (\cite{KK} Corollary 4.12), there exists a parabolic subgroup $Q=LU$ of $G$ such that

\begin{itemize}
\item $\mathcal{U}=HQ$ is open in $G$;
\item $H\cap Q=H\cap L$ and this subgroup contains the non-anisotropic factors of the derived subgroup of $L$.
\end{itemize}

\noindent Obviously, to show that the restriction of $p$ to $\mathcal{U}$ has the norm descent property it is sufficient to establish that $L\to H\cap L\backslash L$ has the norm descent property. We are thus reduced to the case where $H$ contains all the non-anisotropic factors of the derived subgroup of $G$. Let $G_{der}$ denote the derived subgroup of $G$, $G_{der,nc}$ the product of the non-anisotropic factors of $G_{der}$, $G_{der,c}$ the product of the anisotropic factors of $G_{der}$ and set $G'=G_{der,c}Z(G)^0$, $H'=H\cap G'$. Then, we have $H=G_{der,nc}H'$ and the multiplication map $G_{der,nc}\times G'\to G$ is an isogeny. It follows that there exists a finite set $\{\gamma_i;i\in I \}$ of elements of $G(F)$ such that

$$\displaystyle G(F)=\bigsqcup_{i\in I} G_{der,nc}(F)G'(F)\gamma_i$$

\noindent and

$$\displaystyle H(F)\backslash G(F)=\bigcup_{i\in I} H'(F)\backslash G'(F)\gamma_i$$

\noindent From these decompositions, we infer that we only need to prove the norm descent property for $G'\to H'\backslash G'$ i.e. we may assume that $G_{der,nc}=1$. Let $G_c$ be the product of $G_{der,c}$ with the maximal anisotropic subtorus of $Z(G)^0$ and consider the projection

$$\displaystyle p': X:=H\backslash G\to X':=HG_{c}\backslash G$$

\noindent We claim that

\begin{align}\label{eq 1.8.4}
\displaystyle \sigma_X(x)\sim \sigma_{X'}(p'(x))
\end{align}

\noindent for all $x\in X(F)$. As $H$ is reductive (since $G(F)$ contains no unipotent element), $X$ is affine and the claim follows from the Lemma \ref{lem anisotropic action}.

\vspace{2mm}

\noindent Now because of \ref{eq 1.8.4}, we may replace $X$ by $X'$ i.e. we may assume that $H$ contains $G_{c}$. As the multiplication map $G_{c}\times A_G\to G$ is an isogeny, by a similar argument as before we are reduced to the case where $G$ is a split torus for which the proposition is easy to establish directly. $\blacksquare$

\section{Definition of a distribution for all symmetric pairs}

\subsection{The statement}\label{statement convergence}

Let $G$ be a connected reductive group over $F$, $H$ be a symmetric subgroup of $G$ and $\theta$ be the involution of $G$ associated to $H$ (see \S \ref{symmetric subgroups}). Set $A_G^H=(A_G\cap H)^0$, $\overline{G}:=G/A_G$, $\overline{H}:=H/A_H$, $X:=A_G(F)H(F)\backslash G(F)$, $\mathbf{X}:=HA_G\backslash G$, $\sigma_X:=\sigma_{\mathbf{X}}$ and $\overline{\sigma}:=\sigma_{\overline{G}}$. Note that $X$ is an open subset of $\mathbf{X}(F)$. Let $\chi$ and $\omega$ be continuous unitary characters of $H(F)$ and $A_G(F)$ respectively such that $\chi_{\mid A_G^H(F)}=\omega_{\mid A_G^H(F)}$. Then, for all $f\in \mathcal{S}_{\omega}(G(F))$ we define a function $K^\chi_f$ on $X$ by

$$\displaystyle K^\chi_f(x):=\int_{A_G^H(F)\backslash H(F)} f(x^{-1}hx)\chi(h)^{-1}dh$$

\noindent If the pair $(G,H)$ is tempered then the expression defining $K^\chi_f$ makes sense for all $f\in \mathcal{C}_{\omega}(G(F))$ (by Lemma \ref{lem tempered}). The goal of this chapter is to show that if $f$ is strongly cuspidal then the expression

$$\displaystyle J^\chi(f):=\int_{X} K^\chi_f(x) dx$$

\noindent is convergent. More precisely we will prove the following

\begin{theo}\label{theo conv}
For all $f\in \mathcal{S}_{\omega,\scusp}(G(F))$, the expression defining $J^\chi(f)$ is absolutely convergent. Moreover, if the pair $(G,H)$ is tempered then the expression defining $J^\chi(f)$ is also absolutely convergent for all $f\in \mathcal{C}_{\omega,\scusp}(G(F))$.
\end{theo}

\subsection{Some estimates}

Let $A\subset G$ be a $(\theta,F)$-split subtorus, set $M:=\Cent_G(A)$ and let $Q=LU_Q\in \mathcal{F}^\theta(M)$ where $L:=Q\cap\theta(Q)$. Let $\overline{Q}=\theta(Q)=LU_{\overline{Q}}$ be the opposite parabolic subgroup and set
$$\displaystyle A_{\overline{Q}}^{+}:=\{a\in A;\; \lv \alpha(a)\rv \geqslant 1\; \forall \alpha\in R(A,U_{\overline{Q}}) \}$$
and
$$\displaystyle A_{\overline{Q}}^{+}(\delta):=\{a\in A;\; \lv \alpha(a)\rv \geqslant e^{\delta\overline{\sigma}(a)}\; \forall \alpha\in R(A,U_{\overline{Q}}) \}$$
for all $\delta>0$. Recall that if $Y$ is an algebraic variety over $F$ and $M>0$ we denote by $Y[<M]$ the subset of $y\in Y(F)$ with $\sigma_Y(y)<M$. We also recall that we have fixed a (classical) norm $\lvert .\rvert_{\mathfrak{g}}$ and that for any $R>0$, $B(0,R)$ denotes the closed ball of radius $R$ centered at the origin for this norm (see \S \ref{log-norms}).

\begin{lem}\label{lem estimates}
\begin{enumerate}[(i)]
\item\label{lem estimates i} Let $\epsilon>0$ and $\delta>0$. Then, we have
$$\displaystyle \overline{\sigma}(a)\ll \sup(\overline{\sigma}(g),\overline{\sigma}(a^{-1}ga))$$
for all $a\in A_{\overline{Q}}^{+}(\delta)$ and all $g\in G(F)\backslash \left(Q(F)aU_{Q}[<\epsilon \overline{\sigma}(a)]a^{-1} \right)$;

\item\label{lem estimates ii} Let $0<\delta'<\delta$ and $c_0>0$. Then, if $\epsilon>0$ is sufficiently small we have
$$\displaystyle aU_Q[<\epsilon \overline{\sigma}(a)]a^{-1}\subseteq \exp\left(B(0,c_0e^{-\delta'\overline{\sigma}(a)})\cap \mathfrak{u}_Q(F) \right)$$
for all $a\in A^+_{\overline{Q}}(\delta)$.

\item\label{lem estimates iii} We have
$$\displaystyle \overline{\sigma}(h)\ll\overline{\sigma}(a^{-1}ha)$$
and
$$\displaystyle \overline{\sigma}(h)+\overline{\sigma}(a)\ll \overline{\sigma}(ha)$$
for all $a\in A$ and all $h\in H$.

\item\label{lem estimates iv} Set $L:=Q\cap \theta(Q)$, $H_L:=H\cap L$ and $H^Q:=H_L\ltimes U_Q$ where $U_Q$ denotes the unipotent radical of $Q$. Then, we have
$$\displaystyle \overline{\sigma}(h^Q)\ll\overline{\sigma}(a^{-1}h^Qa)$$
for all $a\in A_{\overline{Q}}^+$ and all $h^Q\in H^Q$.
\end{enumerate}
\end{lem}

\vspace{2mm}

\noindent\ul{Proof}: (\ref{lem estimates i}) and (\ref{lem estimates ii}) are essentially \cite{B1} Lemma 1.3.1 (i) and (ii) applied to the group $\overline{G}:=A_G\backslash G$. To prove (\ref{lem estimates iii}), we first observe that
$$\displaystyle \theta(a^{-1}ha)=aha^{-1} \mbox{ and } \theta(ha)^{-1}ha=a^2$$
for all $a\in A$, $h\in H$. Hence, we have
$$\displaystyle \sup(\overline{\sigma}(aha^{-1}),\overline{\sigma}(a^{-1}ha))\ll \overline{\sigma}(a^{-1}ha) \mbox{ and } \overline{\sigma}(a^2)\ll \overline{\sigma}(ha)$$
for all $a\in A$, $h\in H$. Since $\overline{\sigma}(a)\sim \overline{\sigma}(a^2)$ for all $a\in A$ and $\overline{\sigma}(h)\ll \overline{\sigma}(ha)+\overline{\sigma}(a)$, this already suffices to establish the second inequality. To prove the first one, it only remains to show the following

\vspace{3mm}

\begin{num}
\item\label{eq 2.2.1} We have
$$\displaystyle \overline{\sigma}(g)\ll \sup(\overline{\sigma}(aga^{-1}),\overline{\sigma}(a^{-1}ga))$$
for all $a\in A$ and all $g\in G$.
\end{num}

\vspace{3mm}

\noindent Fix an embedding $\iota: \overline{G}\hookrightarrow SL_n$ for some $n\geqslant 1$ which sends the torus $A$ into the standard maximal torus $A_n$ of $SL_n$. Then, we are reduced to proving \ref{eq 2.2.1} in the particular case where $G=SL_n$ and $A=A_n$. For every matrix $g\in SL_n$, denote by $g_{i,j}$ ($1\leqslant i,j\leqslant n$) the $(i,j)$th-entry of $g$ and for $a\in A_n$, set $a_i=a_{i,i}$. Then, we have

$$\displaystyle \overline{\sigma}(g)\sim \sup_{i,j}\log(2+ \lv g_{i,j}\rvert)$$

\noindent Hence, it suffices to show that for all $1\leqslant i,j\leqslant n$ we have

$$\displaystyle \log(2+ \lv g_{i,j}\rvert)\leqslant \sup\left(\log(2+\lv (aga^{-1})_{i,j}\rvert),\log(2+\lv (a^{-1}ga)_{i,j}\rvert)\right)$$

\noindent for all $g\in SL_n$ and $a\in A_n$. However, $(aga^{-1})_{i,j}=a_ia_j^{-1}g_{i,j}$, $(a^{-1}ga)_{i,j}=a_ja_i^{-1}g_{i,j}$ and at least one of the quotients $a_ia_j^{-1}$, $a_ja_i^{-1}$ is of absolute value greater than $1$. The result follows. 

\vspace{2mm}

We now prove (\ref{lem estimates iv}). Every $h^Q\in H^Q(F)$ can be written $h^Q=h_Lu_Q$ where $h_L\in H_L\subset L$ and $u_Q\in U_Q$. Moreover, we have $\overline{\sigma}(lu_Q)\sim \overline{\sigma}(l)+\overline{\sigma}(u_Q)$ and $\overline{\sigma}(u_Q)\ll\overline{\sigma}(a^{-1}u_Qa)$ for all $l\in L$, $u_Q\in U_Q$ and $a\in A_{\overline{Q}}^+$. Besides, as $H_L\subset H$, by (\ref{lem estimates iii}) we have $\overline{\sigma}(h_L)\ll\overline{\sigma}(a^{-1}h_La)$ for all $h_L\in H_L$ and $a\in A$. It follows that
\[\begin{aligned}
\displaystyle \overline{\sigma}(h^Q)\sim \overline{\sigma}(h_L)+\overline{\sigma}(u_Q)\ll \overline{\sigma}(a^{-1}h_La)+\overline{\sigma}(a^{-1}u_Qa)\sim \overline{\sigma}(a^{-1}h^Qa)
\end{aligned}\]
for all $h^Q=h_Lu_Q\in H^Q=H_L\ltimes U_Q$ and all $a\in A_{\overline{Q}}^+$. This proves (\ref{lem estimates iv}). $\blacksquare$

\subsection{Weak Cartan decompositions and Harish-Chandra-Schwartz space of $X$}\label{weak Cartan and HCS for X}

Let $A_{0,j}$, $j\in J$, be representatives of the $H(F)$-conjugacy classes of maximal $(\theta,F)$-split tori of $G$. There are a finite number of them and by a result of Benoist-Oh and Delorme-S\'echerre (\cite{BO} and \cite{DS}), there exists a compact subset $\mathcal{K}_G\subset G(F)$ such that

\begin{align}\label{eq 2.3.1}
\displaystyle G(F)=\bigcup_{j\in J} H(F)A_{0,j}(F)\mathcal{K}_G
\end{align}

\noindent This decomposition is called a weak Cartan decomposition.

\vspace{2mm}

\noindent Let $C\subset G(F)$ be a compact subset with nonempty interior and set

$$\displaystyle \Xi^{X}_C(x)=\vol_X(xC)^{-1/2}$$

\noindent for all $x\in X$ and where $\vol_X$ refers to a $G(F)$-invariant measure on $X$ (which exists as $H$ is reductive hence unimodular). If $C'\subset G(F)$ is another compact subset with nonempty interior, the functions $\Xi^{X}_C$ and $\Xi^{X}_{C'}$ are equivalent and we will denote by $\Xi^X$ any such function (for some choice of $C$).

\begin{prop}\label{prop HCS for X}
\begin{enumerate}[(i)]

\item\label{prop HCS for X i} For every compact subset $\mathcal{K}\subseteq G(F)$, we have the following equivalences

\begin{align}\label{eq 2.3.2}
\displaystyle \Xi^{X}(xk)\sim\Xi^{X}(x)
\end{align}

\begin{align}\label{eq 2.3.3}
\displaystyle \sigma_{X}(xk)\sim \sigma_{X}(x)
\end{align}

\noindent for all $x\in X$ and all $k\in \mathcal{K}$.

\item\label{prop HCS for X ii} Let $A_0$ be a $(\theta,F)$-split subtorus of $G$. Then, there exists $d>0$ such that

\begin{align}\label{eq 2.3.4}
\displaystyle \Xi^G(a)\overline{\sigma}(a)^{-d}\ll \Xi^X(a)\ll \Xi^G(a)
\end{align}

\begin{align}\label{eq 2.3.5}
\displaystyle \sigma_X(a)\sim \overline{\sigma}(a)
\end{align}

\noindent for all $a\in A_0(F)$.

\item\label{prop HCS for X iii} There exists $d>0$ such that the integral

$$\displaystyle \int_{X} \Xi^X(x)^2\sigma_X(x)^{-d}dx$$

\noindent is absolutely convergent.
\end{enumerate}

\vspace{2mm}

\noindent Assume moreover that the pair $(G,H)$ is tempered, then we have

\begin{enumerate}[(i)]
\setcounter{enumi}{3}
\item\label{prop HCS for X iv} For all $d>0$ there exists $d'>0$ such that

$$\displaystyle \int_{A^H_G(F)\backslash H(F)} \Xi^G(hx) \overline{\sigma}(hx)^{-d'}dh\ll \Xi^X(x)\sigma_X(x)^{-d}$$

\noindent for all $x\in X$.

\item\label{prop HCS for X v} There exist $d>0$ and $d'>0$ such that

$$\displaystyle \int_{A^H_G(F)\backslash H(F)} \Xi^G(x^{-1}hx)\overline{\sigma}(x^{-1}hx)^{-d}dh\ll \Xi^X(x)^2 \sigma_X(x)^{d'}$$

\noindent for all $x\in X$.

\item\label{prop HCS for X vi} More generally, let $Q$ be a $\theta$-split parabolic subgroup of $G$ and set $L:=Q\cap \theta(Q)$, $H_L:=H\cap L$ and $H^Q:=H_L\ltimes U_Q$ where $U_Q$ denotes the unipotent radical of $Q$. Let $A_0\subset L$ be a maximal $(\theta,F)$-split subtorus and set
$$\displaystyle A_{\overline{Q}}^+:=\{a\in A_0(F);\; \lvert \alpha(a)\rvert\geqslant 1\; \forall \alpha\in R(A_0,U_{\overline{Q}}) \}$$
where $U_{\overline{Q}}$ denotes the unipotent radical of $\overline{Q}$. Then, $H^Q$ is a unimodular algebraic group, $(G,H^Q)$ is a tempered pair and, fixing a Haar measure $dh^Q$ on $H^Q(F)$, there exists $d>0$ and $d'>0$ such that
$$\displaystyle \int_{A^H_G(F)\backslash H^Q(F)} \Xi^G(a^{-1}h^Qa)\overline{\sigma}(a^{-1}h^Qa)^{-d}dh^Q\ll \Xi^X(a)^2\sigma_X(a)^{d'}$$
for all $a\in A_{\overline{Q}}^+$.
\end{enumerate}
\end{prop}

\vspace{2mm}

\noindent\ul{Proof}:
\begin{enumerate}[(i)]
\item is easy and left to the reader.

\item Let $M_0$ be the centralizer of $A_0$ in $G$. For all $P_0\in \mathcal{P}^\theta(M_0)$ set

$$\displaystyle A_{P_0}^+:=\{a\in A_0(F); \lv \alpha(a)\rvert\geqslant 1\; \forall \alpha\in R(A_0,P_0) \}$$

\noindent Then

\begin{align}\label{eq 2.3.6}
\displaystyle A_0(F)=\bigcup_{P_0\in \mathcal{P}^\theta(M_0)} A_{P_0}^+
\end{align}

\noindent Thus, we may fix $P_0\in \mathcal{P}^\theta(M)$, set $A^+:=A_{P_0}^+$ and prove \ref{eq 2.3.4} and \ref{eq 2.3.5} for those $a$ belonging to $A^+$. Since $P_0$ is $\theta$-split, $HP_0$ is open in $G$. The proof of (\ref{prop HCS for X ii}) is now the same as Proposition 6.7.1(ii) of \cite{B1} after replacing Proposition 6.4.1(iii) of {\it loc. cit.} by Lemma \ref{lem estimates} (\ref{lem estimates iii}).

\item The proof is exactly the same as for Proposition 6.7.1 (iii) of \cite{B1}: we use the weak Cartan decomposition \ref{eq 2.3.1} to show that $X$ has polynomial growth in the sense of \cite{Ber} and then we conclude as in loc. cit.

\item Let $A_0$ be a maximal $(\theta,F)$-split subtorus of $G$, $M_0:=\Cent_G(A_0)$ and $P_0\in \mathcal{P}^\theta(M_0)$. By the decompositions \ref{eq 2.3.1} and \ref{eq 2.3.6}, points (\ref{prop HCS for X i}) and (\ref{prop HCS for X ii}) and Lemma \ref{lem estimates}(\ref{lem estimates iii}) it suffices to show the existence of $d>0$ such that

$$\displaystyle \int_{A^H_G(F)\backslash H(F)} \Xi^G(ha)\overline{\sigma}(h)^{-d}dh\ll \Xi^G(a)$$

\noindent for all $a\in A_{P_0}^+$. Since $P_0$ is $\theta$-split, $H(F)P_0(F)$ is open in $G(F)$. It follows that if $K$ is a maximal compact subgroup of $G(F)$ by which $\Xi^G$ is right invariant there exists an open-compact subgroups $J\subset G(F)$ and $J_H\subset H(F)$ such that $J\subset J_HaKa^{-1}$ for all $a\in A_{P_0}^+$. Hence, for all $d>0$, all $k\in J$ and all $a\in A_{P_0}^+$, writing $k=k_Hak_Ga^{-1}$ with $k_H\in J_H$ and $k_G\in K$,  we have

$$\displaystyle \int_{A^H_G(F)\backslash H(F)} \Xi^G(hka)\overline{\sigma}(h)^{-d}dh=\int_{A^H_G(F)\backslash H(F)} \Xi^G(ha)\overline{\sigma}(hk_H^{-1})^{-d}dh$$

\noindent It follows that for all $d>0$ we have

$$\displaystyle \int_{A^H_G(F)\backslash H(F)} \Xi^G(ha)\overline{\sigma}(h)^{-d}dh\ll \int_{A^H_G(F)\backslash H(F)}\int_K \Xi^G(hka)dk\overline{\sigma}(h)^{-d}dh$$

\noindent for all $a\in A_{P_0}^+$ and we conclude by the `doubling principle' (see \cite{B1} Proposition 1.5.1) and the fact that the pair $(G,H)$ is tempered.

\item Once again, the proof is very similar to the proof of Proposition 6.7.1 (v) of \cite{B1} so that we shall only sketch the argument. Let $A_0$ be a maximal $(\theta,F)$-split torus of $G$ and $P_0\in \mathcal{P}^\theta(M_0)$ where $M_0:=\Cent_G(A_0)$. By the weak Cartan decomposition, (\ref{prop HCS for X i}), (\ref{prop HCS for X ii}), the inequality $\overline{\sigma}(h)\ll \overline{\sigma}(a^{-1}ha)\overline{\sigma}(a)$ and the decomposition \ref{eq 2.3.6}, we are reduced to proving the existence of $d,d'>0$ such that
$$\displaystyle \int_{A^H_G(F)\backslash H(F)} \Xi^G(a^{-1}ha)\overline{\sigma}(h)^{-d}dh\ll \Xi^G(a)^2 \overline{\sigma}(a)^{d'}$$
for all $a\in A_{P_0}^+$. Using the fact that $H(F)P_0(F)$ is open in $G(F)$ we show as in the proof of (\ref{prop HCS for X iv}) that if $K$ is a maximal compact subgroup of $G(F)$ we have
$$\displaystyle \int_{A^H_G(F)\backslash H(F)} \Xi^G(a^{-1}ha)\overline{\sigma}(h)^{-d}dh\ll \int_{A^H_G(F)\backslash H(F)} \int_{K\times K} \Xi^G(a^{-1}k_1hk_2a)\overline{\sigma}(h)^{-d}dk_1dk_2dh$$
for all $a\in A_{P_0}^+$ and then we conclude again by the `doubling principle' (see \cite{B1} Proposition 1.5.1 (vi)) and the fact that the pair $(G,H)$ is tempered.

\item The proof that $H^Q$ is unimodular is similar to the proof of Proposition 6.8.1 (ii) of \cite{B1} noticing that $\overline{Q}$ is a good parabolic subgroup with respect to $H$ (that is $\overline{Q}H$ is open in $G$) and that $H_L=H\cap \overline{Q}$. Moreover, the fact that the pair $(G,H^Q)$ is tempered and the estimate can be proved in much the same way as Proposition 6.8.1 (iv)-(vi) of \cite{B1}. Indeed, if we denote by $M_0$ the centralizer of $A_0$ in $G$, we have
$$\displaystyle A_{\overline{Q}}^+=\bigcup_{\substack{P_0\in \mathcal{P}^\theta(M_0) \\ P_0\subset Q}} A_{\overline{P}_0}^+$$
and thus, fixing $P_0\in \mathcal{P}^\theta(M_0)$ with $P_0\subset Q$, it suffices to show the existence of $d,d'>0$ such that

\vspace{3mm}

\begin{num}
\item\label{eq 2.3.7} The integral $\displaystyle \int_{A_G^H(F)\backslash H^Q(F)} \Xi^G(h^Q)\overline{\sigma}(h^Q)^{-d}dh^Q$ converges;
\end{num}

\vspace{3mm}
 
\noindent and

\vspace{3mm}

\begin{num}
\item\label{eq 2.3.8} $\displaystyle \int_{A_G^H(F)\backslash H^Q(F)} \Xi^G(a^{-1}h^Qa)\overline{\sigma}(a^{-1}h^Qa)^{-d}dh^Q\ll \Xi^X(a)^2 \sigma_X(a)^{d'}$ for all $a\in A_{\overline{P}_0}^+$.
\end{num}

\vspace{3mm}

\noindent Then, \ref{eq 2.3.7} can be proved exactly as Proposition 6.8.1 (iv) of \cite{B1} where the only inputs used are the facts that $\overline{Q}H$ is open in $G$ and that the pair $(G,H)$ is tempered. Also, \ref{eq 2.3.8} can be proved exactly as Proposition 6.8.1 (vi) of {\it loc. cit.} where this time the only inputs are the estimates of (\ref{prop HCS for X ii}) and of Lemma \ref{lem estimates} (\ref{lem estimates iv}), the convergence of \ref{eq 2.3.7} and the fact that $\overline{P}_0H$ is open in $G$. $\blacksquare$

\end{enumerate}

\subsection{Proof of theorem \ref{theo conv}}

By Proposition \ref{prop HCS for X} (\ref{prop HCS for X iii}), it suffices to establish the two following claims

\vspace{3mm}

\begin{num}
\item\label{eq 2.4.1} For all $f\in \mathcal{S}_{\omega,\scusp}(G(F))$ the function  $x\mapsto K^\chi_f(x)$ is compactly supported.
\end{num}

\vspace{3mm}

\begin{num}
\item\label{eq 2.4.2} Assume that the pair $(G,H)$ is tempered. Then, for all $d>0$ and all $f\in \mathcal{C}_{\omega,\scusp}(G(F))$, we have
$$\displaystyle \lv K^\chi_f(x)\rvert\ll \Xi^X(x)^2\sigma_X(x)^{-d}$$
for all $x\in X$.
\end{num}

\vspace{3mm}

\noindent We will only show \ref{eq 2.4.2}, the proof of \ref{eq 2.4.1} being similar and actually slightly easier. Moreover, the proof of \ref{eq 2.4.2} is also very similar to the proof of Theorem 8.1.1 (ii) of \cite{B1}. We will thus content ourself with outlining the main steps. Let $A_0$ be a maximal $(\theta,F)$-split subtorus of $G$, $M_0:=\Cent_G(A_0)$ and $P_0=M_0U_0\in \mathcal{P}^\theta(M_0)$. Let $\overline{P}_0=\theta(P_0)$ be the opposite parabolic subgroup and set
$$\displaystyle A_{\overline{P}_0}^+:=\{a\in A_0(F);\;\lvert \alpha(a)\rvert\geqslant 1\; \forall \alpha\in R(A_0,\overline{P}_0) \}$$
By the weak Cartan decomposition \ref{eq 2.3.1} as well as Proposition \ref{prop HCS for X}(\ref{prop HCS for X i}), it suffices to show \ref{eq 2.4.2} only for $x=a\in A_{\overline{P}_0}^+$. For all $Q\in \mathcal{F}^\theta(M_0)$ and $\delta>0$ set
$$\displaystyle A_{\overline{Q}}^+(\delta):=\{a\in A_0(F);\; \lv \alpha(a)\rvert\geqslant e^{\delta\overline{\sigma}(a)}\; \forall \alpha\in R(A_0,U_{\overline{Q}}) \}$$
where $\overline{Q}:=\theta(Q)$ and $U_{\overline{Q}}$ denotes the unipotent radical of $\overline{Q}$. Then, if $\delta$ is sufficiently small we have
$$\displaystyle A_{\overline{P}_0}^+=\bigcup_{Q\in \mathcal{F}^\theta(M_0)-\{G\}, P_0\subset Q} A_{\overline{Q}}^+(\delta)\cap A_{\overline{P}_0}^+$$
Thus, fixing $Q\in \mathcal{F}^\theta(M_0)-\{ G\}$ with $P_0\subset Q$ and $\delta>0$, it suffices to prove the estimate \ref{eq 2.4.2} only for $x=a\in A_{\overline{Q}}^+(\delta)\cap A_{\overline{P}_0}^+$. We fix such a $Q$ and such a $\delta$ henceforth. Let $U_Q$ be the unipotent radical of $Q$ and set $L:=Q\cap \overline{Q}$, $H_L:=H\cap L$ and $H^Q:=H_L\ltimes U_Q$. We define a unitary character $\chi^Q$ of $H^Q(F)$ by setting $\chi^Q(h_Lu_Q)=\chi(h_L)$ for all $h_L\in H_L(F)$ and $u_Q\in U_Q(F)$.  Then by Proposition \ref{prop HCS for X}(\ref{prop HCS for X vi}) , $H^Q$ is a unimodular algebraic group and the pair $(G,H^Q)$ is tempered. Thus, fixing a Haar measure $dh^Q$ on $H^Q(F)$, we can define
$$\displaystyle K_{f}^{\chi,Q}(x):=\int_{A_G^H(F)\backslash H^Q(F)} f(x^{-1}h^Qx)\chi^Q(h^Q)^{-1}dh^Q$$
for all $f\in \mathcal{C}_{\omega}(G(F))$ and all $x\in G(F)$. As $U_Q\subset H^Q\subset Q$, for any strongly cuspidal function $f\in \mathcal{C}_{\omega,scusp}(G(F))$, the function $K_{f}^{\chi,Q}$ vanishes identically. Therefore, it is sufficient to show the existence of $c>0$ such that for every $f\in \mathcal{C}_{\omega}(G(F))$ and $d>0$ we have

\begin{align}\label{eq 2.4.a}
\displaystyle \left\lvert K^{\chi}_f(a)-cK_{f}^{\chi,Q}(a)\right\rvert\ll \Xi^X(a)^2\sigma_X(a)^{-d}
\end{align}

\noindent for all $a\in A_{\overline{Q}}^+(\delta)\cap A_{\overline{P}_0}^+$. We prove this following closely the proof of Proposition 8.1.4 of \cite{B1}. We henceforth fix $f\in \mathcal{C}_{\omega}(G(F))$ and $d>0$.

\vspace{2mm}

\noindent The $F$-analytic map
$$H(F)\cap \overline{P}_0(F)U_0(F)\to U_0(F)$$
$$h=\overline{p}u\mapsto u$$
is submersive at the origin (this follows from the fact that $\overline{P}_0(F)H(F)$ is open in $G(F)$). Therefore, we may find a compact-open neighborhood $\mathcal{U}_0$ of $1$ in $U_0(F)$ together with an $F$-analytic map
$$h:\mathcal{U}_0\to H(F)$$
such that $h(u)\in \overline{P}_0(F)u$ for all $u\in \mathcal{U}_0$ and $h(1)=1$. Set $\mathcal{U}_{Q}:=\mathcal{U}_0\cap U_{Q}(F)$, $\mathcal{H}:=H_L(F)h(\mathcal{U}_{Q})$ and fix Haar measures $dh_L$ and $du_Q$ on $H_L(F)$ and $U_Q(F)$ whose product is the fixed Haar measure on $H^Q(F)$. The following fact is easy and can be proved exactly the same way as point (8.1.7) of \cite{B1} (note that here $H\cap \overline{Q}=H\cap L$):

\vspace{3mm}

\begin{num}
\item\label{eq 2.4.3} The map $H_L(F)\times \mathcal{U}_Q\to H(F)$, $(h_L,u_Q)\mapsto h_Lh(u_Q)$, is an $F$-analytic open embedding with image $\mathcal{H}$ and there exists a smooth function $j\in C^\infty(\mathcal{U}_Q)$ such that
$$\displaystyle \int_{\mathcal{H}}\varphi(h)dh=\int_{H_{L}(F)}\int_{\mathcal{U}_Q} \varphi(h_{L}h(u_Q)) j(u_Q)du_Qdh_{L}$$
for all $\varphi \in L^1(\mathcal{H})$.
\end{num}

\vspace{3mm}

\noindent Fix $\epsilon>0$ that we will assume sufficiently small in what follows. By Lemma \ref{lem estimates}(\ref{lem estimates ii}), for $\epsilon$ small enough we have
$$\displaystyle aU_Q\left[<\epsilon \overline{\sigma}(a)\right]a^{-1}\subseteq \mathcal{U}_Q$$
for all $a\in A_{\overline{Q}}^+(\delta)$. This allows us to define
$$\displaystyle H^{<\epsilon,a}:=H_{L}(F) h\left(a U_Q\left[<\epsilon \overline{\sigma}(a)\right]a^{-1}\right),$$
$$\displaystyle H^{Q,<\epsilon,a}:=H_{L}(F)a U_Q\left[<\epsilon \sigma(a)\right]a^{-1}$$
and the following expressions
$$\displaystyle K_f^{\chi,<\epsilon}(a):=\int_{A_G^H(F)\backslash H^{<\epsilon,a}} f(a^{-1}ha)\chi(h)^{-1}dh$$
$$\displaystyle K_f^{\chi,Q,<\epsilon}(a):=\int_{A_G^H(F)\backslash H^{Q,<\epsilon,a}} f(a^{-1}h^Qa)\chi^Q(h^Q)^{-1}dh^Q$$
for all $a\in A_{\overline{Q}}^+(\delta)$. Set $c=j(1)$ (where the function $j(.)$ is the one appearing in \ref{eq 2.4.3}). Obviously \ref{eq 2.4.a} will follows if we can show that for $\epsilon$ sufficiently small we have:

\begin{align}\label{eq 2.4.b}
\displaystyle \left\lvert K_f^\chi(a)-K_f^{\chi,<\epsilon}(a)\right\rvert\ll \Xi^{X}(a)^2\sigma_{X}(a)^{-d}
\end{align}

\begin{align}\label{eq 2.4.c}
\displaystyle \left\lvert K_f^Q(a)-K_f^{\chi,Q,<\epsilon}(a)\right\rvert\ll \Xi^{X}(a)^2\sigma_{X}(a)^{-d}
\end{align}

\begin{align}\label{eq 2.4.d}
\displaystyle K_f^{\chi,<\epsilon}(a)=cK_f^{\chi,Q,<\epsilon}(a)
\end{align}

\noindent for all $a\in A_{\overline{Q}}^+(\delta) \cap A_{\overline{P}_0}^+$.

\vspace{3mm}

\noindent First we prove \ref{eq 2.4.b} and \ref{eq 2.4.c}. By Lemma \ref{lem estimates} (\ref{lem estimates i})-(\ref{lem estimates iii})-(\ref{lem estimates iv}) and noticing that 
$$\displaystyle H^{<\epsilon,a}=H(F)\cap \overline{Q}(F) a U_Q\left[<\epsilon \overline{\sigma}(a)\right]a^{-1} \mbox{ and } H^{Q,<\epsilon,a}=H^Q(F)\cap \overline{Q}(F) a U_Q\left[<\epsilon \overline{\sigma}(a)\right]a^{-1},$$
we see that $\overline{\sigma}(a)\ll \overline{\sigma}(a^{-1}ha)$ and $\overline{\sigma}(a)\ll \overline{\sigma}(a^{-1}h^Qa)$ for all $a\in A_{\overline{Q}}^+(\delta)$, $h\in H(F)\setminus H^{<\epsilon,a}$ and $h^Q\in H^Q(F)\setminus H^{Q,<\epsilon,a}$. Thus, by definition of $\mathcal{C}_{\omega}(G(F))$, for every $d_1>0$ and $d_2>0$, the left hand sides of \ref{eq 2.4.b} and \ref{eq 2.4.c} are essentially bounded by
$$\displaystyle \overline{\sigma}(a)^{-d_1}\int_{A_G^H(F)\backslash H(F)} \Xi^G(a^{-1}ha)\overline{\sigma}(a^{-1}ha)^{-d_2} dh$$
and
$$\displaystyle \overline{\sigma}(a)^{-d_1}\int_{A_G^H(F)\backslash H^Q(F)} \Xi^G(a^{-1}h^Qa)\overline{\sigma}(a^{-1}h^Qa)^{-d_2} dh^Q$$
\noindent for all $a\in A_{\overline{Q}}^+(\delta)$ respectively. By Proposition \ref{prop HCS for X} (\ref{prop HCS for X v})-(\ref{prop HCS for X vi}), there exists $d_3>0$ such that for $d_2$ sufficiently large the last two expressions above are essentially bounded by $\Xi^X(a)^2 \sigma_X(a)^{d_3}\overline{\sigma}(a)^{-d_1}$ for all $a\in A_{\overline{Q}}^+$. Finally, by Proposition \ref{prop HCS for X} (\ref{prop HCS for X ii}), if $d_1$ is sufficiently large we have $\sigma_X(a)^{d_3}\overline{\sigma}(a)^{-d_1}\ll \sigma_X(a)^{-d}$ for all $a\in A_0(F)$. This proves \ref{eq 2.4.b} and \ref{eq 2.4.c}.

\vspace{2mm}

\noindent It only remains to prove \ref{eq 2.4.d}. By \ref{eq 2.4.3} and the choice of Haar measures, we have

\begin{align}\label{eq 2.4.7}
\displaystyle K_f^{\chi,<\epsilon}(a)=\int_{A^H_G(F)\backslash H_L(F)} \int_{aU_Q[<\epsilon \overline{\sigma}(a)]a^{-1}} f(a^{-1}h_Lh(u_Q)a)\chi(h_Lh(u_Q))^{-1} j(u_Q) du_Qdh_L
\end{align}

\noindent and

\begin{align}\label{eq 2.4.7bis}
\displaystyle K_f^{\chi,Q,<\epsilon}(a)=\int_{A^H_G(F)\backslash H_L(F)} \int_{aU_Q[<\epsilon \overline{\sigma}(a)]a^{-1}} f(a^{-1}h_Lu_Qa) \chi(h_L)^{-1} du_Q dh_L
\end{align}

\noindent for all $a\in A_{\overline{Q}}^+(\delta)$. Since the function $u_Q\mapsto \chi(h(u_Q))^{-1} j(u_Q)$ is smooth, by Lemma \ref{lem estimates}(\ref{lem estimates ii}) for $\epsilon$ sufficiently small we have

\begin{align}\label{eq 2.4.8}
\displaystyle \chi(h(u_Q))^{-1} j(u_Q)=j(1)=c
\end{align}

\noindent for all $a\in A_{\overline{Q}}^+(\delta)$ and all $u_Q\in aU_Q[<\epsilon \overline{\sigma}(a)]a^{-1}$. Let $J\subset G(F)$ be a compact-open subgroup by which $f$ is right invariant. By definition, the map $u_Q\mapsto h(u_Q)u_Q^{-1}$ is $F$-analytic, sends $1$ to $1$ and takes values in $\overline{P}_0(F)$. By Lemma \ref{lem estimates}(\ref{lem estimates ii}) again, it follows that for all $0<\delta'<\delta$ and all $c_0>0$ if $\epsilon$ is sufficiently small we have
$$\displaystyle a^{-1}h(u_Q)u_Q^{-1}a\in \exp\left(B(0,c_0e^{-\delta'\overline{\sigma}(a)})\cap\overline{\mathfrak{p}}_0(F) \right)$$
for all $a\in A_{\overline{Q}}^+(\delta)\cap A_{\overline{P}_0}^+$ and all $u_Q\in aU_Q[<\epsilon \overline{\sigma}(a)]a^{-1}$. Moreover, there exists $\alpha>0$ such that $\lvert \Ad(g^{-1})X\rvert_{\mathfrak{g}}\leqslant e^{\alpha \overline{\sigma}(g)}\lvert X\rvert_{\mathfrak{g}}$ for all $g\in G(F)$ and all $X\in \mathfrak{g}(F)$. Hence, if $\epsilon$ is sufficiently small we have

\begin{align}\label{eq 2.4.9}
\displaystyle a^{-1}u_Q^{-1}h(u_Q)a=(a^{-1}u_Qa)^{-1}(a^{-1}h(u_Q)u_Q^{-1}a)(a^{-1}u_Qa)\in J
\end{align}

\noindent for all $a\in A_{\overline{Q}}^+(\delta)\cap A_{\overline{P}_0}^+$ and all $u_Q\in aU_Q[<\epsilon \overline{\sigma}(a)]a^{-1}$. It is clear that \ref{eq 2.4.d} follows from \ref{eq 2.4.7}, \ref{eq 2.4.7bis}, \ref{eq 2.4.8} and \ref{eq 2.4.9} and this ends the proof of Theorem \ref{theo conv}. $\blacksquare$

\section{The spectral side}

\subsection{The statement}\label{section statement spectral side}

In this chapter $G$ is a connected reductive group over $F$, $\overline{G}:=G/A_G$ and $H$ is a symmetric subgroup of $G$. As in the previous chapter, we let $A^H_G$ as the connected component of $H\cap A_G$. Let $\omega$ and $\chi$ be continuous unitary characters of $A_G(F)$ and $H(F)$ respectively such that $\omega_{\mid A^H_G(F)}=\chi_{\mid A^H_G(F)}$. Set

$$\displaystyle \nu(H):=\left[H(F)\cap A_G(F):A_H(F) \right]$$

In \S \ref{statement convergence}, we have defined a linear form

$$\displaystyle f\in \mathcal{S}_{\omega,\scusp}(G(F))\mapsto J^\chi(f)$$  

\noindent which, if $(G,H)$ is a tempered pair, extends to a continuous linear form

$$\displaystyle f\in \mathcal{C}_{\omega,\scusp}(G(F))\mapsto J^\chi(f)$$

For all $\pi\in \Irr(G)$ we define a multiplicity $m(\pi,\chi)$ by

$$\displaystyle m(\pi,\chi):=\dim \Hom_H(\pi,\chi)$$

\noindent where $\Hom_H(\pi,\chi)$ denotes the space of linear forms $\ell:\pi\to\mathbf{C}$ such that $\ell\circ\pi(h)=\chi(h)\ell$ for all $h\in H(F)$. By Theorem 4.5 of \cite{Del}, we know that this space is always finite dimensional so that the multiplicity $m(\pi,\chi)$ is well-defined.

\vspace{2mm}

Recall that $\Irr_{\omega,\cusp}(G)$, resp. $\Irr_{\omega,\sqr}(G)$, denote the sets of equivalence classes of irreducible supercuspidal, resp. square-integrable, representations of $G(F)$ with central character $\omega$. Define the following linear forms

$$\displaystyle f\in \mathcal{S}_{\omega,\scusp}(G(F))\mapsto J_{spec,\cusp}^\chi(f):=\nu(H)\sum_{\pi\in \Irr_{\omega,\cusp}(G)}m(\pi,\chi) \Tr(\pi^\vee(f))$$

\noindent and

$$\displaystyle f\in \mathcal{C}_{\omega,\scusp}(G(F))\mapsto J_{spec,disc}^\chi(f):=\nu(H)\sum_{\pi\in \Irr_{\omega,\sqr}(G)}m(\pi,\chi) \Tr(\pi^\vee(f))$$

\noindent Notice that the sums defining these linear forms are always finite by the result of Harish-Chandra that for every compact-open subgroup $J\subset G(F)$ the set of $\pi\in \Irr_{\omega,\sqr}(G)$ with $\pi^J\neq 0$ is finite (\cite{Wal1} Th\'eor\`eme VIII.1.2). Recall that in \S \ref{strongly cuspidal} we have introduced certain spaces ${}^0\mathcal{S}_\omega(G(F))$, ${}^0\mathcal{C}_{\omega}(G(F))$ of cusp forms. The goal of this chapter is to prove the following

\begin{theo}\label{theo spec}
For all $f\in {}^0\mathcal{S}_{\omega}(G(F))$ we have

$$\displaystyle J^\chi(f)=J^\chi_{spec,\cusp}(f)$$

\noindent Moreover, if $(G,H)$ is a tempered pair, for all $f\in {}^0\mathcal{C}_{\omega}(G(F))$ we have

$$\displaystyle J^\chi(f)=J^\chi_{spec,disc}(f)$$
\end{theo}

\subsection{Explicit description of the intertwinings}

For all $\pi\in \Irr_{\omega,\cusp}(G(F))$ we define a bilinear form

$$\displaystyle \mathcal{B}_\pi:\pi \times \pi^\vee\to \mathbf{C}$$

\noindent by

$$\displaystyle \mathcal{B}_\pi(v,v^\vee):=\int_{A^H_G(F)\backslash H(F)} \langle \pi(h)v,v^\vee\rangle \chi(h)^{-1} dh$$

\noindent for all $(v,v^\vee)\in \pi\times\pi^\vee$. If $(G,H)$ is a tempered pair and $\pi\in \Irr_{\omega,\sqr}(G(F))$ then the above integral is also absolutely convergent and thus also defines a bilinear form $\mathcal{B}_\pi:\pi\times \pi^\vee\to \mathbf{C}$. In all cases, we have

$$\displaystyle \mathcal{B}_\pi(\pi(h_1)v,\pi^\vee(h_2)v^\vee)=\chi(h_1)\chi^{-1}(h_2)\mathcal{B}_\pi(v,v^\vee)$$

\noindent for all $(v,v^\vee)\in \pi\times \pi^\vee$ and all $h_1,h_2\in H(F)$. Thus $\mathcal{B}_\pi$ factorizes through a bilinear form

$$\displaystyle \mathcal{B}_\pi: \pi_\chi\times \pi^\vee_{\chi^{-1}}\to\mathbf{C}$$

\noindent where $\pi_\chi$ and $\pi^\vee_{\chi^{-1}}$ denote the spaces of $(H(F),\chi)$- and $(H(F),\chi^{-1})$-coinvariants in $\pi$ and $\pi^\vee$ respectively i.e. the quotients of $\pi$ and $\pi^\vee$ by the subspaces generated by vectors of the form $\pi(h)v-\chi(h)v$ ($h\in H(F)$, $v\in \pi$) and $\pi^\vee(h)v^\vee-\chi(h)^{-1}v^\vee$ ($h\in H(F)$, $v^\vee\in \pi^\vee$) respectively. The following proposition has been proved in more generality in \cite{SV} Theorem 6.4.1 when the subgroup $H$ is {\it strongly tempered} (in the sense of {\it loc. cit.}). The same kind of idea already appears in \cite{Wal3} Proposition 5.6.

\begin{prop}\label{prop intertwinings}
$\mathcal{B}_\pi$ induces a nondegenerate pairing between $\pi_\chi$ and $\pi_{\chi^{-1}}$.
\end{prop}

\vspace{2mm}

\noindent\ul{Proof}: We will prove the proposition when $(G,H)$ is a tempered pair and $\pi\in \Irr_{\omega,\sqr}(G)$, the case where $\pi\in \Irr_{\omega,\cusp}(G)$ being similar and easier. Fix a $G(F)$-invariant scalar product $(.,.)$ on $\pi$. We can define the following sesquilinear version of $\mathcal{B}_\pi$

$$\displaystyle \mathcal{L}_\pi: \pi\times \pi\to \mathbf{C}$$

$$\displaystyle (v,v')\mapsto \int_{A_G^H(F)\backslash H(F)}  (\pi(h)v,v') \chi(h)^{-1} dh$$

\noindent which factorizes through a sesquilinear pairing $\mathcal{L}_\pi:\pi_\chi\times\pi_\chi\to \mathbf{C}$. Obviously, it suffices to show that this pairing is non degenerate. Since $\pi_\chi$ is finite dimensional, this is equivalent to saying that the map

$$\displaystyle v\in \pi\mapsto \mathcal{L}_\pi(.,v)\in \Hom_H(\pi,\chi)$$

\noindent is surjective. To continue we need the following lemma, a consequence of the weak Cartan decomposition (\ref{eq 2.3.1}):

\begin{lem}\label{lemma intertwinings}
For all $\ell\in \Hom_H(\pi,\chi)$ and all $v\in \pi$ we have

$$\displaystyle \int_X \lvert \ell(\pi(x)v)\rvert^2 dx<\infty$$

\noindent and moreover for all $f\in \mathcal{C}_{\omega^{-1}}(G(F))$ the integral

$$\displaystyle \int_{A_G(F)\backslash G(F)} f(g) \ell(\pi(g)v) dg$$

\noindent is absolutely convergent and equals $\ell(\pi(f)v)$.
\end{lem}

\noindent\ul{Proof}: For every compact-open subgroup $J\subset G(F)$ we will denote by $e_J\ast \ell$ the smooth linear form (i.e. an element of $\pi^\vee$) $v\in \pi\mapsto \frac{1}{\vol(J)}\int_J\ell(\pi(k)v)dk$. Let $A_0$ be a maximal $(\theta,F)$-split subtorus of $G$, $M_0:=\Cent_G(A_0)$ (a minimal $\theta$-split Levi subgroup) and $P_0\in \mathcal{P}^\theta(M_0)$. Set

$$\displaystyle A_{P_0}^+:=\{a\in A_0(F); \lvert \alpha(a)\rvert\geqslant 1 \; \forall \alpha\in \Delta_{P_0} \}$$

\noindent Then, by the weak Cartan decomposition \ref{eq 2.3.1} and Proposition \ref{prop HCS for X}(\ref{prop HCS for X i})-(\ref{prop HCS for X ii})-(\ref{prop HCS for X iii})-(\ref{prop HCS for X iv}), in order to prove that the two integrals of the proposition are convergent it suffices to show that for all $d>0$ we have

\begin{align}\label{eq 3.2.1}
\displaystyle \lvert\ell(\pi(a)v)\rvert\ll \Xi^G(a)\overline{\sigma}(a)^{-d}
\end{align}

\noindent for all $a\in A_{P_0}^+$. Let $J\subset G(F)$ be a compact-open subgroup such that $v\in \pi^J$. Since $P_0$ is $\theta$-split, $H(F)P_0(F)$ is open in $G(F)$ and consequently there exists a compact-open subgroup $J'\subset G(F)$ such that

$$\displaystyle J'\subset H(F)a\left(J\cap P_0(F)\right) a^{-1}$$

\noindent for all $a\in A_{P_0}^+$. Thus, for all $a\in A_{P_0}^+$ we have $\ell(\pi(a)v)=\langle e_{J'}\ast \ell, \pi(a)v\rangle$ and the inequality \ref{eq 3.2.1} now follows from the known asymptotics of smooth coefficients of square-integrable representations.

\vspace{2mm}

\noindent To prove the last part of the proposition, choose $J\subset G(F)$ a compact-open subgroup by which $f$ is invariant on the left. Then, we have

\[\begin{aligned}
\displaystyle \int_{A_G(F)\backslash G(F)} f(g) \ell(\pi(g)v) dg & =\frac{1}{\vol(J)} \int_{A_G(F)\backslash G(F)} f(g) \int_J\ell(\pi(kg)v)dk dg \\
 & =\int_{A_G(F)\backslash G(F)} f(g) \langle e_J\ast\ell,\pi(g)v\rangle dg \\
 & =\langle e_J\ast\ell,\pi(f)v\rangle=\ell(\pi(f)v)\; \blacksquare
\end{aligned}\]

\vspace{2mm}

By the lemma we can define a scalar product, also denoted $(.,.)$, on $\Hom_H(\pi,\chi)$ characterized by

$$\displaystyle \int_X \ell(\pi(x)v) \overline{\ell'(\pi(x)v')}dx=(\ell,\ell')(v,v')$$

\noindent for all $\ell,\ell'\in \Hom_H(\pi,\chi)$ and all $v,v'\in \pi$. Let $\ell\in \Hom_H(\pi,\chi)$ which is orthogonal for this scalar product to all the forms $\mathcal{L}_\pi(.,v)$ for $v\in \pi$. To conclude it suffices to show that this implies $\ell=0$. Since for all $v,v'\in \pi$ we have $(v,\pi(.)v')\in\mathcal{C}_{\omega^{-1}}(G(F))$, by the lemma we have

\[\begin{aligned}
\displaystyle 0=\int_X \ell(\pi(x)v')\overline{\mathcal{L}_\pi(\pi(x)v',v)} dx=\nu(H)\int_{A_G(F)\backslash G(F)} \ell(\pi(g)v') (v,\pi(g)v')dg=\frac{\nu(H)}{d(\pi)} \ell(v) (v',v')
\end{aligned}\]

\noindent for all $v,v'\in \pi$ and where $d(\pi)$ stands for the formal degree of $\pi$. Hence, $\ell=0$. $\blacksquare$

\subsection{Proof of Theorem \ref{theo spec}}

Once again we will prove the theorem in the case where $(G,H)$ is a tempered pair and $f\in {}^0\mathcal{C}_{\omega}(G(F))$. The case where $f\in {}^0\mathcal{S}_{\omega}(G(F))$ is completely similar since compactly supported cusp forms are linear combinations of matrix coefficients of supercuspidal representations whereas a general cusp form $f\in {}^0\mathcal{C}_{\omega}(G(F))$ is a linear combination of matrix coefficients of square-integrable representations.

More precisely, let $f\in {}^0\mathcal{C}_{\omega}(G(F))$ and for all $\pi\in \Irr_{\omega,\sqr}(G(F))$, set $f_\pi(g):=\Tr(\pi^\vee(g^{-1})\pi^\vee(f))$ for all $g\in G(F)$. Then, we have $f_\pi\in {}^0\mathcal{C}_{\omega}(G(F))$ for all $\pi\in \Irr_{\omega,\sqr}(G(F))$ and by the Harish-Chandra-Plancherel formula for cusp forms \ref{eq 1.6.1} we have

$$\displaystyle J^\chi(f)=\sum_{\pi\in \Irr_{\omega,\sqr}(G(F))}d(\pi)J^\chi(f_\pi)$$

\noindent Thus, it suffices to show that 

$$J^\chi(f_\pi)=\nu(H)d(\pi)^{-1}m(\pi,\chi)\Tr(\pi^\vee(f))$$

\noindent for all $\pi\in \Irr_{\omega,\sqr}(G(F))$. Fix $\pi\in \Irr_{\omega,\sqr}(G(F))$. As $f_\pi$ is a sum of coefficients of $\pi$, the equality above is equivalent to

$$\displaystyle J^\chi(f_{v,v^\vee})=\nu(H)d(\pi)^{-1}m(\pi,\chi)f_{v,v^\vee}(1)$$

\noindent for all $(v,v^\vee)\in \pi\times \pi^\vee$ where $f_{v,v^\vee}(g):=\langle \pi(g)v,v^\vee\rangle$ for all $g\in G(F)$. Fix $(v,v^\vee)\in \pi\times \pi^\vee$. Then, we have

$$\displaystyle K_{f_{v,v^\vee}}^\chi(x)=\mathcal{B}_\pi(\pi(x)v,\pi^\vee(x)v^\vee)$$

\noindent for all $x\in X$. Choose a basis $(\overline{v}_1,\ldots,\overline{v}_N)$ of $\pi_\chi$ (where $N=m(\pi,\chi)$) and let $(\overline{v}^\vee_1,\ldots,\overline{v}^\vee_N)$ be the dual basis of $\pi^\vee_{\chi^{-1}}$ with respect to $\mathcal{B}_\pi$ (such a dual basis exists thanks to Proposition \ref{prop intertwinings}). Let $(v_1,\ldots,v_N)$ and $(v^\vee_1,\ldots,v^\vee_N)$ be any lifting of these basis to $\pi$ and $\pi^\vee$ respectively. Then we have

$$\displaystyle \mathcal{B}_\pi(\pi(x)v,\pi^\vee(x)v^\vee)=\sum_{i=1}^N \mathcal{B}_\pi(\pi(x)v,v_i^\vee)\mathcal{B}_\pi(v_i,\pi^\vee(x)v^\vee)$$

\noindent for all $x\in X$. Now by Lemma \ref{lemma intertwinings}, we have

\[\begin{aligned}
\displaystyle J^\chi(f_{v,v^\vee}) & =\int_X \mathcal{B}_\pi(\pi(x)v,\pi^\vee(x)v^\vee) dx=\sum_{i=1}^N \int_X \mathcal{B}_\pi(\pi(x)v,v_i^\vee)\mathcal{B}_\pi(v_i,\pi^\vee(x)v^\vee) dx \\
 & =\nu(H)\sum_{i=1}^N \int_{A_G(F)\backslash G(F)}\langle \pi(g)v,v_i^\vee\rangle \mathcal{B}_\pi(v_i,\pi^\vee(g)v^\vee) dg \\
 & =\nu(H)\sum_{i=1}^N \frac{\langle v,v^\vee\rangle}{d(\pi)} \mathcal{B}_\pi(v_i,v_i^\vee) \\
 & =\nu(H) N \frac{\langle v,v^\vee\rangle}{d(\pi)}=\nu(H)d(\pi)^{-1}m(\pi,\chi) f_{v,v^\vee}(1) \; \blacksquare
\end{aligned}\]

\section{The geometric side}

\subsection{The statement}\label{statement geom side}

In this chapter $E/F$ is a quadratic extension, $H$ is a connected reductive group over $F$ and $G:=R_{E/F}H_E$. We have a natural inclusion $H\hookrightarrow G$ and we shall denote by $\theta$ involution of $G$ induced by the nontrivial element of $\Gal(E/F)$. Hence $H=G^\theta$. Note that in this case, with the notations of \S \ref{statement convergence} we have $A^H_G=A_H$. Set

$$\displaystyle \nu(H):=\left[H(F)\cap A_G(F):A_H(F) \right]$$

\noindent As in \S \ref{statement convergence}, we let $\overline{G}:=G/A_G$, $\overline{H}:=H/A_H$, $X:=A_G(F)H(F)\backslash G(F)$, $\mathbf{X}:=HA_G\backslash G$, $\sigma_X:=\sigma_{\mathbf{X}}$ and $\overline{\sigma}:=\sigma_{\overline{G}}$. Note that $X$ is an open subset of $\mathbf{X}(F)$. We have the following identity between Weyl discriminants

\begin{align}\label{eq 4.1.1}
\displaystyle D^H(h)=D^G(h)^{1/2},\;\;\;h\in H_{\reg}(F)
\end{align}

\noindent which will be crucial in what follows.

\vspace{2mm}

\noindent Let $\omega$ and $\chi$ be continuous unitary characters of $A_G(F)$ and $H(F)$ respectively such that $\omega_{\mid A_H(F)}=\chi_{\mid A_H(F)}$. In \S \ref{statement convergence}, we have defined a continuous linear form $f\in \mathcal{C}_{\omega, \scusp}(G(F))\mapsto J^\chi(f)$. We define a second continuous linear form $f\in \mathcal{C}_{\omega, \scusp}(G(F))\mapsto J_{\geom}^\chi(f)$ by setting

$$\displaystyle J^\chi_{\geom}(f):=\nu(H)\sum_{T\in \mathcal{T}_{\elli}(H)} \lv W(H,T)\rv^{-1} \int_{\overline{T}(F)} D^H(t) \Theta_f(t)\chi(t)^{-1}dt$$

\noindent for all $f\in \mathcal{C}_{\omega, \scusp}(G(F))$, where $\mathcal{T}_{\elli}(H)$ denotes a set of representatives of the $H(F)$-conjugacy classes of maximal elliptic tori in $H$, we have set $\overline{T}:=T/A_H$ for all $T\in \mathcal{T}_{\elli}(H)$ and we recall that $\overline{T}(F)$ is equipped with the Haar measure of total mass $1$. Since for all $f\in \mathcal{C}_{\omega,\scusp}(G(F))$ the function $(D^G)^{1/2}\Theta_f$ is locally bounded, by \ref{eq 4.1.1} we see that the expression defining $J^\chi_{\geom}(f)$ is absolutely convergent. The goal of this chapter is to show the following

\begin{theo}\label{theo geom side}
For all $f\in \mathcal{C}_{\omega, \scusp}(G(F))$, we have

$$\displaystyle J^\chi(f)=J^\chi_{\geom}(f)$$
\end{theo}

\vspace{2mm}

\noindent We fix a function $f\in \mathcal{C}_{\omega,\scusp}(G(F))$ until the end of this chapter.
 
\subsection{Truncation and first decomposition}\label{truncation}

\noindent We fix a sequence $(\kappa_N)_{N\geqslant 1}$ of functions $\kappa_N:X(F)\to \{0,1\}$ satisfying the two following conditions:

\vspace{3mm}

\begin{num}
\item\label{eq 4.2.1} There exist $C_1,C_2>0$ such that for all $x\in X(F)$ and all $N\geqslant 1$, we have:

$$\sigma_{X}(x)\leqslant C_1N\Rightarrow \kappa_N(x)=1$$

$$\kappa_N(x)\neq 0\Rightarrow \sigma_{X}(x)\leqslant C_2N$$
\end{num}

\vspace{3mm}

\begin{num}
\item\label{eq 4.2.2} There exists an open-compact subgroup $K'\subset G(F)$ such that the function $\kappa_N$ is right-invariant by $K'$ for all $N\geqslant 1$.
\end{num}

\vspace{3mm}

\noindent Such a sequence of truncation functions is easy to construct (see \cite{B1} \S 10.9). Set

$$\displaystyle J^\chi_N(f)=\int_{X} K^\chi(f,x) \kappa_N(x) dx$$

\noindent for all $N\geqslant 1$. Then we have

\begin{align}\label{eq 4.2.3}
\displaystyle J^\chi(f)=\lim\limits_{N\to \infty} J^\chi_N(f)
\end{align}

\noindent Let $\mathcal{T}(H)$ be a set of representatives of the $H(F)$-conjugacy classes of maximal tori in $H$. By the Weyl integration formula for $H$, we have

$$\displaystyle K^\chi(f,x)=\sum_{T\in \mathcal{T}(H)} \lv W(H,T)\rvert^{-1}\int_{\overline{T}(F)} D^H(t)\int_{T(F)\backslash H(F)} f(x^{-1}h^{-1}thx)dh\chi(t)^{-1}dt$$

\noindent where we have set $\overline{T}:=T/A_H$ for all $T\in \mathcal{T}(H)$. At least formally, it follows that for all $N\geqslant 1$

\begin{align}\label{eq 4.2.4}
\displaystyle J^\chi_N(f)=\sum_{T\in \mathcal{T}(H)} \lv W(H,T)\rvert^{-1}\int_{\overline{T}(F)} D^H(t)\int_{T^G(F)\backslash G(F)}f(g^{-1}tg)\kappa_{N,T}(g)dg \chi(t)^{-1}dt
\end{align}

\noindent where for all $T\in \mathcal{T}(H)$ we have denoted by $T^G$ the centralizer of $T$ in $G$ (a maximal torus in $G$) and we have set

$$\displaystyle \kappa_{N,T}(g):=\int_{A_G(F)T(F)\backslash T^G(F)} \kappa_N(ag)da$$

\noindent for all $g\in G(F)$. Define

$$\displaystyle J^\chi_{N,T}(f):=\int_{\overline{T}(F)} D^H(t)\int_{T^G(F)\backslash G(F)}f(g^{-1}tg)\kappa_{N,T}(g)dg \chi(t)^{-1}dt$$

\noindent for all $N\geqslant 1$ and all $T\in \mathcal{T}(H)$. The equality \ref{eq 4.2.4} can thus be restated as

\begin{align}\label{eq 4.2.5}
\displaystyle J^\chi_N(f)=\sum_{T\in \mathcal{T}(H)} \lv W(H,T)\rvert^{-1}J^\chi_{N,T}(f)
\end{align}

\noindent The previous formal manipulations are justified a posteriori by the following lemma:

\begin{lem}\label{lemma kappaN}
\begin{enumerate}[(i)]
\item\label{lemma kappaN i} There exists $k\geqslant 1$ such that

$$\displaystyle \kappa_{N,T}(g)\ll N^k\sigma_{T^G\backslash G}(g)^k$$

\noindent for all $N\geqslant 1$ and all $g\in G(F)$.

\item\label{lemma kappaN ii} For all $T\in \mathcal{T}(H)$, the expression defining $J^\chi_{N,T}(f)$ is absolutely convergent and the identity \ref{eq 4.2.5} is valid.

\end{enumerate}
\end{lem}

\vspace{2mm}

\noindent\ul{Proof}:

\begin{enumerate}[(i)]
\item Since the natural inclusion $A_GT\backslash T^G\subset \mathbf{X}$ is a closed immersion (essentially because $T^G$ is $\theta$-stable), we have $\sigma_X(a)\sim \sigma_{
A_GT\backslash T^G}(a)$ for all $a\in T^G(F)$. As $\sigma_X(x)\ll \sigma_X(xg)\sigma(g)$ for all $(x,g)\in X\times G(F)$, it follows from \ref{eq 4.2.1} that there exists $c_1>0$ such that for all $N\geqslant 1$, all $a\in T^G(F)$ and all $g\in G(F)$ we have

$$\displaystyle \kappa_N(ag)\neq 0\Rightarrow \sigma_{
A_GT\backslash T^G}(a)< c_1N\sigma(g)$$

\noindent Hence, since the function $\kappa_N$ is nonnegative and bounded by $1$, by \ref{eq 1.2.3} there exists $k>0$ such that

$$\displaystyle \kappa_{N,T}(g)\leqslant meas\left((A_GT\backslash T^G)[< c_1 N\sigma(g)]\right)\ll N^k \sigma(g)^k$$

\noindent for all $N\geqslant 1$ and all $g\in G(F)$. The function $g\mapsto \kappa_{N,T}(g)$ being left invariant by $T^G(F)$ we may replace $\sigma(g)$ in the inequality above by $\inf_{a\in T^G(F)}\sigma(ag)$ which by \ref{eq 1.2.1} is equivalent to $\sigma_{T^G\backslash G}(g)$. This proves (\ref{lemma kappaN i}).

\item Since $f$ belongs to the Harish-Chandra-Schwartz space $\mathcal{C}_{\omega}(G(F))$, this follows from a combination of (\ref{lemma kappaN i}), \ref{eq 4.1.1} and Lemma \ref{lemma conv abs woi}. $\blacksquare$
\end{enumerate}

\noindent From now on and until the end of \S \ref{proof geo side}, we fix a torus $T\in \mathcal{T}(H)$.

\subsection{Change of truncation}

Set $\overline{T}^G:=T^G/A_G$ (a maximal torus of $\overline{G}$) and let $\overline{A}$ be the maximal $(\theta,F)$-split subtorus of $\overline{T}^G$. Let $A\subset T^G$ be the inverse image of $\overline{A}$ and set

$$\displaystyle \kappa_{N,A}(g):=\int_{A_G(F)\backslash A(F)} \kappa_N(ag) da$$

\noindent for all $g\in G(F)$ and all $N\geqslant 1$. We define the following quantity

$$\displaystyle \nu(T):=[H(F)\cap A_G(F):A_H(F)] \times [A(F)\cap A_T(F):A_H(F)]^{-1}$$

\noindent Then, we have

\begin{align}\label{eq 4.3.1}
\displaystyle  J^\chi_{N,T}(f)=\nu(T)\int_{\overline{T}(F)} D^H(t)\int_{A_T(F)A(F)\backslash G(F)} f(g^{-1}tg) \kappa_{N,A}(g)dg \chi(t)^{-1}dt
\end{align}

\noindent Indeed, by our choices of Haar measures on tori (see \S \ref{Groups, measures, notations}) and noting that $T(F)\cap A_G(F)=H(F)\cap A_G(F)$, we have

\[\begin{aligned}
\displaystyle \kappa_{N,T}(g) & =\int_{A_G(F)T(F)\backslash T^G(F)} \kappa_N(ag) da = [(T\cap A_G)(F):(A_T\cap A_G)(F)] \int_{A_G(F)A_T(F)\backslash T^G(F)} \kappa_N(ag) da \\
 & =[(H\cap A_G)(F): (A_T\cap A_G)(F)] \int_{A_T(F)A(F)\backslash T^G(F)} \int_{A_G(F) (A\cap A_T)(F)\backslash A(F)} \kappa_N(at^Gg)da dt^G \\
 & =[(H\cap A_G)(F): (A_G\cap A_T)(F)] [(A\cap A_T)(F):(A_G\cap A_T)(F)]^{-1}\int_{A_T(F)A(F)\backslash T^G(F)}  \kappa_{N,A}(t^Gg) dt^G \\
 & =\nu(T) \int_{A_T(F)A(F)\backslash T^G(F)} \kappa_{N,A}(t^Gg) dt^G
\end{aligned}\]

\noindent for all $g\in G(F)$ and all $N\geqslant 1$, hence the result.

\vspace{2mm}

\noindent Since $A(F)A_T(F)\backslash T^G(F)$ is compact by \ref{eq 1.2.1} we have

\begin{align}\label{eq 4.3.2}
\displaystyle \sigma_{T^G\backslash G}(g)\sim \inf_{a\in A(F)A_T(F)}\sigma(ag)
\end{align}

\noindent for all $g\in G(F)$ and hence the same proof as that of Lemma \ref{lemma kappaN}(\ref{lemma kappaN i}) shows that there exists $k>0$ such that

\begin{align}\label{eq 4.3.3}
\displaystyle \kappa_{N,A}(g)\ll N^k \sigma_{T^G\backslash G}(g)^k
\end{align}

\noindent for all $N\geqslant 1$ and all $g\in G(F)$.

\vspace{2mm}

\noindent Let $M$ be the centralizer in $G$ of $A$. It is a $\theta$-split Levi subgroup with $A_{M,\theta}=A_\theta$. Indeed, the inclusion $A_\theta\subset A_{M,\theta}$ is obvious and $T^G$ is a maximal torus of $M$ hence $A_{M,\theta}$ is included in $A_{T^G,\theta}=A_{\theta}$. Let $\overline{A}_0$ be a maximal $(\theta,F)$-split subtorus of $\overline{G}$ containing $\overline{A}$ and denote by $A_0$ its inverse image in $G$. Let $M_0$ be the centralizer in $G$ of $A_0$. It is a minimal $\theta$-split Levi subgroup, we again have $A_{0,\theta}=A_{M_0,\theta}$ and we set $\mathcal{A}_{0,\theta}:=\mathcal{A}_{M_0,\theta}$. Let $K$ be a special maximal compact subgroup of $G(F)$. We use $K$ to define the functions $H_{Q,\theta}$ for all $Q\in \mathcal{F}^\theta(M_0)$ (see \S \ref{symmetric subgroups}). Fix $P_0\in \mathcal{P}^\theta(M_0)$ and let $\Delta_0$ be the set of simple roots of $A_0$ in $P_0$. To every $Y\in \mathcal{A}^+_{P_0,\theta}$ we associate a positive $(G,M_0,\theta)$-orthogonal set $(Y_{P_0'})_{P_0'\in \mathcal{P}^\theta(M_0)}$ by setting $Y_{P_0'}=wY$ where $w$ is the unique element in the little Weyl group $W(G,A_0)$ such that $wP_0=P_0'$. By the general constructions of \S \ref{GMT families}, this also induces a positive $(G,M,\theta)$-orthogonal set $(Y_P)_{P\in \mathcal{P}^\theta(M)}$. For all $g\in G(F)$, we define another $(G,M,\theta)$-orthogonal set $\mathcal{Y}(g)=(\mathcal{Y}(g)_P)_{P\in \mathcal{P}^\theta(M)}$ by setting

$$\displaystyle \mathcal{Y}(g)_P:= Y_P-H_{\overline{P},\theta}(g)$$

\noindent for all $P\in \mathcal{P}^\theta(M)$ where $\overline{P}:=\theta(P)$. Recall that this $(G,M,\theta)$-orthogonal set induces a function $\Gamma_{M,\theta}^{G}(.,\mathcal{Y}(g))$ on $\mathcal{A}_{M,\theta}$ (see \S \ref{GMT families}). If $\mathcal{Y}(g)$ is a positive $(G,M,\theta)$-orthogonal set then this is just the characteristic function of the convex hull of $\{\mathcal{Y}(g)_P;  P\in \mathcal{P}^\theta(M)\}$. Define

$$\displaystyle \widetilde{v}_{M,\theta}(Y,g):=\int_{A_G(F)\backslash A(F)} \Gamma_{M,\theta}^{G}(H_{M,\theta}(a),\mathcal{Y}(g))da$$

\noindent for all $Y\in \mathcal{A}_{P_0,\theta}^+$ and all $g\in G(F)$. Fixing a norm $\lvert .\rvert$ on $\mathcal{A}_{0,\theta}$, by \ref{eq 1.2.3} there exists $k>0$ such that we have an inequality

\begin{align}\label{eq 4.3.4}
\displaystyle \lvert \widetilde{v}_{M,\theta}(Y,g)\rvert\ll (1+\lvert Y\rvert)^k \sigma_{T^G\backslash G}(g)^k
\end{align}

\noindent for all $Y\in \mathcal{A}_{P_0,\theta}^+$ and all $g\in G(F)$. Define the following expression

$$\displaystyle J^\chi_{Y,T}(f):=\nu(T)\int_{\overline{T}(F)} D^H(t)\int_{A_T(F)A(F)\backslash G(F)} f(g^{-1}tg) \widetilde{v}_{M,\theta}(Y,g)dg \chi(t)^{-1}dt$$

\noindent for all $Y\in \mathcal{A}_{P_0,\theta}^+$. Using \ref{eq 4.3.4} and reasoning as in the proof of Lemma \ref{lemma kappaN}(\ref{lemma kappaN ii}), we can show that this expression is absolutely convergent.

\begin{prop}\label{prop change of truncation}
Let $0<\epsilon_1<\epsilon_2<1$. Then, for all $k>0$ we have

$$\displaystyle \left\lvert J^\chi_{N,T}(f)-J^\chi_{Y,T}(f)\right\rvert\ll N^{-k} $$

\noindent for all $N\geqslant 1$ and all $Y\in \mathcal{A}_{P_0,\theta}^+$ satisfying the two inequalities

\begin{align}\label{eq 4.3.5}
\displaystyle N^{\epsilon_1}\leqslant \inf_{\alpha\in \Delta_0} \alpha(Y)
\end{align}

\begin{align}\label{eq 4.3.6}
\displaystyle \sup_{\alpha\in \Delta_0} \alpha(Y)\leqslant N^{\epsilon_2}
\end{align}
\end{prop}

\vspace{2mm}

\noindent\ul{Proof}: Let $0<\epsilon_1<\epsilon_2<1$. For $M>0$ we will denote by $\mathbf{1}_{<M}$ and $\mathbf{1}_{\geqslant M}$ the characteristic functions of the sets of $g\in G(F)$ satisfying $\displaystyle \inf_{a\in A(F)A_T(F)}\sigma(ag)<M$ and $\displaystyle\inf_{a\in A(F)A_T(F)}\sigma(ag)\geqslant M$ respectively. For all $M>0$, we can write

$$\displaystyle J^\chi_{N,T}(f)=J^\chi_{N,T,< M}(f)+J^\chi_{N,T,\geqslant M}(f)$$

$$\displaystyle J^\chi_{Y,T}(f)=J^\chi_{Y,T,< M}(f)+J^\chi_{Y,T,\geqslant M}(f)$$

\noindent for all $N\geqslant 1$ and all $Y\in \mathcal{A}_{P_0,\theta}^+$, where

$$\displaystyle J^\chi_{N,T,<M}(f):=\nu(T)\int_{\overline{T}(F)} D^H(t)\int_{A_T(F)A(F)\backslash G(F)}\mathbf{1}_{<M}(g)f(g^{-1}tg)\kappa_{N,A}(g)dg \chi(t)^{-1}dt$$

$$\displaystyle J^\chi_{N,T,\geqslant M}(f):=\nu(T)\int_{\overline{T}(F)} D^H(t)\int_{A_T(F)A(F)\backslash G(F)}\mathbf{1}_{\geqslant M}(g)f(g^{-1}tg)\kappa_{N,A}(g)dg \chi(t)^{-1}dt$$

\noindent and $J^\chi_{Y,T,<M}(f)$, $J^\chi_{Y,T,\geqslant M}(f)$ are defined by similar expressions. First we show

\vspace{3mm}

\begin{num}
\item\label{eq 4.3.7} For all $\epsilon>0$ and all $k>0$ we have
$$\displaystyle \lv J^\chi_{N,T,\geqslant N^\epsilon}(f)\rvert\ll N^{-k}$$
and
$$\displaystyle \lv J^\chi_{Y,T,\geqslant N^\epsilon}(f)\rvert\ll N^{-k}$$
for all $N\geqslant 1$ and all $Y\in \mathcal{A}_{P_0,\theta}^+$ satisfying inequality \ref{eq 4.3.6}.
\end{num}

\vspace{3mm} 

\noindent By \ref{eq 4.3.3} and \ref{eq 4.3.4} and the fact that $f$ belongs to the Harish-Chandra-Schwartz space $\mathcal{C}_{\omega}(G(F))$, we only need to show that for all $k,k'>0$ and $\epsilon>0$ there exists $d>0$ such that

$$\displaystyle \int_{\overline{T}(F)} D^H(t)\int_{A(F)A_T(F)\backslash G(F)} \mathbf{1}_{\geqslant N^\epsilon}(g) \Xi^G(g^{-1}tg)\overline{\sigma}(g^{-1}tg)^{-d} \sigma_{T^G\backslash G}(g)^kdgdt\ll N^{-k'}$$

\noindent for all $N\geqslant 1$. By \ref{eq 4.3.2}, for all $r>0$ this integral is essentially bounded by 

$$\displaystyle N^{-\epsilon r}\int_{\overline{T}(F)} D^H(t)\int_{T^G(F)\backslash G(F)} \Xi^G(g^{-1}tg)\overline{\sigma}(g^{-1}tg)^{-d} \sigma_{T^G\backslash G}(g)^{k+r}dgdt$$

\noindent for all $N\geqslant 1$. By \ref{eq 4.1.1} and Lemma \ref{lemma conv abs woi} for all $k$, $r>0$ there exists $d>0$ making the last integral above convergent. The claim follows.

\vspace{2mm}

\noindent Choose $\epsilon>0$ such that $\epsilon<\epsilon_1$. By \ref{eq 4.3.7}, it suffices to show that for all $k>0$ we have

\begin{align}\label{eq 4.3.8}
\displaystyle \left\lvert J^\chi_{N,T,<N^\epsilon}(f)-J^\chi_{Y,T,<N^\epsilon}(f)\right\rvert\ll N^{-k}
\end{align}

\noindent for all $N\geqslant 1$ and all $Y\in \mathcal{A}_{P_0,\theta}^+$ satisfying inequalities \ref{eq 4.3.5} and \ref{eq 4.3.6}. For $Q\in \mathcal{F}^\theta(M)$, $Y\in \mathcal{A}_{P_0,\theta}^+$ and $N\geqslant 1$, we set

$$\displaystyle \kappa^{Y,Q}_{N,A}(g):=\int_{A_G(F)\backslash A(F)} \Gamma_{M,\theta}^{Q}(H_{M,\theta}(a),\mathcal{Y}(g)) \tau_{Q,\theta}^{G}(H_{M,\theta}(a)-\mathcal{Y}(g)_Q) \kappa_N(ag) da$$

\noindent for all $g\in G(F)$ where the functions $\Gamma_{M,\theta}^{Q}(.,\mathcal{Y}(g))$ and $\tau_{Q,\theta}^{G}$ have been defined in \S \ref{GMT families}. Note that

$$\displaystyle \Gamma_{M,\theta}^{Q}(.,\mathcal{Y}(ag))=\Gamma_{M,\theta}^{Q}(.+H_{M,\theta}(a),\mathcal{Y}(g))$$

\noindent for all $a\in A(F)$. Hence the functions $\kappa_{N,A}^{Y,Q}$ are $A(F)$-invariant on the left and this allows us to define the following expressions

$$\displaystyle J_{N,T,<N^\epsilon}^{\chi,Y,Q}(f):=\int_{\overline{T}(F)} D^H(t)\int_{A(F)A_T(F)\backslash G(F)} \mathbf{1}_{<N^\epsilon}(g) f(g^{-1}tg)\kappa_{N,A}^{Y,Q}(g)dg\chi(t)^{-1}dt$$

\noindent for all $N\geqslant 1$, all $Y\in \mathcal{A}_{P_0,\theta}^+$ and all $Q\in \mathcal{F}^\theta(M)$. By \ref{eq 1.7.2}, we have

$$\displaystyle J^\chi_{N,T,<N^\epsilon}(f)=\sum_{Q\in \mathcal{F}^\theta(M)} J_{N,T,<N^\epsilon}^{\chi,Y,Q}(f)$$

\noindent for all $N\geqslant 1$ and all $Y\in \mathcal{A}_{P_0,\theta}^+$. Thus, to show \ref{eq 4.3.8} it suffices to establish the two following facts

\vspace{3mm}

\begin{num}
\item\label{eq 4.3.9} There exists $N_0\geqslant 1$ such that
$$\displaystyle J_{N,T,<N^\epsilon}^{\chi,Y,G}(f)=J_{Y,T,<N^\epsilon}^{\chi}(f)$$
for all $N\geqslant N_0$ and all $Y\in \mathcal{A}_{P_0,\theta}^+$ satisfying inequality \ref{eq 4.3.6}.
\end{num}

\vspace{3mm}

\begin{num}
\item\label{eq 4.3.10} For all $Q\in \mathcal{F}^\theta(M)$, $Q\neq G$, and all $k>0$ we have
$$\displaystyle \left\lvert J_{N,T,<N^\epsilon}^{\chi,Y,Q}(f)\right\rvert\ll N^{-k}$$
for all $N\geqslant 1$ and all $Y\in \mathcal{A}_{P_0,\theta}^+$ satisfying inequality \ref{eq 4.3.5}.
\end{num} 

\vspace{3mm}

\noindent First we prove \ref{eq 4.3.9}. By definition of $J_{N,T,<N^\epsilon}^{\chi,Y,G}(f)$ and $J_{Y,T,<N^\epsilon}^{\chi}(f)$ it suffices to show that there exists $N_0\geqslant 1$ such that

$$\displaystyle \kappa_{N,A}^{Y,G}(g)=\widetilde{v}_{M,\theta}(Y,g)$$

\noindent for all $N\geqslant N_0$, all $Y\in \mathcal{A}_{P_0,\theta}^+$ satisfying inequality \ref{eq 4.3.6} and all $g\in G(F)$ with $\sigma(g)<N^\epsilon$. Unraveling the definitions of $\kappa_{N,A}^{Y,Q}(g)$ and $\widetilde{v}_{M,\theta}(Y,g)$, we see that it would follow if we can show the implication

\begin{align}\label{eq 4.3.11}
\displaystyle \Gamma_{M,\theta}^{G}\left(H_{M,\theta}(a),\mathcal{Y}(g) \right)\neq 0\Longrightarrow \kappa_N(ag)=1
\end{align}

\noindent for all $N\gg 1$, all $Y\in \mathcal{A}_{P_0,\theta}^+$ satisfying inequality \ref{eq 4.3.6}, all $g\in G(F)$ with $\sigma(g)<N^\epsilon$ and all $a\in A(F)$. By \ref{eq 1.7.1}, there exists $C>0$ such that for all $Y\in \mathcal{A}_{P_0,\theta}^+$, all $g\in G(F)$ and all $a\in A(F)$ we have

$$\displaystyle \Gamma_{M,\theta}^{G}\left(H_{M,\theta}(a),\mathcal{Y}(g) \right)\neq 0\Longrightarrow \sigma_{A_G\backslash A}(a)\leqslant C\left(\sup_{\alpha\in \Delta_0}\alpha(Y)+\sigma(g) \right)$$

\noindent As $\sigma_{A_G\backslash A}(a)\sim \sigma_X(a)$ (this is a consequence of the facts that $A^\theta A_G\backslash A$ is closed in $\mathbf{X}$ and $A_G\backslash A^\theta A_G$ is finite) and $\sigma_X(ag)\ll \sigma_X(a)+\sigma(g)$ for all $a\in A(F)$ and all $g\in G(F)$, it follows that there exists $C'>0$ such that for all $N\geqslant 1$, all $Y\in \mathcal{A}_{P_0,\theta}^+$ satisfying inequality \ref{eq 4.3.6}, all $g\in G(F)$ with $\sigma(g)<N^\epsilon$ and all $a\in A(F)$ we have

$$\displaystyle \Gamma_{M,\theta}^{G}\left(H_{M,\theta}(a),\mathcal{Y}(g) \right)\neq 0\Longrightarrow \sigma_{X}(ag)\leqslant C'\left(N^{\epsilon_1}+N^\epsilon \right)$$

\noindent Since $\epsilon,\epsilon_1<1$, by property \ref{eq 4.2.1} of our sequence of truncation functions the last inequality above implies $\kappa_N(ag)=1$ whenever $N\gg 1$. This shows \ref{eq 4.3.11} and ends the proof of \ref{eq 4.3.9}.

\vspace{2mm}

\noindent It only remains to prove claim \ref{eq 4.3.10}. Fix $Q\in \mathcal{F}^\theta(M)$, $Q\neq G$, with Levi decomposition $Q=LU_Q$ where $L:=Q\cap\theta(Q)$. Let $\overline{Q}=\theta(Q)=L\overline{U}_Q$ be the opposite parabolic subgroup. We have the Iwasawa decomposition $G(F)=L(F)\overline{U}_Q(F)K$ and accordingly we can decompose the integral

\[\begin{aligned}
\displaystyle \int_{A(F)A_T(F)\backslash G(F)} \mathbf{1}_{<N^\epsilon}(g) f(g^{-1}tg)  & \kappa_{N,A}^{Y,Q}(g)dg= \\
 & \int_{A(F)A_T(F)\backslash L(F)\times \overline{U}_Q(F)\times K}  \mathbf{1}_{<N^\epsilon}(luk) f(k^{-1}u^{-1}l^{-1}tluk) \kappa_{N,A}^{Y,Q}(luk)dl du dk
\end{aligned}\]

\noindent for all $N\geqslant 1$, all $Y\in \mathcal{A}_{P_0,\theta}^+$ and all $t\in \overline{T}_{\reg}(F)$. To continue we need the following fact which we will establish after we finish the proof of \ref{eq 4.3.10}:

\vspace{3mm}

\begin{num}
\item\label{eq 4.3.12} There exists $N_0\geqslant 1$ such that for all $N\geqslant N_0$, all $Y\in \mathcal{A}_{P_0,\theta}^+$ satisfying inequality \ref{eq 4.3.5} and all $l\in L(F)$, $u\in \overline{U}_Q(F)$, $k\in K$ with $\sigma(luk)<N^\epsilon$ we have

$$\displaystyle \kappa_{N,A}^{Y,Q}(luk)=\kappa_{N,A}^{Y,Q}(lk)$$
\end{num}

\vspace{3mm}

\noindent Taking \ref{eq 4.3.12} for granted we get

\[\begin{aligned}
\displaystyle \int_{A(F)A_T(F)\backslash G(F)} \mathbf{1}_{<N^\epsilon}(g) f(g^{-1}tg)  & \kappa_{N,A}^{Y,Q}(g)dg= \\
 & \int_{A(F)A_T(F)\backslash L(F)\times \overline{U}_Q(F)\times K}  \mathbf{1}_{<N^\epsilon}(luk) f(k^{-1}u^{-1}l^{-1}tluk) \kappa_{N,A}^{Y,Q}(lk)dl du dk
\end{aligned}\]

\noindent for all $N\gg 1$, all $t\in \overline{T}_{\reg}(F)$ and all $Y\in \mathcal{A}_{P_0,\theta}^+$ satisfying inequality \ref{eq 4.3.5}. As $f$ is strongly cuspidal if we forget the term $\mathbf{1}_{<N^\epsilon}(luk)$ in the integrand of the last integral above we get zero. Thus, we have

\[\begin{aligned}
\displaystyle \int_{A(F)A_T(F)\backslash G(F)} \mathbf{1}_{<N^\epsilon}(g) f(g^{-1}tg)  & \kappa_{N,A}^{Y,Q}(g)dg= \\
 & -\int_{A(F)A_T(F)\backslash L(F)\times \overline{U}_Q(F)\times K}  \mathbf{1}_{\geqslant N^\epsilon}(luk) f(k^{-1}u^{-1}l^{-1}tluk) \kappa_{N,A}^{Y,Q}(lk)dl du dk
\end{aligned}\]

\noindent for all $N\gg 1$, all $t\in \overline{T}_{\reg}(F)$ and all $Y\in \mathcal{A}_{P_0,\theta}^+$ satisfying inequality \ref{eq 4.3.5}. Hence, to get \ref{eq 4.3.10} it only remains to show that for all $k>0$ we have

\begin{align}\label{eq 4.3.13}
\displaystyle \int_{\overline{T}(F)}D^H(t)\int_{A(F)A_T(F)\backslash L(F)\times \overline{U}_Q(F)\times K}  \mathbf{1}_{\geqslant N^\epsilon}(luk) \left\lvert f(k^{-1}u^{-1}l^{-1}tluk)\right\rvert \left\lvert\kappa_{N,A}^{Y,Q}(lk)\right\rvert dl du dk dt\ll N^{-k}
\end{align}

\noindent for all $N\geqslant 1$ and all $Y\in \mathcal{A}_{\theta,P_0}^+$. Since for all $(G,M,\theta)$-orthogonal set $\mathcal{Z}=(\mathcal{Z}_P)_{P\in \mathcal{P}^\theta(M)}$ the function $X\in \mathcal{A}_\theta\mapsto \Gamma^{Q}_{M,\theta}(X,\mathcal{Z})$ is uniformly bounded independently of $\mathcal{Z}$ (this follows from the definition of this function) and $\tau_{Q,\theta}^{G}$ is a characteristic function, we have

$$\displaystyle \left\lvert\kappa_{N,A}^{Y,Q}(g)\right\rvert\ll \kappa_{N,A}(g)$$

\noindent for all $N\geqslant 1$, all $Y\in \mathcal{A}_{P_0,\theta}^+$ and all $g\in G(F)$. Hence, by \ref{eq 4.3.3} and the fact that $f$ belongs to the Harish-Chandra-Schwartz space $\mathcal{C}_\omega(G(F))$ there exists $k>0$ such that for all $d>0$ the left hand side of \ref{eq 4.3.13} is essentially bounded by the product of $N^k$ with

\[\begin{aligned}
\displaystyle \int_{\overline{T}(F)}D^H(t)\int_{A(F)A_T(F)\backslash L(F)\times \overline{U}_Q(F)\times K}  \mathbf{1}_{\geqslant N^\epsilon}(luk) \Xi^G(k^{-1}u^{-1}l^{-1}tluk) & \overline{\sigma}(k^{-1}u^{-1}l^{-1}tluk)^{-d} \\
 & \sigma_{T^G\backslash G}(lk)^k dldudkdt
\end{aligned}\]

\noindent for all $N\geqslant 1$ and all $Y\in \mathcal{A}_{P_0,\theta}^+$. Since $T^G\subset L$ we have $\sigma_{T^G\backslash G}(lk)\ll \sigma_{T^G\backslash G}(luk)$ for all $l\in L(F)$, $u\in \overline{U}_Q(F)$ and all $k\in K$. Thus, the last expression above is essentially bounded by

$$\displaystyle \int_{\overline{T}(F)}D^H(t)\int_{A(F)A_T(F)\backslash G(F)}  \mathbf{1}_{\geqslant N^\epsilon}(g) \Xi^G(g^{-1}tg) \overline{\sigma}(g^{-1}tg)^{-d}\sigma_{T^G\backslash G}(g)^k dg$$

\noindent for all $N\geqslant 1$. We already saw that for all $k'>0$ we can find $d>0$ such that this last integral is essentially bounded by $N^{-k'}$ for all $N\geqslant 1$. This shows \ref{eq 4.3.13} and ends the proof of \ref{eq 4.3.10} granting \ref{eq 4.3.12}.

\vspace{2mm}

\noindent We now prove \ref{eq 4.3.12}. For $g\in G(F)$, the function $\Gamma_{M,\theta}^{Q}(.,\mathcal{Y}(g)) \tau_{Q,\theta}^{G}(.-\mathcal{Y}(g)_Q)$ depends only on the points $\mathcal{Y}(g)_{P}$ for all $P\in \mathcal{P}^\theta(M)$ with $P\subset Q$ and those points remain invariant by left translation of $g$ by $\overline{U}_Q(F)$. Hence, it suffices to show the following:

\vspace{3mm}

\begin{num}
\item\label{eq 4.3.14} There exists $N_0\geqslant 1$ such that for all $N\geqslant N_0$, all $Y\in \mathcal{A}_{P_0,\theta}^+$ satisfying inequality \ref{eq 4.3.5}, all $l\in L(F)$, $u\in \overline{U}_Q(F)$, $k\in K$ with $\sigma(luk)<N^\epsilon$ and all $a\in A(F)$ we have

$$\displaystyle \Gamma_{M,\theta}^{Q}(H_{M,\theta}(a),\mathcal{Y}(l)) \tau_{Q,\theta}^{G}(H_{M,\theta}(a)-\mathcal{Y}(l)_Q)\neq 0\Rightarrow \kappa_N(aluk)=\kappa_N(alk)$$
\end{num}

\vspace{3mm}

\noindent Let $N\geqslant 1$ and $Y$, $l$, $u$, $k$ be as above (in particular $Y$ satisfies condition \ref{eq 4.3.5} and $\sigma(luk)<N^\epsilon$). We will show that the conclusion of \ref{eq 4.3.14} holds provided $N$ is sufficiently large. Let $a\in A(F)$ be such that

$$\displaystyle \Gamma_{M,\theta}^{Q}(H_{M,\theta}(a),\mathcal{Y}(l)) \tau_{Q,\theta}^{G}(H_{\theta,M}(a)-\mathcal{Y}(l)_Q)\neq 0$$

\noindent We need to show that $\kappa_N(aluk)=\kappa_N(alk)$.

\vspace{2mm}

\noindent There exists $C>0$ such that for all $g\in G(F)$ and all $Y\in \mathcal{A}_{P_0,\theta}^+$ if

$$\displaystyle \sigma(g)\leqslant C\inf_{\alpha\in \Delta_0} \alpha(Y)$$

\noindent then $\mathcal{Y}(g)_P\in \mathcal{A}_{P,\theta}^+$ for all $P\in \mathcal{P}^\theta(M)$ and thus $\mathcal{Y}(g)$ is a positive $(G,M,\theta)$-orthogonal set. As $\epsilon<\epsilon_1$ it follows that for $N$ sufficiently large the $(G,M,\theta)$-orthogonal set $\mathcal{Y}(l)$ is positive. In particular, again for $N$ sufficiently large, the function

$$\displaystyle X\in \mathcal{A}_{M,\theta}\mapsto \Gamma_{M,\theta}^{Q}(X,\mathcal{Y}(l)) \tau_{Q,\theta}^{G}(X-\mathcal{Y}(l)_Q)$$

\noindent is the characteristic function of the sum of $\mathcal{A}_{Q,\theta}^+$ with the convex hull of the family $(\mathcal{Y}(l)_P)_{P\subset Q}$. As $\epsilon<\epsilon_1$ and $\sigma(l)\ll N^\epsilon$, it follows that

\begin{align}\label{eq 4.3.15}
\displaystyle \log \lvert \beta(a)\rvert\geqslant \inf_{P\subset Q} \beta\left(Y_P-H_{\overline{P},\theta}(l)\right)\gg \inf_{\alpha\in \Delta_0} \alpha(Y)-\sigma(l)\gg N^{\epsilon_1}
\end{align}

\noindent for all $\beta\in R(A,U_Q)$. Fix a norm $\lvert .\rvert$ on $\mathfrak{g}(F)$ and let us denote by $B(0,r)$ the open ball of radius $r$ centered at the origin for all $r>0$. Since $\sigma(lul^{-1})\ll N^\epsilon$ and $\epsilon<\epsilon_1$, we deduce from \ref{eq 4.3.15} that there exists a constant $c_1>0$ such that for $N$ big enough we have

$$\displaystyle alul^{-1}a^{-1}\in \exp\left(B(0,e^{-c_1N^{\epsilon_1}})\right)$$

\noindent Let $P_a\in \mathcal{P}^\theta(M)$ be such that $H_{M,\theta}(a)\in \overline{\mathcal{A}^+_{P_a,\theta}}$ (the closure of the positive chamber associated to $P_a$). Since $P_a$ is $\theta$-split, the multiplication map $H(F)\times P_a(F)\to G(F)$ is submersive at the origin and hence this map admits an $F$-analytic section defined on a neighborhood of $1$ in $G(F)$. It follows that there exists $c_2>0$ (independent of $a$ since there is only a finite number of possibilities for $P_a$) so that for $N$ large enough

$$\displaystyle alul^{-1}a^{-1}\in H(F)\exp\left(B(0,e^{-c_2N^{\epsilon_1}})\cap \mathfrak{p}_a(F)\right)$$

\noindent Choose $X\in B(0,e^{-c_2N^{\epsilon_1}})\cap \mathfrak{p}_a(F)$ with $alul^{-1}a^{-1}\in H(F)\exp(X)$. Since $\kappa_N$ is left invariant by $H(F)$ we have

$$\displaystyle \kappa_N(aluk)=\kappa_N(\exp(X)alk)$$

\noindent As $a\in \overline{\mathcal{A}^+_{P_a,\theta}}$, $\epsilon<\epsilon_1$ and $\sigma(l)\ll N^\epsilon$, there exists a constant $c_3>0$ such that

$$\displaystyle k^{-1}l^{-1}a^{-1}Xalk\in B(0, e^{-c_3N^{\epsilon_1}})$$

\noindent By property \ref{eq 4.2.2} of our sequence of truncation functions, we deduce that for $N$ sufficiently large $\kappa_N$ is right invariant by $\exp(k^{-1}l^{-1}a^{-1}Xalk)$. Hence,

$$\displaystyle \kappa_N(\exp(X)alk)=\kappa_N(alk)$$

\noindent This proves claim \ref{eq 4.3.14} and ends the proof of the proposition. $\blacksquare$

\subsection{First computation of the limit}

Recall the function $g\in G(F)\mapsto v_{M,\theta}(g)$ introduced in \S \ref{GMT families}. Its value at $g\in G(F)$ is given by the volume of the convex hull of the set $\{H_{P,\theta}(g);\; P\in \mathcal{P}^\theta(M) \}$.

\begin{prop}\label{prop first computation}
We have

$$\displaystyle \lim\limits_{N\to \infty} J_{N,T}^\chi(f)=(-1)^{a^G_{M,\theta}}\nu(T)\int_{\overline{T}(F)} D^H(t)\int_{A_T(F)A(F)\backslash G(F)} f(g^{-1}tg)v_{M,\theta}(g)dg \chi(t)^{-1}dt$$
\end{prop}

\vspace{2mm}

\noindent\ul{Proof}: Let $0<\epsilon_1<\epsilon_2<1$, $0<\delta<1$ and set

$$\displaystyle \mathcal{A}_{P_0,\theta}^+(\delta):=\{Y\in \mathcal{A}_{P_0,\theta}^+;\; \inf_{\alpha\in \Delta_0}\alpha(Y)\geqslant \delta\sup_{\alpha\in \Delta_0}\alpha(Y) \}$$

\noindent Then $\mathcal{A}_{P_0,\theta}^+(\delta)$ is a cone in $\mathcal{A}_{0,\theta}$ with nonempty interior. By Proposition \ref{prop change of truncation} for all $k>0$ we have

$$\displaystyle \left\lvert J^\chi_{N,T}(f)-J^\chi_{Y,T}(f)\right\rvert\ll N^{-k}$$

\noindent for all $N\geqslant 1$ and all $Y\in \mathcal{A}_{P_0,\theta}^+(\delta)$ with $\displaystyle N^{\epsilon_1}\leqslant \inf_{\alpha\in \Delta_0}\alpha(Y)\leqslant \delta^{-1}N^{\epsilon_2}$. As for $N$ sufficiently large the two sets

$$\displaystyle \left\{ Y\in \mathcal{A}_{P_0,\theta}^+(\delta);\; N^{\epsilon_1}\leqslant \inf_{\alpha\in \Delta_0}\alpha(Y)\leqslant \delta^{-1} N^{\epsilon_2} \right\}$$

\noindent and

$$\displaystyle \left\{ Y\in \mathcal{A}_{P_0,\theta}^+(\delta);\; (N+1)^{\epsilon_1}\leqslant \inf_{\alpha\in \Delta_0}\alpha(Y)\leqslant \delta^{-1} (N+1)^{\epsilon_2} \right\}$$

\noindent intersect, it follows that the two limits 

$$\displaystyle \lim\limits_{N\to \infty} J^\chi_{N,T}(f),\; \lim\limits_{Y\in \mathcal{A}_{P_0,\theta}^+(\delta)\to \infty} J^\chi_{Y,T}(f)$$

\noindent exist and are equal. We will denote by $J_{\infty,T}^\chi(f)$ this common limit.

\vspace{2mm}

\noindent Let $\mathcal{A}_{0,\theta,F}$ denote the image of $A_{0,\theta}(F)$ by $H_{M_0,\theta}$. Then by Lemma \ref{lem GMT families}, we know that for every lattice $\mathcal{R}\subset \mathcal{A}_{0,\theta,F}\otimes \mathbf{Q}$ and all $g\in G(F)$ the function $Y\in \mathcal{R}\cap \mathcal{A}_{P_0,\theta}^+(\delta)\mapsto \widetilde{v}_{M,\theta}(Y,g)$ coincides with the restriction of an exponential-polynomial of bounded degree and with exponents in a fixed finite set (both independent of $g$). Let us denote by $\widetilde{v}_{M,\theta,0}(\mathcal{R},g)$ the constant term of the purely polynomial part of this exponential-polynomial. Then by Lemma \ref{lemma polynomial-exponential}, we have

\begin{align}\label{eq 4.4.1}
\displaystyle J_{\infty,T}^\chi(f)=\nu(T)\int_{\overline{T}(F)} D^H(t)\int_{A_T(F)A(F)\backslash G(F)} f(g^{-1}tg) \widetilde{v}_{M,\theta,0}(\mathcal{R},g)dg \chi(t)^{-1}dt
\end{align}

\noindent for every lattice $\mathcal{R}\subset \mathcal{A}_{0,\theta,F}\otimes \mathbf{Q}$. Fix such a lattice $\mathcal{R}$. By Lemma \ref{lem GMT families}, there exists $r>0$ such that

$$\displaystyle \left\lvert \widetilde{v}_{M,\theta,0}(\frac{1}{k}\mathcal{R},g)-(-1)^{a^G_{M,\theta}}v_{M,\theta}(g)\right\rvert\ll \sigma(g)^r k^{-1}$$

\noindent for all $k\geqslant 1$ and all $g\in G(F)$. Since the left hand side is invariant by left translation of $g$ by $T^G(F)$, we also have

$$\displaystyle \left\lvert \widetilde{v}_{M,\theta,0}(\frac{1}{k}\mathcal{R},g)-(-1)^{a^G_{M,\theta}}v_{M,\theta}(g)\right\rvert\ll \sigma_{T^G\backslash G}(g)^r k^{-1}$$

\noindent for all $k\geqslant 1$ and all $g\in G(F)$. By Lemma \ref{lemma conv abs woi}, this implies that

\[\begin{aligned}
\displaystyle \lim\limits_{k\to \infty} \int_{\overline{T}(F)} D^H(t)\int_{A_T(F)A(F)\backslash G(F)} & f(g^{-1}tg) \widetilde{v}_{M,\theta,0}(\frac{1}{k}\mathcal{R},g)dg \chi(t)^{-1}dt= \\
 & (-1)^{a^G_{M,\theta}}\int_{\overline{T}(F)} D^H(t)\int_{A_T(F)A(F)\backslash G(F)} f(g^{-1}tg) v_{M,\theta}(g)dg \chi(t)^{-1}dt
\end{aligned}\]

\noindent From this and \ref{eq 4.4.1} (which is of course also true if we replace $\mathcal{R}$ by $\frac{1}{k}\mathcal{R}$, $k\geqslant 1$) we deduce the proposition. $\blacksquare$

\subsection{End of the proof}\label{proof geo side}

By the descent formula \ref{eq 1.7.4} and Proposition \ref{prop first computation} we have

$$\displaystyle \lim\limits_{N\to\infty} J^\chi_{N,T}(f)= (-1)^{a^G_{M,\theta}}\nu(T)\left[A_{T^G}(F):A_T(F)A(F) \right]\sum_{L\in \mathcal{L}(M)} d_{M,\theta}^{G}(L) \int_{\overline{T}(F)} \Phi_M^Q(t,f) \chi(t)^{-1}dt$$

\noindent Since $f$ is strongly cuspidal only the term corresponding to $L=G$ can contribute to the sum above so that

\begin{align}\label{eq 4.5.1}
\displaystyle \lim\limits_{N\to\infty} J^\chi_{N,T}(f)= (-1)^{a^G_{M,\theta}}\nu(T)\left[A_{T^G}(F):A_T(F)A(F) \right] d_{M,\theta}^{G}(G) \int_{\overline{T}(F)} \Phi_M(t,f) \chi(t)^{-1}dt
\end{align}

\noindent Assume first that $T$ is not elliptic in $H$. We distinguish two cases:

\begin{itemize}
\item If $M\neq \Cent_G(A_{T^G})$ we have $\Phi_M(t,f)=0$ for all $t\in \overline{T}_{\reg}(F)$ as $M\neq M(t)$ so that the limit \ref{eq 4.5.1} vanishes.

\item If $M=\Cent_G(A_{T^G})$, we have $\mathcal{A}_{M}^{G,\theta}=\mathcal{A}_T^H\neq 0$ (as $T$ is not elliptic in $H$), thus $d_{M,\theta}^G(G)=0$ and the limit \ref{eq 4.5.1} also vanishes in this case.
\end{itemize}

\noindent Hence, in both cases the limit \ref{eq 4.5.1} equals zero for $T$ nonelliptic in $H$. Now, if $T$ is elliptic in $H$ we have $A_T=A_H$, $A_{T^G}=A$ and $\mathcal{A}_{M}^{G,\theta}=0$ so that $d_{M,\theta}^{G}(G)=1$, $\left[A_{T^G}(F):A_T(F)A(F) \right]=1$ and $\nu(T)=\nu(H)$ and we get

$$\lim\limits_{N\to\infty} J^\chi_{N,T}(f)= (-1)^{a^G_{M}}\nu(H) \int_{\overline{T}(F)} \Phi_M(t,f) \chi(t)^{-1}dt=\nu(H)\int_{\overline{T}(F)} D^H(t)\Theta_f(t) \chi(t)^{-1}dt$$

\noindent Theorem \ref{theo geom side} now follows from the above equality, \ref{eq 4.2.3} and \ref{eq 4.2.5}. $\blacksquare$

\section{Applications to a conjecture of Prasad}

In this chapter $E/F$ is a quadratic extension, $H$ is a connected reductive group over $F$ and $G:=R_{E/F}H_E$. We will denote by $\theta$ the involution of $G$ induced by the nontrivial element of $\Gal(E/F)$. Hence $H=G^\theta$. As before we set $\overline{G}:=G/A_G$ and $\overline{H}:=H/A_H$. If $Q$ is an algebraic subgroup of $H$ then $R_{E/F}Q_E$ is an algebraic subgroup of $G$. Also, note that if $P$ (resp. $M$) is a parabolic (resp. Levi) subgroup of $G$ which is $\theta$-stable (i.e. $\theta(P)=P$, resp. $\theta(M)=M$) then there exists a parabolic (resp. Levi) subgroup $\mathcal{P}$ (resp. $\mathcal{M}$) of $H$ such that $P=R_{E/F}\mathcal{P}_E$ (resp. $M=R_{E/F}\mathcal{M}_E$). We will denote by $\mathcal{R}(G)$ the space of {\it virtual} representations of $G(F)$ that is the complex vector spaces with basis $\Irr(G)$. Similarly, if $A$ is an abelian group we will denote by $\mathcal{R}(A)$ the space of virtual characters of $A$. We will write $H^i(F,.)$ for the functors of Galois cohomology and if $\mathcal{H}$, $\mathcal{G}$ are algebraic groups over $F$ with $\mathcal{H}$ a subgroup of $\mathcal{G}$ we will set

$$\displaystyle \ker1^1(F;\mathcal{H},\mathcal{G}):=\Ker\left(H^1(F,\mathcal{H})\to H^1(F,\mathcal{G}) \right)$$

\noindent By \cite{KottSt} Theorem 1.2, for every connected reductive group over $F$ there exists a natural structure of abelian group on $H^1(F,\mathcal{G})$ which is uniquely characterized by the fact that for every elliptic maximal torus $T\subset \mathcal{G}$ the natural map $H^1(F,T)\to H^1(F,\mathcal{G})$ is a group morphism. Moreover, for every connected reductive groups $\mathcal{H}$, $\mathcal{G}$ and every morphism $\mathcal{H}\to \mathcal{G}$ the induced map $H^1(F,\mathcal{H})\to H^1(F,\mathcal{G})$ is a group morphism. Indeed, if $\mathcal{T}$ is an elliptic maximal torus in $\mathcal{H}$ (whose existence is guaranteed by \cite{KnII} p.271) and $\mathcal{T}'$ is a maximal torus of $\mathcal{G}$ containing the image of $\mathcal{T}$, then we have a commuting square

$$\xymatrix{H^1(F,\mathcal{T}) \ar@{->}[r] \ar@{->}[d] & H^1(F,\mathcal{T}') \ar@{->}[d] \\ H^1(F,\mathcal{H}) \ar@{->}[r] & H^1(F,\mathcal{G})}$$

\noindent where the upper, left and right arrows are morphisms of groups and moreover the map $H^1(F,\mathcal{T})\to H^1(F,\mathcal{H})$ is surjective (\cite{KottSt} Lemma 10.2). From these, it easily follows that $H^1(F,\mathcal{H})\to H^1(F,\mathcal{G})$ is a morphism of abelian groups.

\vspace{2mm}

\noindent Finally, for every $\pi\in \Irr(G)$ and every continuous character $\chi$ of $H(F)$ we recall that in \ref{section statement spectral side} we have defined a multiplicity

$$\displaystyle m(\pi,\chi):=\dim \Hom_H(\pi,\chi)$$

\noindent which is always finite by \cite{Del}, Theorem 4.5. The function $\pi\in \Irr(G)\mapsto m(\pi,\chi)$ extends by linearity to $\mathcal{R}(G)$.

\subsection{A formula for the multiplicity}\label{section formula for the multiplicity}

We will denote by $\Gamma_{\elli}(\overline{H})$ the set of regular elliptic conjugacy classes in $\overline{H}(F)$ and we equip this set with a topology and a measure characterized by the fact that for all $x\in \overline{H}_{\reg}(F)$ the map $t\in G_x(F)\mapsto tx\in \Gamma_{\elli}(\overline{H})$, which is well-defined in a neighborhood of the identity, is a local isomorphism preserving measures near $1$ (recall that in \S \ref{Groups, measures, notations} we have fixed Haar measures on the $F$-points of any torus and in particular on $G_x(F)$). More concretely, if we fix a set $\mathcal{T}_{\elli}(H)$ of representatives of the $H(F)$-conjugacy classes of elliptic maximal tori in $H$, then for every integrable function $\varphi$ on $\Gamma_{\elli}(\overline{H})$ we have the following integration formula

$$\displaystyle \int_{\Gamma_{\elli}(\overline{H})}\varphi(x)dx=\sum_{T\in \mathcal{T}_{\elli}(H)} \lvert W(H,T)\rvert^{-1} \int_{\overline{T}(F)} \varphi(t) dt$$

\noindent where we have set $\overline{T}:=T/A_H$ for all $T\in \mathcal{T}_{\elli}(H)$ and we recall since $\overline{T}$ is anisotropic, the Haar measure on $\overline{T}(F)$ is of total mass $1$. For all $\pi\in \Irr(G)$ and every continuous character $\chi$ of $H(F)$ with $\omega_{\pi \mid A_H(F)}=\chi_{\mid A_H(F)}$, set

$$\displaystyle m_{\geom}(\pi,\chi):=\int_{\Gamma_{\elli}(\overline{H})}D^H(x)\Theta_\pi(x)\chi(x)^{-1}dx$$

\noindent This expression makes sense since semisimple regular elements of $H$ are also semisimple and regular in $G$ and the function $x\in H_{\reg}(F)\mapsto D^H(x)^{1/2}\Theta_\pi(x)$ is locally bounded on $H(F)$ by \cite{HCDeBS} Theorem 16.3 and the identity $D^H(x)=D^G(x)^{1/2}$.

Recall that we are denoting by $\Irr_{\sqr}(G)$ the set of (equivalence classes of) irreducible essentially square-integrable representations of $G(F)$. The main theorem of this section is the following.

\begin{theo}\label{theo formula multiplicity}
For all $\pi\in \Irr_{\sqr}(G)$ and every continuous character $\chi$ of $H(F)$ with $\omega_{\pi\mid A_H(F)}=\chi_{\mid A_H(F)}$ we have

$$\displaystyle m(\pi,\chi)=m_{\geom}(\pi,\chi)$$
\end{theo}

\noindent\ul{Proof}: Up to twisting $\pi$ and $\chi$ by real unramified characters, we may assume that $\omega_\pi$ and $\chi$ are unitary. Set $\omega:=\omega_{\pi\mid A_G(F)}$. By Theorem \ref{theo spec} and Theorem \ref{theo geom side}, for all $f\in {}^0\mathcal{C}_\omega(G)$ we have

$$\displaystyle \sum_{\sigma\in \Irr_{\omega,\sqr}(G)}m(\sigma,\chi)\Tr(\sigma^\vee(f))=\int_{\Gamma_{\elli}(\overline{H})}D^H(x)\Theta_f(x)\chi(x)^{-1}dx$$

\noindent By \ref{eq 1.4.1} and Proposition \ref{prop coeff char}, when we apply this equality to a coefficient of $\pi$ we get the identity of the theorem. $\blacksquare$ 

\subsection{Galoisian characters and Prasad's character $\omega_{H,E}$}\label{section galoisian chars}

Let $\check{H}$ denote the complex dual group of $H$, $Z(\check{H})$ be its center and $W_F$ be the Weil group of $F$. Denoting by $\Hom_{cont}(H(F),\mathbf{C}^\times)$ the group of continuous characters of $H(F)$, Langlands has defined an homomorphism

$$\displaystyle \alpha_H:H^1(W_F,Z(\check{H}))\to \Hom_{cont}(H(F),\mathbf{C}^\times)$$

\noindent which is injective since $F$ is $p$-adic but is not always surjective although it is most of the time (e.g. if $H$ is quasi-split). We refer the reader to \cite{LL} for discussion of these matters. We will call the image of $\alpha_H$ the set of {\it Galoisian characters} (of $H(F)$). Assume that $H$ is semi-simple. Let $H_{sc}$ be the simply connected cover of $H$ and $\pi_1(H)$ be the kernel of the projection $H_{sc}\to H$. By Tate-Nakayama duality, we have an isomorphism $H^1(W_F,Z(\check{H}))\simeq H^1(F,\pi_1(H))^D$, where $(.)^D$ denotes duality for finite abelian groups, and the morphism $\alpha_H$ is the composition of this isomorphism with the (dual of the) connecting map $\displaystyle H(F)\to H^1(F,\pi_1(H))$. In particular, in this case, a character of $H(F)$ is Galoisian if and only if it factorizes through $H^1(F,\pi_1(H))$.

\vspace{2mm}

In \cite{Pras15}, Prasad has defined a quadratic character $\omega_{H,E}:H(F)\to \{\pm 1 \}$ which depends not only on $H$ but also on the quadratic extension $E/F$. It is a Galoisian character whose simplest definition is as the image by $\alpha_H$ of the cocycle $c$ defined by $c(w)=1$ if $w\in W_E$ (the Weil group of $E$) and $c(w)=z$ if $w\in W_F\backslash W_E$, where $z$ denotes the image of the central element $\begin{pmatrix} -1 & \\ & -1\end{pmatrix}\in SL_2(\mathbf{C})$ by any principal $SL_2$-morphism $SL_2(\mathbf{C})\to \check{H}$. In what follows we shall need another description of Prasad's character (see \cite{Pras15} \S 8 for the equivalence between the two definitions). First of all, $\omega_{H,E}$ is the pullback by $H(F)\to H_{ad}(F)$ of $\omega_{H_{ad},E}$, where $H_{ad}$ denotes the adjoint group of $H$. Thus, to describe $\omega_{H,E}$ we may assume that $H$ is semisimple. We introduce notations as before: $H_{sc}$ is the simply connected cover of $H$ and $\pi_1(H)$ stands for the kernel of the projection $H_{sc}\to H$. Let $B\subset H_{sc,\overline{F}}$ and $T\subset B$ be a Borel subgroup and a maximal torus thereof (both a priori only defined over $\overline{F}$). Let $\rho\in X^*_{\overline{F}}(T)$ be the half sum of the positive roots of $T$ with respect to $B$ (this belongs to the character lattice of $T$ since $H_{sc}$ is simply-connected). Then, it can be easily shown that the restriction of $\rho$ to $\pi_1(H)$ induces a morphism $\pi_1(H)\to \mu_2$ defined over $F$. Pushing this through the inclusion $\mu_2\hookrightarrow \Ker N_{E/F}$ we get a morphism $\pi_1(H)\to \Ker N_{E/F}$ and $\omega_{H,E}$ is simply the composition of the connecting map $H_{ad}(F)\to H^1(F,\pi_1(H))$ with the corresponding homomorphism between Galoisian $H^1$'s:

$$\displaystyle H^1(F,\pi_1(H))\to H^1(F,\Ker N_{E/F})\simeq \{\pm 1 \}$$

\subsection{First application: comparison between inner forms}\label{section first application}

Let $H'$ be another connected reductive group over $F$ and let $\psi^H:H_{\overline{F}}\simeq H'_{\overline{F}}$ be an inner twisting. Set $G':=R_{E/F}H'$. Then $\psi^H$ induces an inner twisting $\psi^G:G_{\overline{F}}\simeq G'_{\overline{F}}$ (actually, there are natural isomorphisms $G_{\overline{F}}\simeq H_{\overline{F}}\times H_{\overline{F}}$, $G'_{\overline{F}}\simeq H'_{\overline{F}}\times H'_{\overline{F}}$ and using these as identifications we just have $\psi^G=\psi^H\times \psi^H$). 

Recall that two regular elements $x\in G_{\reg}(F)$ and $x'\in G'_{\reg}(F)$ are {\it stably conjugate} if there exists $g\in G(\overline{F})$ such that $x'=\psi^G(gxg^{-1})$ and the isomorphism $\psi^G\circ \Ad(g): G_{x,\overline{F}}\simeq G'_{x',\overline{F}}$ is defined over $F$. Similarly, two regular elements of $G(F)$ are {\it stably conjugate} if they are conjugate by an element of $G(\overline{F})$ which induces an isomorphism defined over $F$ between their connected centralizers.

We say that a virtual representation $\Pi\in \mathcal{R}(G)$ (or $\Pi'\in \mathcal{R}(G')$) is {\it stable} if its character $\Theta_\Pi$ (or $\Theta_{\Pi'}$) is constant on regular stable conjugacy classes in $G(F)$ (resp. in $G'(F)$). Two stable virtual representations $\Pi\in \mathcal{R}(G)$ and $\Pi'\in \mathcal{R}(G')$ are said to be {\it transfer of each other} if for all pairs $(x,x')\in G_{\reg}(F)\times G'_{\reg}(F)$ of stably conjugate regular elements we have $\Theta_\Pi(x)=\Theta_{\Pi'}(x')$. By the main results of \cite{Art4}, every stable virtual representation $\Pi\in \mathcal{R}(G)$ is the transfer of a stable virtual representation $\Pi'\in \mathcal{R}(G')$ and conversely.

We define similarly the notion of stable conjugacy for regular elements in $H(F)$ and $H'(F)$ and of transfer between (virtual) representations of $H(F)$ and $H'(F)$. The inner twist $\psi^H$ allows to identify the $L$-groups of $H$ and $H'$ and thus to get an identification $H^1(W_F,Z(\check{H}))=H^1(W_F,Z(\check{H}'))$. We say that two Galoisian characters $\chi$, $\chi'$ of $H(F)$ and $H'(F)$ {\it correspond to each other} if they originate from the same element of $H^1(W_F,Z(\check{H}))$. Galoisian characters are always stable and if $\chi$, $\chi'$ are Galoisian characters of $H(F)$, $H'(F)$ respectively that correspond to each other then they are also transfer of each other.

\begin{theo}\label{theo inner forms}
Let $\Pi$ and $\Pi'$ be stable virtual essentially square-integrable representations of $G(F)$ and $G'(F)$ respectively. Let $\chi$ and $\chi'$ be Galoisian characters of $H(F)$ and $H'(F)$ respectively. Then, if $\Pi$, $\Pi'$ are transfer of each other and $\chi$, $\chi'$ correspond to each other, we have

$$\displaystyle m(\Pi,\chi)=m(\Pi',\chi')$$
\end{theo}

\noindent\ul{Proof}: Let $\Gamma_{\elli}(\overline{H})/stab$ be the set of stable conjugacy classes in $\Gamma_{\elli}(\overline{H})$. It is easy to see that we can equip $\Gamma_{\elli}(\overline{H})/stab$ with a unique topology and a unique measure such that the natural projection $p:\Gamma_{\elli}(\overline{H})\twoheadrightarrow \Gamma_{\elli}(\overline{H})/stab$ is a local isomorphism preserving measures locally. We define $\Gamma_{\elli}(\overline{H}')/stab$ and equip it with a topology and a measure in a similar way. Let $p':\Gamma_{\elli}(\overline{H}')\twoheadrightarrow \Gamma_{\elli}(\overline{H}')/stab$ be the natural projection. Since $\Pi$ and $\Pi'$ are stable and essentially square-integrable, by Theorem \ref{theo formula multiplicity} we have

$$\displaystyle m(\Pi,\chi)=\int_{\Gamma_{\elli}(\overline{H})/stab} \lvert p^{-1}(x)\rvert D^H(x) \Theta_\Pi(x) \chi(x)^{-1} dx$$

\noindent and

$$\displaystyle m(\Pi',\chi')=\int_{\Gamma_{\elli}(\overline{H}')/stab} \lvert {p'}^{-1}(y)\rvert D^{H'}(y)\Theta_{\Pi'}(y) \chi'(y)^{-1} dy$$

Since we can always transfer elliptic regular elements to all inner forms (see \cite{KottSt} \S 10), there is a bijection

$$\displaystyle \Gamma_{\elli}(\overline{H})/stab \simeq \Gamma_{\elli}(\overline{H}')/stab$$

\noindent characterized by: $x\mapsto y$ if and only if $x$ and $y$ are stably conjugate. It is not hard to see that this bijection preserves measures locally and hence globally. Let $x\in H_{\reg}(F)$ and $y\in H'_{\reg}(F)$ be two stably conjugate elements. As $\Theta_\Pi$, $\Theta_{\Pi'}$ on the one hand and $\chi$, $\chi'$ on the other hand are transfer of each other, we have $\Theta_\Pi(x)=\Theta_{\Pi'}(y)$ and $\chi(x)=\chi'(y)$. Moreover, we also have $D^H(x)=D^{H'}(y)$. Therefore, to get the theorem it only remains to show that $\lvert p^{-1}(x)\rvert=\lvert {p'}^{-1}(y)\rvert$. By standard cohomological arguments we have $\lvert p^{-1}(x)\rvert =\left\lvert \ker1^1(F;T,G)\right\rvert$ and $\lvert {p'}^{-1}(y)\rvert=\left\lvert \ker1^1(F;T',G')\right\rvert$ where $T:=G_x$ and $T':=G'_y$. Since $F$ is $p$-adic, by \cite{KottSt} Theorem 1.2 there exist structures of abelian groups on $H^1(F,G)$ and $H^1(F,G')$ such that the natural maps $H^1(F,T)\to H^1(F,G)$ and $H^1(F,T')\to H^1(F,G')$ are morphisms of groups. Moreover, by \cite{KottSt} Lemma 10.2 these are surjective. It follows that $\left\lvert \ker1^1(F;T,G)\right\rvert=\lvert H^1(F,T)\rvert \lvert H^1(F,G)\rvert^{-1}$ and $\left\lvert \ker1^1(F;T',G')\right\rvert=\lvert H^1(F,T')\rvert \lvert H^1(F,G')\rvert^{-1}$. As $T$ and $T'$ are $F$-isomorphic we have $H^1(F,T)\simeq H^1(F,T')$ and by \cite{KottSt} Theorem 1.2 again we have $H^1(F,G)\simeq H^1(F,G')$ (since $G$ and $G'$ have isomorphic $L$-groups). This suffices to conclude that $\left\lvert \ker1^1(F;T,G)\right\rvert=\left\lvert \ker1^1(F;T',G')\right\rvert$ and therefore that $\lvert p^{-1}(x)\rvert=\lvert {p'}^{-1}(y)\rvert$. $\blacksquare$

\subsection{Elliptic twisted Levi subgroups}\label{section elliptic twisted Levi}

\textbf{In this section we assume for simplicity that $H$ is semi-simple and quasi-split.}(All the results presented in this section are still true, with obvious modifications, in general. However, the assumption that $H$ is semi-simple and quasi-split simplifies a lot the proofs and, in any case, we will only need to apply them for such groups.)

\vspace{2mm}

We say that an algebraic subgroup $\mathcal{M}$ of $H$ is a {\it twisted Levi subgroup} if $R_{E/F}\mathcal{M}_E$ is a Levi subgroup of $G$. If $\mathcal{M}$ is a twisted Levi subgroup of $H$, we say that it is {\it elliptic} if $A_{\mathcal{M}}=\{ 1\}$.

\begin{lem}\label{lemma twisted Levi}
Let $M$ be a Levi subgroup of $G$ and set $\mathcal{M}:=M\cap H$. Then, the following assertions are equivalent: 
\begin{enumerate}[(i)]
\item\label{lemma twisted Levi i} $\mathcal{M}$ is an elliptic twisted Levi subgroup of $H$;

\item\label{lemma twisted Levi ii} $\mathcal{M}$ contains an elliptic maximal torus of $H$;

\item\label{lemma twisted Levi iii} $A_M$ is $\theta$-split;

\item\label{lemma twisted Levi iv} $M$ is $\theta$-split and $\mathcal{P}(M)=\mathcal{P}^\theta(M)$.
\end{enumerate}
\end{lem}

\noindent\ul{Proof}:(\ref{lemma twisted Levi i})$\Rightarrow$ (\ref{lemma twisted Levi ii}): By \cite{KnII} p.271, $\mathcal{M}$ contains a maximal torus $T$ such that $A_T=A_{\mathcal{M}}=\{ 1\}$. Thus, $T$ is elliptic and since $\mathcal{M}$ is of the same (absolute) rank as $H$, it is also maximal in $H$. This proves the first implication.

\noindent (\ref{lemma twisted Levi ii})$\Rightarrow$ (\ref{lemma twisted Levi iii}): Assume that $\mathcal{M}$ contains an elliptic maximal torus $\mathcal{T}$ of $H$ and set $T:=R_{E/F}\mathcal{T}_E$. Then we have $A_M\subset A_{T}$ and $A_{T}$ is $\theta$-split (as $(A_{T}^\theta)^0$ is a split torus contained in $\mathcal{T}$ and so is trivial) from which it follows that $A_M$ is also $\theta$-split.

\noindent (\ref{lemma twisted Levi iii})$\Rightarrow$ (\ref{lemma twisted Levi iv}): Assume that $A_M$ is $\theta$-split. Then, $M$ is $\theta$-split since it is the centralizer of $A_M$. Moreover $\theta$ acts on $\mathcal{A}_M$ as $-Id$ thus sending any positive chamber $\mathcal{A}_P^+$ corresponding to $P\in \mathcal{P}(M)$ to its opposite. This shows that $\mathcal{P}(M)=\mathcal{P}^\theta(M)$.

\noindent (\ref{lemma twisted Levi iv})$\Rightarrow$ (\ref{lemma twisted Levi i}): Assume that $M$ is $\theta$-split and $\mathcal{P}(M)=\mathcal{P}^\theta(M)$. In particular $M$ is $\theta$-stable and since $M$ comes by restriction of scalars from a subgroup of $H_E$ (as is any evi subgroup of $G$), it follows that $M=R_{E/F}\mathcal{M}_E$ i.e. $\mathcal{M}$ is a twisted Levi. It only remains to show that $\mathcal{M}$ is elliptic. Assume, by way of contradiction, that it is not the case i.e. $A_{\mathcal{M}}\neq \{1\}$. Then there exists a parabolic $P\in \mathcal{P}(M)$ such that $H_M(A_{\mathcal{M}}(F))$ contains a nonzero element of the closure of the positive chamber associated to $P$. Since this element is $\theta$-fixed and the intersection of the closures of the positive chambers associated to $P$ and $\overline{P}$ (the parabolic subgroup opposite to $P$) is reduced to $\{0\}$ we cannot have $\theta(P)=\overline{P}$ thus contradicting the fact that $\mathcal{P}(M)=\mathcal{P}^\theta(M)$. $\blacksquare$

\vspace{2mm}

Let $P_0$ be a minimal $\theta$-split parabolic subgroup of $G$. We claim that $P_0$ is also a minimal parabolic subgroup of $G$ (hence a Borel subgroup since $G$ is quasi-split).

\noindent\ul{Proof}: Let $B$ be a Borel subgroup of $G$. Since two minimal $\theta$-split parabolic subgroups are always $G(F)$-conjugate (\cite{HW} Proposition 4.9), it suffices to show the existence of $g\in G(F)$ such that $gBg^{-1}$ is $\theta$-split. Over the algebraic closure we have $G_{\overline{F}}\simeq H_{\overline{F}}\times H_{\overline{F}}$ with $\theta$ exchanging the two copies. Since $B$ is in the same class as its opposite Borel this shows the existence of $g_1\in G(\overline{F})$ such that $g_1B_{\overline{F}}g_1^{-1}$ is $\theta$-split. Then, the set $\mathcal{U}:=H_{\overline{F}}g_1B_{\overline{F}}$ is a Zariski open subset of $G_{\overline{F}}$ with the property that for all $g\in \mathcal{U}$ the Borel $gB_{\overline{F}}g^{-1}$ is $\theta$-split. Since $G(F)$ is dense in $G$ for the Zariski topology (\cite{SGA3}, Exp XIV, 6.5, 6.7) we can find $g\in G(F)\cap \mathcal{U}$ and it has the desired property. $\blacksquare$

\vspace{2mm}

\noindent Set $T_0:=P_0\cap \theta(P_0)$, $A_{min}:=A_{T_0}$, $A_0:=A_{min,\theta}$ and denote by $\Delta_{min}$ and $\Delta_0$ the sets of simple roots of $A_{min}$ and $A_0$ in $P_0$. All the parabolic subgroups that we will consider in this section will be standard with respect to $P_0$ (i.e. contain $P_0$) and when we write $P=MU$ for such a parabolic subgroup we always mean that $U$ is the unipotent radical and $M$ the unique Levi component containing $T_0$. We have a natural projection $\Delta_{min}\twoheadrightarrow \Delta_0$ and $\theta$ naturally acts on $R(A_{min},G)$ sending $\Delta_{min}$ to $-\Delta_{min}$. Let $\Delta_-$ be the set of simple roots $\alpha\in \Delta_{min}$ such that $\theta(\alpha)=-\alpha$. It can be identified with a subset of $\Delta_0$ through he projection $\Delta_{min}\twoheadrightarrow \Delta_0$ (i.e. the restriction of this projection to $\Delta_-$ is injective). Let $I\subset \Delta_-$. We will denote by $P_I=M_IU_I$ the unique parabolic subgroup containing $P_0$ such that the set of simple roots of $A_{min}$ in $M_I\cap P_0$ is precisely $\Delta_{min}-I$. Then $P_I$ is $\theta$-split and we have $M_I=P_I\cap \theta(P_I)$. Moreover, as $\mathcal{A}^*_{M_I}$ is generated by the restrictions of the roots in $I$, $\theta$ acts as $-Id$ on this space showing that $A_{M_I}$ is $\theta$-split and thus by point (\ref{lemma twisted Levi ii}) of the previous lemma that $\mathcal{M}_I:=M_I\cap H$ is an elliptic twisted Levi subgroup of $H$. Let $H_{ab}$ be the quotient of $H(F)$ by the common kernel of all the Galoisian characters $\chi$ of $H(F)$. It is an abelian group and for all $I\subset \Delta_-$ we will denote by $\mathcal{M}_{I,ab}$ the image of $\mathcal{M}_I(F)$ in $H_{ab}$. 

\vspace{2mm}

Define $\underline{\mathcal{C}}$ as the set of pairs $(\mathcal{M},P)$ with $\mathcal{M}$ an elliptic twisted Levi subgroup of $H$ and $P\in \mathcal{P}(M)$ where $M:=R_{E/F}\mathcal{M}_E$ and $\underline{\mathcal{D}}$ as the set of triples $(T,\mathcal{M},P)$ with $(\mathcal{M},P)\in \underline{\mathcal{C}}$ and $T\subset \mathcal{M}$ an elliptic maximal torus. We let $\mathcal{C}$ and $\mathcal{D}$ denote the $H(F)$-conjugacy classes in $\underline{\mathcal{C}}$ and $\underline{\mathcal{D}}$ respectively. We say that two pairs $(\mathcal{M},P)$, $(\mathcal{M}',P')\in \mathcal{C}$ are {\it stably conjugate}, and we will write $(\mathcal{M},P)\sim_{stab}(\mathcal{M}',P')$, if there exists $h\in H(\overline{F})$ such that $h\mathcal{M}_{\overline{F}}h^{-1}=\mathcal{M}'_{\overline{F}}$ and $hP_{\overline{F}}h^{-1}=P'_{\overline{F}}$. Write $\mathcal{C}/stab$ for the set of stable conjugacy classes in $\mathcal{C}$. We let $\mathcal{T}_{\elli}(H)$ be a set of representatives of the $H(F)$-conjugacy classes of elliptic maximal tori in $H$. Finally for all $P=MU\supset P_0$ and all $T\in \mathcal{T}_{\elli}(H)$, we define 

$$\displaystyle \Gamma_M(T):=\left\{\gamma\in G(F); \gamma^{-1}T\gamma\subset M \right\}/M(F)$$

\begin{prop}\label{prop twisted Levi}
\begin{enumerate}[(i)]
\item\label{prop twisted Levi i} For all $T\in \mathcal{T}_{\elli}(H)$, all $P=MU\supseteq P_0$ and all $\gamma\in \Gamma_M(T)$ we have $(T,\gamma M\gamma^{-1}\cap H,\gamma P\gamma^{-1})\in \mathcal{D}$;

\item\label{prop twisted Levi ii} The map

$$\displaystyle \bigsqcup_{T\in \mathcal{T}_{\elli}(H)} \bigsqcup_{P_0\subseteq P=MU} \Gamma_M(T)\to\mathcal{D}$$
$$\displaystyle \gamma\in \Gamma_M(T)\mapsto (T,\gamma M\gamma^{-1}\cap H,\gamma P\gamma^{-1})$$

\noindent (which is well-defined by (\ref{prop twisted Levi i})) is surjective and the fiber over $(T,\mathcal{M},P)\in \mathcal{D}$ is of cardinality $\lvert W(H,T)\rvert \lvert W(\mathcal{M},T)\rvert^{-1}$.

\item\label{prop twisted Levi iii} For all $(\mathcal{M},P)\in \mathcal{C}$ the fiber of the map $\mathcal{C}\to \mathcal{C}/stab$ containing $(\mathcal{M},P)$ is of cardinality $\lvert \ker1^1(F;\mathcal{M},H)\rvert$.

\item\label{prop twisted Levi iv} For all $(\mathcal{M},P)$, $(\mathcal{M}',P')\in \mathcal{C}$ we have $(\mathcal{M},P)\sim_{stab}(\mathcal{M}',P')$ if and only if $P$ and $P'$ are in the same class and the map $I\subseteq \Delta_-\mapsto (\mathcal{M}_I,P_I)\in \mathcal{C}/stab$ is a bijection.

\item\label{prop twisted Levi v} Let $(\mathcal{M},P),(\mathcal{M}',P')\in \mathcal{C}$ be such that $(\mathcal{M},P)\sim_{stab}(\mathcal{M}',P')$. Then, for every Galoisian character $\chi$ of $H(F)$ we have $\chi_{\mid \mathcal{M}}=1$ if and only if $\chi_{\mid \mathcal{M}'}=1$, where we have denoted by $\chi_{\mid \mathcal{M}}$ and $\chi_{\mid \mathcal{M}'}$ the restrictions of $\chi$ to $\mathcal{M}(F)$ and $\mathcal{M}'(F)$ respectively.

\item\label{prop twisted Levi vi} We have the following identity in $\mathcal{R}(H_{ab})$:

$$\displaystyle \sum_{I\subset \Delta_-} (-1)^{\lvert \Delta_- -I\rvert} \lvert \ker1^1(F;\mathcal{M}_I,H)\rvert Ind_{\mathcal{M}_{I,ab}}^{H_{ab}}(\mathbf{1})=\omega_{H,E}$$

\noindent where $\omega_{H,E}:H_{ab}\to \{\pm 1 \}$ is Prasad's character (see \S \ref{section galoisian chars}).
\end{enumerate}
\end{prop}

\noindent\ul{Proof}:
\begin{enumerate}[(i)]
\item For all $T\in \mathcal{T}_{\elli}(H)$, all $P=MU\supset P_0$ and all $\gamma\in \Gamma_M(T)$, the subgroup $\gamma M\gamma^{-1}\cap H$ contains an elliptic maximal torus of $H$ (namely $T$) and thus by Lemma \ref{lemma twisted Levi} is an elliptic twisted Levi subgroup. This shows that $(T,\gamma M\gamma^{-1}\cap H,\gamma P\gamma^{-1})\in \mathcal{D}$.

\item Let $(T,\mathcal{M},P)\in \mathcal{D}$. Up to conjugation we may assume that $T\in \mathcal{T}_{\elli}(H)$. As $P_0$ is a minimal parabolic subgroup of $G$, there exist $\gamma\in G(F)$ such that $\gamma^{-1}P\gamma\supset P_0$ and $\gamma^{-1} R_{E/F}\mathcal{M}_E\gamma\supset T_0$. This shows the surjectivity. The claim about the cardinality of the fibers is a consequence of the two following facts:

\vspace{3mm}

\begin{num}
\item\label{eq 5.4.1} Let $T,T'\in \mathcal{T}_{\elli}(H)$, $P=MU\supset P_0$, $P'=M'U'\supset P_0$, $\gamma\in \Gamma_M(T)$ and $\gamma'\in \Gamma_{M'}(T')$. Then the two triples $(T,\gamma M\gamma^{-1}\cap H,\gamma P\gamma^{-1})$ and $(T',\gamma M'\gamma^{-1}\cap H,\gamma P'\gamma^{-1})$ are $H(F)$-conjugate if and only if $T=T'$, $M=M'$, $P=P'$ and $\gamma'\in \No_{H(F)}(T)\gamma M(F)$.
\end{num}

\vspace{3mm}

\begin{num}
\item\label{eq 5.4.2} Let $T\in \mathcal{T}_{\elli}(H)$, $P=MU\supset P_0$ and $\gamma\in \Gamma_M(T)$. Then, the image of the map
$$\displaystyle \No_{H(F)}(T)\to \Gamma_M(T)$$
$$\displaystyle h\mapsto h\gamma$$
is of cardinality $\lvert W(H,T)\rvert \lvert W(\mathcal{M},T)\rvert^{-1}$ where $\mathcal{M}:= \gamma M\gamma^{-1}\cap H$.
\end{num}

\vspace{3mm}

\noindent Proof of \ref{eq 5.4.1}: Assume that $(T,\gamma M\gamma^{-1}\cap H,\gamma P\gamma^{-1})$ and $(T',\gamma M'\gamma^{-1}\cap H,\gamma P'\gamma^{-1})$ are $H(F)$-conjugate. By definition of $\mathcal{T}_{\elli}(H)$ we have $T=T'$ and since $P$ and $P'$ are both standard with respect to $P_0$ we also have $P=P'$ and $M=M'$. Thus there exists $h\in H(F)$ such that

$$\displaystyle h(T,\mathcal{M}, \gamma P\gamma^{-1})h^{-1}=(T,\mathcal{M}', \gamma P{\gamma'}^{-1})$$

\noindent where $\mathcal{M}:=\gamma M\gamma^{-1}\cap H$ and $\mathcal{M}':=\gamma' M{\gamma'}^{-1}\cap H$. This equality immediately implies that $h\in \No_{H(F)}(T)$. Moreover, since $\gamma M\gamma^{-1}=R_{E/F}\mathcal{M}_E$ and $\gamma' M {\gamma'}^{-1}=R_{E/F}\mathcal{M}'_E$ we also have $h\gamma M \gamma^{-1}h^{-1}=\gamma' M{\gamma'}^{-1}$ and $h\gamma P \gamma^{-1}h^{-1}=\gamma' P{\gamma'}^{-1}$ proving that ${\gamma'}^{-1}h\gamma$ normalizes both $M$ and $P$ i.e. $\gamma'\in h\gamma M(F)$ and hence $\gamma'\in \No_{H(F)}(T)\gamma M(F)$. This proves one direction of the claim the other being obvious.

\vspace{2mm}

\noindent Proof of \ref{eq 5.4.2}: As $\lvert W(H,T)\rvert \lvert W(\mathcal{M},T)\rvert^{-1}$ is the cardinality of the quotient

$$\No_{H(F)}(T)/\No_{\mathcal{M}(F)}(T)$$

\noindent it suffices to show that for all $h,h'\in \No_{H(F)}(T)$ we have $h\gamma=h'\gamma$ in $\Gamma_M(T)$ if and only if ${h'}^{-1}h\in \mathcal{M}(F)$. But as $\mathcal{M}=\gamma M\gamma^{-1}\cap H$ this immediately follows from the definition of $\Gamma_M(T)$.

\item This follows from a standard cohomological argument by noticing that $\mathcal{M}$ is the normalizer of the pair $(\mathcal{M},P)$ in $H$.

\item If $(\mathcal{M},P)\sim_{stab}(\mathcal{M}',P')$ then in particular $P$ and $P'$ are conjugate in $G(\overline{F})$ and thus are in the same class (i.e. are $G(F)$-conjugate). Conversely, assume that $P$ and $P'$ are in the same class and set $M:=R_{E/F}\mathcal{M}_E$, $M':=R_{E/F}\mathcal{M}'_E$. Then, there exists $g\in G(F)$ such that $gPg^{-1}=P'$ and $gMg^{-1}=M'$. Let $\overline{P}$ and $\overline{P}'$ be the parabolic subgroups opposite to $P$ and $P'$ with respect to $M$ and $M'$ respectively. Then $g\overline{P}g^{-1}=\overline{P}'$ and since $\theta(P)=\overline{P}$, $\theta(P')=\overline{P}'$, $\theta(M)=M$ and $\theta(M')=M'$ we also have

$$\displaystyle \theta(g)\overline{P}\theta(g)^{-1}=\theta(gPg^{-1})=\overline{P}'$$

\noindent and

$$\displaystyle \theta(g)M\theta(g)^{-1}=\theta(gMg^{-1})=M'$$

\noindent Therefore $g^{-1}\theta(g)$ normalizes both $M$ and $\overline{P}$ and thus $g^{-1}\theta(g)\in M(F)$. Since the map

$$M(\overline{F})\to \{m\in M(\overline{F}); \theta(m)=m^{-1} \}$$

$$m\mapsto m^{-1}\theta(m)$$

\noindent is surjective, there exist $m\in M(\overline{F})$ and $h\in H(\overline{F})$ such that $g=hm$. We have $hP_{\overline{F}}h^{-1}=P'_{\overline{F}}$ and $h \mathcal{M}_{\overline{F}}h^{-1}=\mathcal{M}'_{\overline{F}}$ showing that $(\mathcal{M},P)$ and $(\mathcal{M}',P')$ are stably conjugate. This shows the first part of (\ref{prop twisted Levi iv}). Notice that for all $P\in \mathcal{P}^\theta(M_0)$ with $P\supset P_0$ the torus $A_M$, where $M:=P\cap \theta(P)$, is $\theta$-split if and only if $P=P_I$ for some $I\subset \Delta_-$. Hence, the second part of (\ref{prop twisted Levi iv}) follows from the first and Lemma \ref{lemma twisted Levi} since every class of parabolic subgroups contains a unique element which is standard with respect to $P_0$.

\item We need to show that $\mathcal{M}_{ab}=\mathcal{M}'_{ab}$. Let $H_{sc}$ be the simply connected cover of $H$ and $\pi_1(H)$ be the kernel of the projection $H_{sc}\to H$. Then we have $H_{ab}=H^1(F,\pi_1(H))$. Let $\mathcal{M}_{sc}$ and $\mathcal{M}'_{sc}$ denote the inverse images of $\mathcal{M}$ and $\mathcal{M}'$ in $H_{sc}$. From the short exact sequences $1\to \pi_1(H)\to \mathcal{M}_{sc}\to\mathcal{M}\to 1$ and $1\to \pi_1(H)\to \mathcal{M}'_{sc}\to\mathcal{M}'\to 1$, we get exact sequences $1\to \mathcal{M}_{ab}\to H_{ab}\to H^1(F,\mathcal{M}_{sc})$ and $1\to \mathcal{M}'_{ab}\to H_{ab}\to H^1(F,\mathcal{M}'_{sc})$. Thus, we need to show that $\Ker\left(H_{ab}\to H^1(F,\mathcal{M}_{sc})\right)=\Ker\left(H_{ab}\to H^1(F,\mathcal{M}'_{sc}) \right)$. By hypothesis, there exists $h\in H(\overline{F})$ such that $h\mathcal{M}_{\overline{F}}h^{-1}=\mathcal{M}'_{\overline{F}}$ and $hP_{\overline{F}}h^{-1}=P'_{\overline{F}}$. Set $M:=R_{E/F} \mathcal{M}_E$ and $M':=R_{E/F} \mathcal{M}'_E$. Then, we also have $hM_{\overline{F}}h^{-1}=M'_{\overline{F}}$ and it follows that for all $\sigma\in \Gamma_F$ we have $h{}^\sigma h^{-1}\in H(\overline{F})\cap M'(\overline{F})=\mathcal{M}'(\overline{F})$. Choose $h_{sc}\in H_{sc}(\overline{F})$ which lifts $h$. Then, there is a bijection $\iota: H^1(F,\mathcal{M}_{sc})\simeq H^1(F,\mathcal{M}'_{sc})$ given at the level of cocycles by

$$\displaystyle c\mapsto \left(\sigma\in \Gamma_F\mapsto h_{sc}c(\sigma){}^\sigma h_{sc}^{-1} \right)$$

\noindent Let $c_0$ be the $1$-cocycle $\sigma\in \Gamma_F\mapsto h_{sc}{}^\sigma h_{sc}^{-1}\in \mathcal{M}'_{sc}(\overline{F})$ and denote by $[c_0]$ its class in $H^1(F,\mathcal{M}'_{sc})$. Then $\iota-[c_0]: H^1(F,\mathcal{M}_{sc})\simeq H^1(F,\mathcal{M}'_{sc})$ is a bijection of pointed sets (and even an isomorphism of abelian groups) making the following square commute

$$\xymatrix{H_{ab} \ar@{->}[r] \ar@{=}[d] & H^1(F,\mathcal{M}_{sc}) \ar@{->}[d]^{\wr} \\ H_{ab} \ar@{->}[r] & H^1(F,\mathcal{M}'_{sc})}$$ 

\noindent This immediately implies that the kernels of the upper and bottom arrows are identical.

\item Let $I\subset \Delta_-$ and denote by $\mathcal{M}_{I,sc}$ the inverse image of $\mathcal{M}_I$ in $H_{sc}$. Then, from the short exact sequence $1\to \pi_1(H)\to \mathcal{M}_{I,sc}\to\mathcal{M}_I\to 1$ we get an exact sequence

$$\displaystyle 1\to \mathcal{M}_{I,ab}\to H_{ab}\to H^1(F,\mathcal{M}_{I,sc})\to H^1(F,\mathcal{M}_I)\to H^2(F,\pi_1(H))$$

\noindent By \cite{KnI} the natural connecting map $H^1(F,H)\to H^2(F,\pi_1(H))$ is an isomorphism and it follows that the previous exact sequence can be rewritten as

$$\displaystyle 1\to \mathcal{M}_{I,ab}\to H_{ab}\to H^1(F,\mathcal{M}_{I,sc})\to \ker1^1(F;\mathcal{M}_I,H)\to 1$$

\noindent From this exact sequence, we deduce the following equality in $\mathcal{R}(H_{ab})$:

\begin{align}\label{eq 5.4.3}
\displaystyle \lvert \ker1^1(F;\mathcal{M}_I,H)\rvert Ind_{\mathcal{M}_{I,ab}}^{H_{ab}}(\mathbf{1})= \Res_I^{H_{ab}}Ind_1^{I}(\mathbf{1})
\end{align}

\noindent where $\Res_I^{H_{ab}}$ denotes the restriction functor with respect to the morphism $H_{ab}\to H^1(F,\mathcal{M}_{I,sc})$ and $Ind_1^I$ denotes the induction functor with respect to the morphism $1\to H^1(F,\mathcal{M}_{I,sc})$. Moreover the morphism $H^1(F,\mathcal{M}_{\Delta_-,sc})\to H^1(F,\mathcal{M}_{I,sc})$ induced by the inclusion $\mathcal{M}_{\Delta_-}\subset \mathcal{M}_I$, makes the following square commute

$$\xymatrix{H_{ab} \ar@{->}[r] \ar@{=}[d] & H^1(F,\mathcal{M}_{\Delta_-,sc}) \ar@{->}[d] \\ H_{ab} \ar@{->}[r] & H^1(F,\mathcal{M}_{I,sc})}$$

\noindent Hence, we have a factorization $\Res^{H_{ab}}_I=\Res^{H_{ab}}_{\Delta_-}\circ \Res_I^{\Delta_-}$ where $\Res_I^{\Delta_-}$ denotes the restriction functor with respect to the morphism $H^1(F,\mathcal{M}_{\Delta_-,sc})\to H^1(F,\mathcal{M}_{I,sc})$. Combining this with \ref{eq 5.4.3}, we get the identity

\begin{align}\label{eq 5.4.4}
\displaystyle  \sum_{I\subset \Delta_-} (-1)^{\lvert \Delta_- -I\rvert} \lvert \ker1^1(F;\mathcal{M}_I,H)\rvert Ind_{\mathcal{M}_{I,ab}}^{H_{ab}}(\mathbf{1})=\Res_{\Delta_-}^{H_{ab}}\left( \sum_{I\subset \Delta_-} (-1)^{\lvert \Delta_- -I\rvert}\Res_I^{\Delta_-}Ind_1^{I}(\mathbf{1})\right)
\end{align}

\noindent To continue, we need to compute the groups $H^1(F,\mathcal{M}_{I,sc})$ and the morphisms 

$$H^1(F,\mathcal{M}_{\Delta_-,sc})\to H^1(F,\mathcal{M}_{I,sc})$$

\noindent explicitly for all $I\subset \Delta_-$. Let $T_{0,sc}$ and $P_{0,sc}$ denote the inverse image of $T_0$ and $P_0$ in $G_{sc}:=R_{E/F}H_{sc,E}$. Since $T_{0,sc}$ is $\theta$-stable, there exists a maximal torus $\mathcal{T}_{0,sc}$ of $H$ such that $T_{0,sc}=R_{E/F}\mathcal{T}_{0,sc,E}$. Moreover, there exists a Borel subgroup $\mathcal{P}_{0,sc}$ of $H_{sc,E}$ such that $P_{0,sc}=R_{E/F} \mathcal{P}_{0,sc}$. In what follows, we fix an algebraic closure $\overline{F}$ of $F$ containing $E$ and we set $\Gamma_E:=\Gal(\overline{F}/E)$. Let $\Delta_{min,\overline{F}}$ be the set of simple roots of $\mathcal{T}_{0,sc,\overline{F}}$ in $\mathcal{P}_{0,sc,\overline{F}}$. It is a subset of $X^*_{\overline{F}}(\mathcal{T}_{0,sc})$ which is $\Gamma_E$-stable (as $\mathcal{P}_{0,sc}$ is defined over $E$) and we have a natural surjection $\Delta_{min,\overline{F}}\twoheadrightarrow \Delta_{min}$ (obtained by restriction to the maximal split subtorus of $\mathcal{T}_{0,sc,E}$) whose fibers are precisely the $\Gamma_E$-orbits in $\Delta_{min,\overline{F}}$. For all $\beta\in \Delta_{min,\overline{F}}$, we will denote by $\varpi_\beta\in X^*_{\overline{F}}(\mathcal{T}_{0,sc})$ the corresponding weight and for all $\alpha\in \Delta_{min}$ we define

$$\displaystyle \varpi_\alpha:=\sum_{\beta\in \Delta_{min,\overline{F}};\beta\mapsto \alpha}\varpi_\beta$$

\noindent where the sum is over the set of simple roots $\beta\in \Delta_{min,\overline{F}}$ mapping to $\alpha$ through the projection $\Delta_{min,\overline{F}}\twoheadrightarrow \Delta_{min}$. We always have $\varpi_\alpha\in X^*_E(\mathcal{T}_{0,sc})$ but we warn the reader that in general $\varpi_\alpha$ is NOT the weight associated to the simple root $\alpha$ in the usual sense (although it is proportional to it). Let $I\subset \Delta_-$. For all $\alpha\in I$ the character $\varpi_\alpha$ extends to $\mathcal{M}_{I,sc,\overline{F}}$ and is defined over $E$ and thus gives rise to a character $R_{E/F}\mathcal{M}_{I,sc,E}\to R_{E/F}\mathbf{G}_{m,E}$. Since $\theta(\varpi_\alpha)=-\varpi_\alpha$, this last character induces a morphism $\mathcal{M}_{I,sc}\to \Ker N_{E/F}$ that we will also denote by $\varpi_\alpha$. Consider the torus $T_I:=\left(\Ker N_{E/F} \right)^I$ and the morphism

$$\displaystyle \kappa_I:=(\varpi_\alpha)_{\alpha\in I}:\mathcal{M}_{I,sc}\to T_I$$

\noindent Then we claim that

\vspace{3mm}

\begin{num}
\item\label{eq 5.4.5} The induced map $H^1(\kappa_I): H^1(F,\mathcal{M}_{I,sc})\to H^1(F,T_I)$ is an isomorphism.
\end{num}

\vspace{3mm}

\noindent Let $\mathcal{M}_{I,sc,der}$ be the derived subgroup of $\mathcal{M}_{I,sc}$ and set 

$$\displaystyle T'_I:=\mathcal{M}_{I,sc}/\mathcal{M}_{I,sc,der}$$

\noindent Then, we have an exact sequence

$$\displaystyle H^1(F,\mathcal{M}_{I,sc,der})\to H^1(F,\mathcal{M}_{I,sc})\to H^1(F,T'_I)$$

\noindent and $H^1(F,\mathcal{M}_{I,sc,der})$ is trivial by \cite{KnI} since $\mathcal{M}_{I,sc,der}$ is simply connected. Moreover the morphism $H^1(F,\mathcal{M}_{I,sc})\to H^1(F,T'_I)$ is surjective. Indeed, if $T\subset \mathcal{M}_{I,sc}$ is a maximal elliptic torus (which exists by \cite{KnII} p.271) then the kernel of the projection $T\to T'_I$ is a maximal anisotropic torus of $\mathcal{M}_{I,sc,der}$ and thus by Tate-Nakayama duality its $H^2$ vanishes and the morphism $H^1(F,T)\to H^1(F,T'_I)$, which obviously factorizes through $H^1(F,\mathcal{M}_{I,sc})$, is surjective. Therefore, the morphism $H^1(F,\mathcal{M}_{I,sc})\to H^1(F,T'_I)$ is an isomorphism and it only remains to show that the natural map $H^1(F,T'_I)\to H^1(F,T_I)$ is also an isomorphism. Let $T''_I$ be the kernel of the projection $T'_I\to T_I$. Then $T''_I$ is connected since it is a quotient of the common kernel of all the $\varpi_\alpha$'s, $\alpha\in I$, which is a connected group (this follows from the fact that $\{\varpi_\alpha;\alpha\in I \}$ generates $X^*_E(\mathcal{M}_{I,sc})$). Therefore, $T''_I$ is an anisotropic torus (since $T'_I$ is) and by Tate-Nakayama duality again we just need to prove the injectivity of $H^1(F,T'_I)\to H^1(F,T_I)$ or, equivalently, that the map $H^1(F,T''_I)\to H^1(F,T'_I)$ has trivial image. We have norm maps $N: R_{E/F}T'_{I,E}\to T'_I$ and $N: R_{E/F} T''_{I,E}\to T''_I$ giving rise to a commuting square

$$\displaystyle \xymatrix{H^1(F,T''_I) \ar@{->}[r] & H^1(F,T'_I) \\ H^1(F,R_{E/F}T''_{I,E}) \ar@{->}[r] \ar@{->}[u]^N & H^1(F,R_{E/F} T'_{I,E})\ar@{->}[u]^N}$$

\noindent Since $T_{I,E}\simeq\mathbf{G}_{m,E}^I$ is the maximal split torus quotient of $\mathcal{M}_{I,sc,E}$, the torus $R_{E/F}T''_{I,E}$ is anisotropic and therefore so is the kernel of the norm map $R_{E/F} T''_{I,E}\to T''_I$ (which is automatically connected). Hence, by Tate-Nakayama again, the map $H^1(F,R_{E/F}T''_{I,E})\to H^1(F,T''_I)$ is surjective and it follows, by the above commuting square, that to conclude we only need to show that $H^1(F,R_{E/F}T'_{I,E})$ is trivial. By an argument similar to what we have done before, the map $H^1(F,M_{I,sc})\to H^1(F,R_{E/F} T'_{I,E})$, where $M_{I,sc}:=R_{E/F}\mathcal{M}_{I,sc,E}$, is surjective. Since $M_{I,sc}$ is a Levi subgroup of $G_{sc}$, the map $H^1(F,M_{I,sc})\to H^1(F,G_{sc})$ has trivial kernel and by \cite{KnI} it follows that $H^1(F,M_{I,sc})=1$. Hence, $H^1(F,R_{E/F}T'_{I,E})=1$ also and this ends the proof of \ref{eq 5.4.5}.

\vspace{2mm}

\noindent By \ref{eq 5.4.5}, we have isomorphisms

\begin{align}\label{eq 5.4.6}
\displaystyle H^1(F,\mathcal{M}_{I,sc})\simeq H^1(F,T_I)\simeq \left(\mathbf{Z}/2\mathbf{Z}\right)^I
\end{align}

\noindent for all $I\subset \Delta_-$ such that the maps $H^1(F,\mathcal{M}_{\Delta_-,sc})\to H^1(F,\mathcal{M}_{I,sc})$ correspond to the natural projections

$$\displaystyle \left(\mathbf{Z}/2\mathbf{Z}\right)^{\Delta_-}\to \left(\mathbf{Z}/2\mathbf{Z}\right)^I$$
$$\displaystyle (e_\alpha)_{\alpha\in \Delta_-}\mapsto (e_\alpha)_{\alpha\in I}$$

\noindent We can now compute the right hand side of \ref{eq 5.4.4}. For all $I\subset \Delta_-$, let $A_I$ denote the kernel of the projection $H^1(F,\mathcal{M}_{\Delta_-,sc})\to H^1(F,\mathcal{M}_{I,sc})$. In particular we have $A_\emptyset=H^1(F,\mathcal{M}_{\Delta_-,sc})\simeq \left(\mathbf{Z}/2\mathbf{Z}\right)^{\Delta_-}$ and $A_I$ is the subgroup $\left(\mathbf{Z}/2\mathbf{Z}\right)^{\Delta_-\backslash I}$. For all $I\subset \Delta_-$ we have

$$\displaystyle \Res_I^{\Delta_-}Ind_1^I(\mathbf{1})=\sum_{\chi\in \Irr(A_\emptyset); \chi_{\mid A_I}=1}\chi$$

\noindent where $\Irr(A_\emptyset)$ denotes the set of characters of $A_\emptyset$. Thus the sum

$$\displaystyle \sum_{I\subset \Delta_-} (-1)^{\lvert \Delta_- -I\rvert}\Res_I^{\Delta_-}Ind_1^{I}(\mathbf{1})$$

\noindent is equal to

$$\displaystyle \sum_{I\subset \Delta_-} (-1)^{\lvert \Delta_--I\rvert} \sum_{\chi\in \Irr(A_\emptyset); \chi_{\mid A_I}=1}\chi=\sum_{\chi\in \Irr(A_\emptyset)}\left(\sum_{I\subset \Delta_-; \chi_{\mid A_I}=1} (-1)^{\lvert \Delta_--I\rvert}\right)\chi$$

\noindent Let $\chi\in \Irr(A_\emptyset)$. By the above description of the subgroups $A_I\subset A_\emptyset$, we see that $A_I+A_J=A_{I\cap J}$ for all $I,J\subset \Delta_-$. Hence there exists a smallest subset $I_\chi\subset \Delta_-$ such that $\chi_{\mid A_{I_\chi}}=1$ and we have

$$\displaystyle \sum_{I\subset \Delta_-; \chi_{\mid A_I}=1} (-1)^{\lvert \Delta_--I\rvert}=\sum_{I_\chi\subset I\subset \Delta_-} (-1)^{\lvert \Delta_--I\rvert}$$

\noindent This sum is zero unless $I_\chi=\Delta_-$ in which case it is equal to $1$. Using again the explicit determination of the subgroups $A_I\subset A_\emptyset$, it is easy to see that there is only one character $\chi\in \Irr(A_\emptyset)$ with $I_\chi=\Delta_-$, namely the character $\omega$ defined (via the isomorphism \ref{eq 5.4.6}) by

$$\displaystyle \left(\mathbf{Z}/2\mathbf{Z}\right)^{\Delta_-}\to\{\pm 1 \}$$

$$\displaystyle (e_\alpha)_{\alpha\in \Delta_-}\mapsto (-1)^{\sum_{\alpha\in \Delta_-}e_\alpha}$$

\noindent Therefore by \ref{eq 5.4.4}, we get

$$\displaystyle \sum_{I\subset \Delta_-} (-1)^{\lvert \Delta_- -I\rvert} \lvert \ker1^1(F;\mathcal{M}_I,H)\rvert Ind_{\mathcal{M}_{I,ab}}^{H_{ab}}(\mathbf{1})=\Res_{\Delta_-}^{H_{ab}}(\omega)$$

\noindent and it only remains to show that $\Res_{\Delta_-}^{H_{ab}}(\omega)=\omega_{H,E}$. Set

$$\displaystyle \rho:=\sum_{\alpha\in \Delta_{min}}\varpi_\alpha=\sum_{\beta\in \Delta_{min,\overline{F}}}\varpi_\beta,\;\;\; \rho_1:=\sum_{\alpha\in \Delta_-}\varpi_\alpha \mbox{ and } \rho_2:=\sum_{\alpha\in \Delta_{min}\backslash \Delta_-}\varpi_\alpha$$

\noindent Then, by restriction $\rho$, $\rho_1$ and $\rho_2$ define three morphisms $\pi_1(H)\to \Ker N_{E/F}$. By definition, $\omega_{H,E}$ and $\Res_{\Delta_-}^{H_{ab}}(\omega)$ are the morphisms

$$H_{ab}=H^1(F,\pi_1(H))\to H^1(F,\Ker N_{E/F})\simeq \{\pm 1 \}$$

\noindent induced by $\rho$ and $\rho_1$ respectively. Thus it suffices to show that the morphism

$$\displaystyle H^1(\rho_2):H^1(F,\pi_1(H))\to H^1(F,\Ker N_{E/F})$$

\noindent induced by $\rho_2$ is trivial. By definition of $\Delta_-$ we can find a subset $S\subset \Delta_{min}\backslash \Delta_-$ such that $\Delta_{min}\backslash \Delta_-=S\sqcup -\theta(S)$ (disjoint union). For all $\alpha\in S$, the character $\varpi_\alpha$ induces a morphism $\pi_1(H)\to R_{E/F}\mathbf{G}_{m,E}\simeq R_{E/F} (\Ker N_{E/F})_E$. By the decomposition $\Delta_{min}\backslash \Delta_-=S\sqcup -\theta(S)$, we see that $H^1(\rho_2)$ is the composition of

$$\displaystyle \sum_{\alpha\in S} H^1(\varpi_\alpha):H^1(F,\pi_1(H))\to H^1(F,R_{E/F} (\Ker N_{E/F})_E)$$

\noindent with the norm map $H^1(F,R_{E/F} (\Ker N_{E/F})_E)\to H^1(F,\Ker N_{E/F})$. By Hilbert 90, we have

$$H^1(F,R_{E/F} (\Ker N_{E/F})_E)=1$$

\noindent and thus $H^1(\rho_2)=0$. This ends the proof of the proposition. $\blacksquare$
\end{enumerate}

\subsection{Reminder on the Steinberg representation}\label{section Steinberg}

Fix a minimal parabolic subgroup $P_0$ of $G$ with Levi decomposition $P_0=M_0U_0$. Then, the Steinberg representation of $G(F)$ is by definition the following virtual representation

$$\displaystyle \St(G):=\sum_{P_0\subset P=MU} (-1)^{a_M-a_{M_0}} i_P^G(\delta_P^{1/2})$$

\noindent where $i_P^G$ denotes the functor of normalized parabolic induction. It follows from \cite{Cass} that $\St(G)$ is in fact a true representation of $G(F)$ which is moreover irreducible and square-integrable. Obviously, the Steinberg representation has trivial central character. Moreover, if $G_{ad}$ denotes the adjoint group of $G$ then $\St(G)$ is the pullback of $\St(G_{ad})$ by the projection $G(F)\to G_{ad}(F)$. For all $x\in G_{\reg}(F)$ and all $P=MU\supset P_0$, let us set

$$\displaystyle \Gamma_M(x):=\{\gamma\in G(F);\; \gamma^{-1}x\gamma\in M(F) \}/M(F)$$

\noindent Then, the character $\Theta_{\St(G)}$ of $\St(G)$ is given by the following formula (\cite{HC} Theorem 30)

\begin{align}\label{eq 5.5.1}
\displaystyle D^G(x)^{1/2}\Theta_{\St(G)}(x)=\sum_{P_0\subset P=MU} (-1)^{a_M-a_{M_0}}\sum_{\gamma\in \Gamma_M(x)}D^M(\gamma^{-1}x\gamma)^{1/2}\delta_P(\gamma^{-1}x\gamma)^{1/2}
\end{align}

\noindent for all $x\in G_{\reg}(F)$. In particular, we have

\begin{align}\label{eq 5.5.2}
\displaystyle \Theta_{\St(G)}(x)=(-1)^{a_G-a_{M_0}}
\end{align}

\noindent for all $x\in G(F)_{\elli}$

\vspace{2mm}

The representation $\St(G)$ is stable: this follows from the fact that parabolic induction sends stable distributions to stable distributions. Let $H'$ be another connected reductive group over $F$ and $\psi^H:H_{\overline{F}}\simeq H'_{\overline{F}}$ be an inner twisting. Following notations of \S \ref{section first application} we define $G':=R_{E/F}H'_E$ and the inner twisting $\psi^G:G_{\overline{F}}\simeq G'_{\overline{F}}$. Let $P_0'=M_0'N_0'$ be a minimal parabolic subgroup of $G'$. Since transfer is compatible with parabolic induction we have that $(-1)^{a_{M_0}}\St(G)$ and $(-1)^{a_{M'_0}}\St(G')$ are transfer of each other. Moreover, in our situation we have $(-1)^{a_{M_0}-a_{M'_0}}=1$. Indeed, by the main result of \cite{Kottsign} we have $(-1)^{a_{M_0}-a_{M'_0}}=e(G)e(G')$ where $e(G)$ and $e(G')$ are the so-called {\it Kottwitz signs} of $G$ and $G'$ respectively and it follows from points (4) and (5) of the Corollary of {\it loc.cit.} that $e(G)=e(H)^2=1$ and $e(G')=e(H')^2=1$. Thus,we have that

\vspace{3mm}

\begin{num}
\item\label{eq 5.5.3} $\St(G)$ and $\St(G')$ are transfer of each other.
\end{num}

\subsection{Harish-Chandra's orthogonality relations for discrete series}\label{section orthogonal relations}

Let $\pi$ and $\sigma$ be essentially square-integrable representations of $H(F)$ with central characters coinciding on $A_H(F)$. Then, we have the following orthogonality relation between the characters of $\pi$ and $\sigma^\vee$ (the smooth contragredient of $\sigma$) which is due to Harish-Chandra (see \cite{Clo} Theorem 3):

$$\displaystyle \int_{\Gamma_{\elli}(\overline{H})}D^H(x)\Theta_{\pi}(x)\Theta_{\sigma^\vee}(x)dx=\left\{
    \begin{array}{ll}
        1 & \mbox{if } \pi\simeq \sigma \\
        0 & \mbox{otherwise}
    \end{array}
\right.
$$

\noindent where the measure on $\Gamma_{\elli}(\overline{H})$ is the one introduced in \S \ref{section formula for the multiplicity}. These relations can be seen as an analog of Theorem \ref{theo formula multiplicity} in the case where $E=F\times F$. Let $\chi$ be a continuous character of $H(F)$. In the particular case where $\pi=\St(H)$ and $\sigma=\St(H)\otimes \chi$, by \S \ref{eq 5.5.2} we get the relation

\begin{align}\label{eq 5.6.1}
\displaystyle \int_{\Gamma_{\elli}(\overline{H})}D^H(x)\chi(x)dx=\left\{
    \begin{array}{ll}
        1 & \mbox{if } \chi=\mathbf{1} \\
        0 & \mbox{otherwise}
    \end{array}
\right.
\end{align}

\noindent Indeed if $\chi$ is nontrivial then $\St(H)\otimes \chi\not\simeq \St(H)$ since $\St(H)$ and $\St(H)\otimes \chi$ have different cuspidal supports.

\subsection{Second application: multiplicity of the Steinberg representation}

\begin{theo}\label{theo Steinberg}
For every Galoisian character $\chi$ of $H(F)$ we have

$$\displaystyle m(\St(G),\chi)=\left\{
    \begin{array}{ll}
        1 & \mbox{if } \chi=\omega_{H,E} \\
        0 & \mbox{otherwise}
    \end{array}
\right.
$$
\end{theo}

\noindent\ul{Proof}: Let $\chi$ be a Galoisian character of $H(F)$. If the restriction of $\chi$ to the center of $H(F)$ is nontrivial then obviously $m(\St(G),\chi)=0$ since $\St(G)$ has trivial central character. We assume now that $\chi$ restricted to the center of $H(F)$ is trivial. Let $H_{ad}$ be the adjoint group of $H$ and set $G_{ad}:=R_{E/F}H_{ad,E}$. Let $H(F)_{ad}$ denote the image of $H(F)$ by the projection $H(F)\to H_{ad}(F)$ (i.e. the quotient $H(F)/Z_H(F)$). Then, since $\St(G)$ is the pullback of $\St(G_{ad})$ to $G(F)$, by Frobenius reciprocity we have

$$\displaystyle m(\St(G),\chi)=\dim \Hom_{H(F)_{ad}}(\St(G),\chi)=\dim \Hom_{H_{ad}(F)}(\St(G_{ad}),Ind_{H(F)_{ad}}^{H_{ad}(F)}\chi)$$

\noindent The representation $Ind_{H(F)_{ad}}^{H_{ad}(F)}\chi$ is a multiplicity-free sum of Galoisian characters containing $\omega_{H_{ad},E}$ if and only if $\chi=\omega_{H,E}$. This shows that the statement of the theorem for $H_{ad}$ implies the statement of the theorem for $H$. Thus, we may assume that $H$ is adjoint. Moreover, by \S \ref{eq 5.5.3} and Theorem \ref{theo inner forms}, up to replacing $H$ by its quasi-split inner form, we may also assume that $H$ is quasi-split. We will now use freely the notations introduced in \S \ref{section elliptic twisted Levi}. By Theorem \ref{theo formula multiplicity} and \S \ref{eq 5.5.1}, we have

$$\displaystyle m(\St(G),\chi)=\sum_{T\in \mathcal{T}_{\elli}(H)} \lvert W(H,T)\rvert^{-1}\sum_{P_0\subseteq P=MU} (-1)^{a_M-a_{M_0}}\sum_{\gamma\in \Gamma_M(T)} \int_{T(F)} D^M(\gamma^{-1} t\gamma)^{1/2}\chi(t)^{-1}dt$$

\noindent (Note that we have $\delta_P(\gamma^{-1}t\gamma)=1$ for all $\gamma$ and $t$ as in the expression above since $\gamma^{-1}t\gamma$ is a compact element). By Proposition \ref{prop twisted Levi}(\ref{prop twisted Levi ii}), it follows that

$$\displaystyle m(\St(G),\chi)=\sum_{(\mathcal{M},P)\in \mathcal{C}}(-1)^{a_{\mathcal{M}}-a_{M_0}} \int_{\Gamma_{\elli}(\mathcal{M})} D^{\mathcal{M}}(x) \chi(x)^{-1}dx$$

\noindent where for all $(\mathcal{M},P)\in \mathcal{C}$ we have set $a_{\mathcal{M}}:=a_M$ with $M:=R_{E/F}\mathcal{M}_E$. Thus, by \ref{eq 5.6.1} we also have

$$\displaystyle m(\St(G),\chi)=\sum_{(\mathcal{M},P)\in \mathcal{C}}(-1)^{a_{\mathcal{M}}-a_{M_0}} (\chi_{\mid \mathcal{M}},\mathbf{1})$$

\noindent where $\chi_{\mid \mathcal{M}}$ denotes the restriction of $\chi$ to $\mathcal{M}(F)$ and $(.,.)$ denotes the natural scalar product on the space of virtual characters of $\mathcal{M}_{ab}$. By Proposition \ref{prop twisted Levi}(\ref{prop twisted Levi iii}), (\ref{prop twisted Levi iv}) and (\ref{prop twisted Levi v}), this can be rewritten as

$$\displaystyle m(\St(G),\chi)=\sum_{I\subset \Delta_-}(-1)^{a_{\mathcal{M}_I}-a_{M_0}} \lvert \ker1^1(F;\mathcal{M}_I,H)\rvert (\chi_{\mid \mathcal{M}_I},\mathbf{1})$$

\noindent By Frobenius reciprocity, it follows that

$$\displaystyle m(\St(G),\chi)=\sum_{I\subset \Delta_-}(-1)^{a_{\mathcal{M}_I}-a_{M_0}} \lvert \ker1^1(F;\mathcal{M}_I,H)\rvert (\chi,Ind_{\mathcal{M}_{I,ab}}^{H_{ab}}\mathbf{1})$$

\noindent It is easy to see that $(-1)^{a_{\mathcal{M}_I}-a_{M_0}}=(-1)^{\lvert \Delta_- - I\rvert}$ for all $I\subset \Delta_-$ and therefore by Proposition \ref{prop twisted Levi}(\ref{prop twisted Levi vi}) the last expression above is equal to $1$ if $\chi=\omega_{H,E}$ and $0$ otherwise. $\blacksquare$

\addcontentsline{toc}{section}{Acknowledgment}
\section*{Acknowledgment}
I am grateful to Jean-Loup Waldspurger for a very careful proofreading of a first version of this paper. I also thank the referee for correcting many inaccuracies and for the numerous comments to make the text more readable. The author has benefited from a grant of Agence Nationale de la Recherche with reference ANR-13-BS01-0012 FERPLAY.

\end{document}